\documentclass[a4paper,fleqn]{cas-sc}

\usepackage[authoryear]{natbib}
\newcommand{\mathbft}[1]{\mathbf{{#1}}}
\usepackage[utf8]{inputenc} 
\usepackage[T1]{fontenc}    
\usepackage{hyperref}       
\usepackage{url}            
\usepackage{booktabs}       
\usepackage{amsfonts}       
\usepackage{nicefrac}       
\usepackage{microtype}      
\usepackage{lipsum}
\usepackage{enumitem}
\usepackage{float}
\usepackage{fancyhdr}       
\usepackage{graphicx}       
\graphicspath{{media/}}     
\usepackage{amsmath}
\usepackage{comment}
\definecolor{mycolor}{rgb}{0.0, 0.5, 0.0}

\usepackage{amsmath,amssymb, amsfonts}

\usepackage{physics}
\usepackage{algorithm}
\usepackage[noend]{algpseudocode}
\usepackage{ulem}

\def\tsc#1{\csdef{#1}{\textsc{\lowercase{#1}}\xspace}}
\tsc{WGM}
\tsc{QE}

\begin{document}
\let\WriteBookmarks\relax
\def\floatpagepagefraction{1}
\def\textpagefraction{.001}

\shorttitle{A spline-based hexahedral mesh generator for coronary arteries}    

\shortauthors{Fabio Marcinn\`{o} et Al}  

\title [mode = title]{A spline-based hexahedral mesh generator for patient-specific coronary arteries}  



%

\author[1]{Fabio Marcinn\`{o}}

\cormark[1]


\ead{fabio.marcinno@epfl.ch}


\credit{}

\affiliation[1]{organization={Institute of Mathematics, École Polytechnique Fédérale de Lausanne (EPFL)},
            city={Lausanne},
            country={Switzerland}}

\author[1]{Jochen Hinz}


\ead{jochen.hinz@epfl.ch}


\credit{}

\author[1,2]{Annalisa Buffa}


\ead{annalisa.buffa@epfl.ch}


\credit{}

\author[1]{Simone Deparis}


\ead{simone.deparis@epfl.ch}


\credit{}

\affiliation[2]{organization={Instituto di Matematica Applicata e Tecnologie Informatiche ’E. Magenes’ (CNR)},
            city={Pavia},
            country={Italy}}

\cortext[1]{Corresponding author}



\begin{abstract}
This paper presents a spline-based hexahedral mesh generator for tubular geometries commonly encountered in h{\ae}modynamics studies, in particular coronary arteries.
We focus on techniques for accurately meshing \textcolor{black}{vessels with stenoses and aneurysms}, as well as non-planar bifurcations.
Our approach incorporates several innovations, including a spline-based description of the vessel geometry in both the radial and the longitudinal directions, the use of Hermite curves for modeling non-planar bifurcations, \textcolor{black}{and a generalization to non-planar $n$ intersecting branches}.
This method eliminates the need for a concrete vessel surface, grid smoothing, and other post-processing. \textcolor{black}{A technique to generate grids with boundary layers is also presented.}
We validate the generated meshes using commonly employed quality indices, \textcolor{black}{compare them against state-of-the-art mesh generators and apply our method to complex coronary trees.} Finally, we present finite element fluid flow simulations with physiological boundary conditions. \textcolor{black}{To validate the proposed framework, a wall-shear-stress-based convergence test and computations of h{\ae}modynamic indices are also presented.}
\end{abstract}


\begin{highlights}
\item Structured grid generator based on a spline-based approach; 
\item Algorithms to generate non-planar bifurcations, non-planar n-furcations and boundary layer grids;
\item Quality mesh indices based analysis to show the efficiency of our method;
\item Fluid dynamics simulations and computation of h{\ae}modynamic indices.
\end{highlights}


\maketitle

\section{Introduction}
Cardiovascular diseases are the leading cause of death in the world \cite{kaptoge2019world}.
In particular, ischaemic heart diseases (IHDs), also called coronary artery diseases (CADs), account for approximately 40$\%$ of the fatalities in European countries \cite{timmis2022european}.
Although in daily clinical practice numerical methods are still not widely adopted, in the last 25 years, computational methods have been established as efficient tools for analyzing the h{\ae}modynamics within arteries \cite{urquiza2006multidimensional, vignon2006outflow, taylor2013computational,quarteroni2017cardiovascular, schwarz2023beyond}. State-of-the-art frameworks for the numerical simulation of blood dynamics (in particular, for coronary arteries) involve the following steps:

\begin{itemize}
    \item Pre-processing: Clinical image segmentation to reconstruct the arterial geometry;
    \item Mesh generation: Generation of an unstructured grid using tetrahedal and/or hexahedral elements;
	\item Numerical simulation: Computational Fluid Dynamics (CFD) simulation using a Finite Element Method (FEM) solver.
\end{itemize}

As discussed in \cite{schwarz2023beyond}, unfortunately, this pipeline is associated with several bottlenecks.
\begin{enumerate}[label=\Roman*., wide=0pt]
    \item  In the pre-processing step, the manual reconstruction of the arterial geometry from the clinical data is time consuming and prone to user errors. 
    \item The use of an unstructured approach in the mesh generation step impedes a precise reconstruction of the vessel geometry and necessitates a high-quality triangulation of the vessel surface \cite{remacle2010high} to warrant acceptable volumetric mesh quality. In addition, the grid may be encoded with numerous parameters. 
    \item While offering a high degree of numerical precision, FEM solvers may require considerable computational resources, in particular for simulating coronary arteries, wherein the presence of cardiac motion and microvasculature at the artery's end points create challenges in obtaining physiological results \cite{schwarz2023beyond}.
\end{enumerate}
In this work, we address the problem regarding the mesh generation presenting a fast and robust approach for the generation of structured mesh exploiting the use of splines.

\subsection{State of the art of the mesh generator}

In numerical h{\ae}modynamics, various mesh generators have been proposed, ranging from completely unstructured approaches to hybrid meshes (combining tetrahedra with unstructured and/or structured hexahedra) to fully structured hexahedral approaches.

The volumetric unstructured mesh is generated using a Delaunay triangulation, starting from a triangulation of the vessel surface. Commonly-employed unstructured mesh generators in h{\ae}modynamics, based on the Delaunay algorithm, are \textit{Gmsh} \cite{geuzaine2009gmsh} and \textit{SimVascular} \cite{updegrove2017simvascular}. 
It should be noted that the standard software for generating unstructured tetrahedal meshes is the \textit{Vascular Model ToolKit} (VMTK) \cite{izzo2018vascular}. 
Despite the aforementioned challenges associated with unstructured tetrahedral/hexahedral approaches, such as the inexact representation of the vessel surface, their wide adoption can be attributed to the relative ease of meshing complex features such as aneurysms, stenoses, and bifurcations. 

Hybrid meshes leverage the advantages of both unstructured and structured techniques. Hybrid approaches typically generate an unstructured grid in the interior of the vessel, while close to the vessel surface a finer structured grid is introduced to capture the boundary layer \cite{sahni2008adaptive}.
In \cite{marchandise2013cardiovascular}, a fully automatic procedure to generate a grid inside tubular geometries based on both unstructured and structured cells is proposed which is compatible with isotropic tetrahedral meshes, anisotropic tetrahedral meshes, and mixed hexahedral/tetrahedral meshes.
Also approaches that exchange the roles of the mesh's structured and unstructured regions have been proposed. For example, \cite{pinto2018numerical} proposes an algorithm that introduces hexahedral cells close to the center of the vessel while a finer grid of tetrahedral cells close to the artery wall is introduced. 
\textcolor{black}{It is worth noting that the} aforementioned fully unstructured approaches are also capable of generating a structured mesh close to the vessel surface to better capture the boundary layer.

To achieve more efficient simulations, hexahedral structured cells are a suitable choice. The concept of using a structured approach to generate meshes is not novel in numerical h{\ae}modynamics. For instance, it has already been numerically demonstrated that structured flow-oriented cells are more efficient than unstructured cells in h{\ae}modynamic simulations. For example, to achieve similar numerical accuracy, unstructured grids necessitate the introduction of more cells compared to structured approaches \cite{de2010patient}, while additionally introducing numerical diffusivity \cite{ghaffari2017large}.

While more challenging than their unstructured counterparts, several strategies have been proposed for structured mesh generation. 
Leveraging the template grid sweeping approach, in \cite{zhang2007patient}, the authors propose a procedure to build hexahedral solid NURBS (Non-Uniform Rational B-Splines) meshes for patient-specific vascular geometries. While still requiring a reconstruction of the vessel surface, a generalization to tri- and n-furcations is discussed which divides the geometry into a master and slave branches. Isogeometric Analysis (IGA) simulations are carried out on a select number of geometries.
In \cite{de2010patient}, the authors propose a semi-automatic procedure to reconstruct a structured grid from coronary angiographies. Using steady simulation, the authors demonstrate that a hexahedral mesh requires fewer cells and less CPU time to achieve the same numerical accuracy as an unstructured grid.
In \cite{de2011patient}, the same group develops another semi-automatic procedure to generate a structured hexahedral mesh starting from the vessel surface. The mesh generator is tested on a patient-specific carotid bifurcation.
Both in \cite{de2010patient,de2011patient}, the bifurcations are modeled as planar entities.
Based on a template grid sweeping method, in \cite{ghaffari2017large}, a procedure to generate structured meshes without the need for vessel surfaces is proposed. The method is applied to six patient-specific cerebral arterial trees. Additionally, the same study demonstrates the convenience of using flow-oriented cells by comparing the performance of structured and unstructured approaches.
\textcolor{black}{In \cite{urick2019review}, the authors propose a novel CAD-integrated pipeline (both image reconstruction and meshing) for performing finite element analysis. 
The use of CAD allows for handling arbitrary vascular topology (both convex and non convex sections), and accurately segment the clinical images.
However, the method still requires the knowledge of the triangulated surface mesh, and as the authors mention the generalization to non planar $n$-furcations is still not possible with their method.
}
We refer to \cite{bovsnjak2023higher}, where a block hexahedral structured approach, specifically designed for tubular geometries, is presented. The mesh generator quality is gauged using a test case involving a patient-specific human aorta and coronary artery. \textcolor{black}{However, in this work,} we note that the grids exhibit $C^0$ continuity along the centerline direction while the representation of bifurcations appears unnatural. Moreover, the authors report that this approach is ineffective in handling pathologies such as stenoses or aneurysms.
To the best of our knowledge, the only study concerning a hexahedral mesh generator applied to a large dataset of patient-specific cerebral arteries has been conducted in \cite{decroocq2023modeling}. 
A pure template grid sweeping method to generate a structured grid from no more than a centerline and the associated radius information is proposed. However, the method still requires the computation of vessel surfaces to generate the mesh of the bifurcation; an extension to planar n-furcations is proposed. While arguably the current state-of-the-art for structured h{\ae}modynamics grid generation, full automation remains a major challenge since, \textcolor{black}{as the authors state}, almost 20\% of the dataset's bifurcations required manual intervention \cite{decroocq2023modeling}.
\begin{table}[h!]
\caption{Differences with respect to the state of the art: \textit{Bifurcation} means whether it is modeled as a planar or non-planar entity. The same holds for \textit{n-furcations}; additionally, "no" indicates that the bifurcation can not be generalized to n-furcations, and master branch indicates that the procedure requires the presence of a master vessel. \textit{Surface} indicates whether the meshing algorithm requires information of the vessel surface as input or in any other step of the mesh generation. \textit{Complex geometries} indicate whether the algorithm can capture the presence of pathological cases (aneurysm or stenosis). Finally, \textit{spline} indicates if the numerical grid is generated using a fully spline-based approach, leaving room for the possibility of performing  Isogeometric Analysis (IGA) in the future. \textit{Large dataset} indicates whether the mesh generator has been tested on a large number of geometries.}
\label{tab: state_of_the_art}
\centering
 \begin{tabular}{||c c c c c c c||} 
 \hline
 Work & \textit{Bifurcation} & \textit{n-furcations}  & \textit{Surface} & \textit{Complex geometries} & \textit{Spline} & \textit{Large dataset}\\ [0.5ex] 
 \hline\hline
 \cite{zhang2007patient}     & non-planar & master branch  & yes 		 & yes            & yes & no\\
 \cite{de2010patient}        & planar     & no                      & yes          & yes                 & no& no\\ 
 \cite{de2011patient}        & planar 	 & no                      & yes 		 & yes				 & no& no\\ 
 \cite{ghaffari2017large}	& non-planar & no                      & no 			 & yes 				 & no& no\\  
 \cite{urick2019review} & non-planar & planar                      & yes 			 & yes 				 & yes& no\\  
 \cite{bovsnjak2023higher} 	& non-planar & no                      & no			 & no 				 & yes& no\\
 \cite{decroocq2023modeling} & non-planar & planar                  & yes 		 & yes                & no& yes \\
 Ours 					& non-planar & non-planar              & no 			 & yes 				 & yes& no\\[1ex] 
 \hline
 \end{tabular}
\end{table}

\subsection{\textcolor{black}{Novelties and datasets}}

\textcolor{black}{To overcome the challenges associated with the vessel surface reconstruction, this paper proposes a structured grid generator that operates entirely on a spline-based approach, without requiring a concrete vessel surface as input, neither for the single branch nor for the bifurcation. The splines allows for a lighter representation of the mesh and for a computationally inexpensive sampling of a finer mesh from the spline mapping operator.
More precisely, the generation of the single branch mesh is based on a template sweeping method which builds a grid of the vessel section by approximating a harmonic map over the underlying geometry model's spline basis, instead of radially expanding the mesh vertices as proposed in \cite{de2010patient, de2011patient, ghaffari2017large, decroocq2023modeling}.
Meanwhile, the non-planar bifurcation mesh is based on a reconstruction of the vessel surface using Hermite curves which generalizes straightforwardly to non-planar $n$-furcations while additionally avoiding the ad-hoc introduction of a master branch \cite{zhang2007patient}). At the best of our knowledge, our approach is new and it has not been previously introduced in the literature. In Tab. \ref{tab: state_of_the_art}, we summarize the main differences with respect the other works.}

\textcolor{black}{
There are two main datasets used in this work.
For all the test cases, except the one presented in Sec.\ \ref{sec: complex_tree}, patient-specific geometries are selected from invasive coronary angiographies performed in the Fractional Flow Reserve versus Angiography for Multivessel Evaluation 2 (FAME 2) trial. 
Details of the trial have been published previously \cite{de2012fractional}. The trial conformed to the ethical guidelines of the 1975 Declaration of Helsinki as reflected in a priori approval by the participant institution's human research committee.
For this dataset, the clinical images have been firstly segmented in \cite{MAHENDIRAN2024132598}, and then the mesh generator's input data (vessel centerline and radius parameterization) is reconstructed following an algorithm presented in \cite{ccimen2016reconstruction}.
The test case presented in Sec.\ \ref{sec: complex_tree}, centerline and radius information are reconstructed from both the left and right coronary tree of a patient of the Vascular Model Repository. Details on this dataset can be found in \cite{wilson2013vascular}.}

\subsection{Paper's outline}

\noindent The work is organized as follows. 
\textcolor{black}{Sec.\ \ref{sec: reconstruction} summarizes the procedure to reconstruct the geometrical information from the clinical images.}
Sec.\ \ref{sec: method} discusses the method for the generation of the numerical grid, namely: the single branch mesh (Sec.\ \ref{sec: single_branch}), including the circular vessel section (Sec.\ \ref{sec: plane_mesh}) and the case of general convex section (Sec.\ \ref{sec: general_convex}). After that, a strategy to reconstruct 
the non-planar bifurcation (Sec.\ \ref{sec: bifurcation}) and non-planar 
$n$-furcation grids (Sec.\ \ref{sec: n_furcations}) are discussed. Finally, a technique to generate grids with boundary layers is presented (Sec.\ \ref{sec: BL}). The section about numerical results (Sec.\ \ref{sec: results}) is divided into two main parts: Sec.\ \ref{sec: meshing_results} discusses the results of the mesh generator and Sec.\ \ref{sec: fluid_results} presents a numerical simulation on a patient-specific coronary artery geometry. Conclusions are reported in Sec.\ \ref{sec: conclusion}.

\section{\textcolor{black}{Geometrical information reconstruction}}\label{sec: reconstruction}


\textcolor{black}{
Nowadays, the gold standard for assessing the severity of coronary artery disease is the anatomical evaluation of stenosis severity during invasive coronary angiography (ICA). In cases of moderate stenoses, this anatomical assessment can be complemented by a haemodynamic evaluation known as fractional flow reserve (FFR) \cite{achenbach2017performing, pijls2012functional,tu2016diagnostic}.
Indeed, given that coronary arteries are vessels with a low to moderate diameter, the ICA by the use of medium contrast provides the most accurate and precise visualization of coronary artery anatomy with respect to $3$D coronary computer tomography angiography (CCTA), which is limited by a low temporal and spatial resolution \cite{stefanini2015can}.}

\begin{figure}[t]
\centering
\includegraphics[height = 14cm, keepaspectratio]{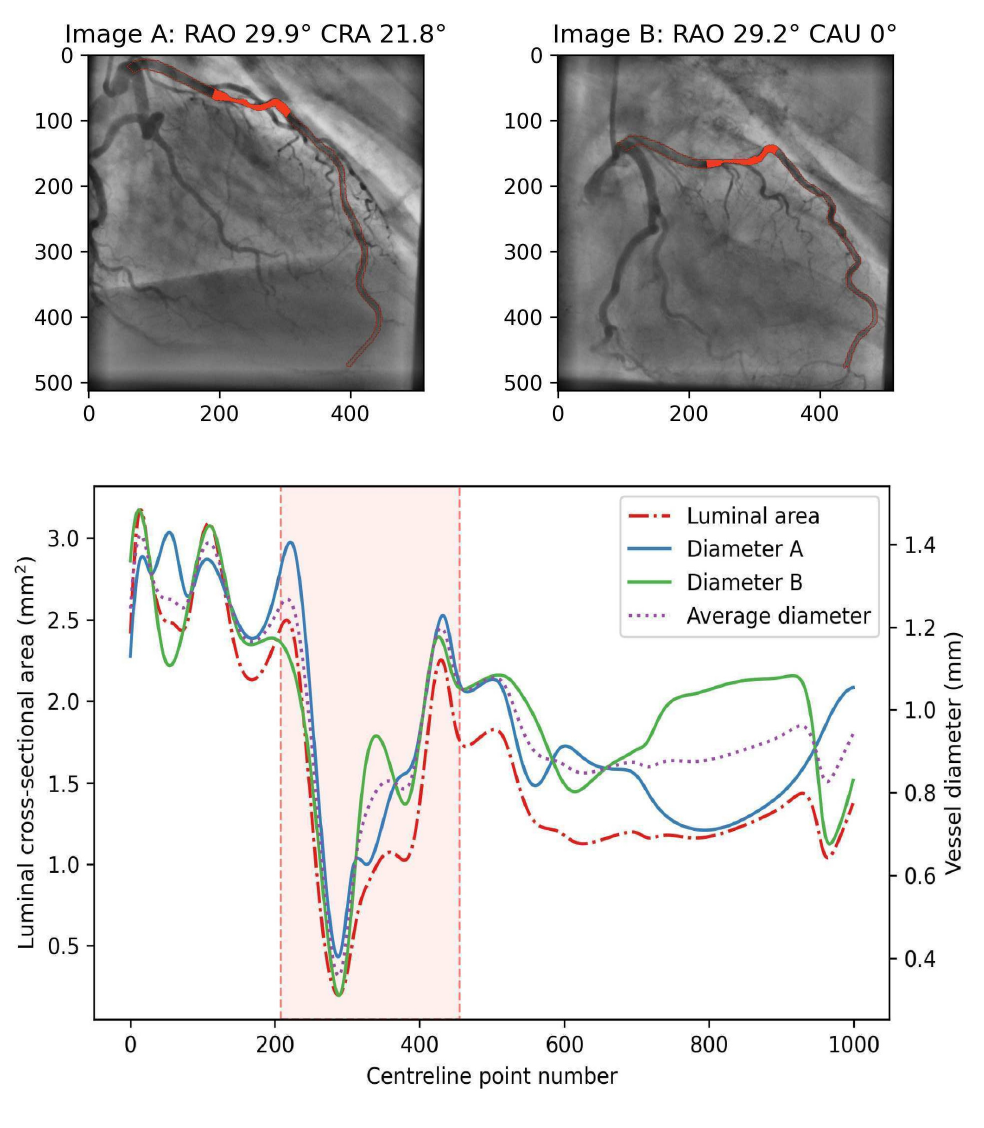}
\caption{\textcolor{black}{Top: Two different projections taken during the ICA are shown. RAO and CRA refer to the right anterior oblique and cranial views, respectively, while CAU refers to the caudal view. The section of the artery with stenosis is highlighted in red. Bottom: The geometrical output is depicted. Cubic splines representing the vessel diameters: diameter A is computed
from the first view (Image A) and diameter B is computed from the second view (Image B).
The average diameter is the mean between diameter A and diameter B, and it is used to compute the luminal cross-sectional area of the vessel. While the centerline information is not explicitly shown, it is also parameterized as a cubic spline. The image is courtesy of Edward Andò (EPFL) and Thabo Mahendiran (CHUV).}}
\label{Fig: imaging_algorithm}
\end{figure}

\textcolor{black}{
The reconstruction algorithm, is fed with two snapshots of almost orthogonal ICA projections (see Fig. \ref{Fig: imaging_algorithm}, top).
From these two projections, a fast and robust deep learning based algorithm (the one presented in \cite{MAHENDIRAN2024132598}) segments the selected vessel tracts in each frame (the one with red contours Fig. \ref{Fig: imaging_algorithm}, top), computing the centerline and vessel diameter for each projection (see Fig. \ref{Fig: imaging_algorithm}, bottom).
The epipolar lines of these two ICA projections are then intersected to generate $3$D points that represent the centerline of the vessel, while the final diameter is computed as the average of the diameters derived from the projections.
In particular, the output of this procedure, as depicted in Fig.\ \ref{Fig: imaging_algorithm}, bottom, is represented by a cubic spline description of both the diameters computed from the two ICA projections, and the luminal cross-sectional area computed using the average of the diameters. In Fig.\ \ref{Fig: imaging_algorithm}, the centerline information is reported under "centerline point number" evaluated in 1000 samples, but we remark that also the centerline of the vessel is output as a cubic spline.}

\textcolor{black}{
Parameterizations based on (cubic) splines constitute a natural choice for reconstructing the vessel centerline thanks to the option of achieving $C^{\geq 1}$ global continuity.
Moreover, only valid parametrization are considered, meaning that if the centerline contains unnatural kinks resulting in an invalid centerline-radius parametrization, i.e.\ a self-intersecting geometry, a penalization to the centerline’s spline curvature can be applied to recover a valid parametrization. }

\section{The mesh generator} \label{sec: method}
This section discusses this paper's meshing strategy.
The inputs to our procedure are a spline-based representation of the centerline curve and a convex parametrization of the contour of the vessel section.
Such an input is provided by the reconstruction algorithm 
summarized in Sec.\ \ref{sec: reconstruction}, \textcolor{black}{or by manual reconstructions.} 
Of course, if the data is provided in the form of discrete data points \textcolor{black}{as in the latter case}, it may always be converted to an appropriate spline-based input using standard spline fitting routines \cite{virtanen2020scipy}.
In the following, for the sake of the explanation, we first describe the methodology for circular vessel sections (see Sec.\ \ref{sec: plane_mesh}), which is then generalized to general convex cross sections (see Sec. \ref{sec: general_convex})

\subsection{The single branch} \label{sec: single_branch}
\begin{figure}[t]
\centering
\includegraphics[width =7cm, keepaspectratio]{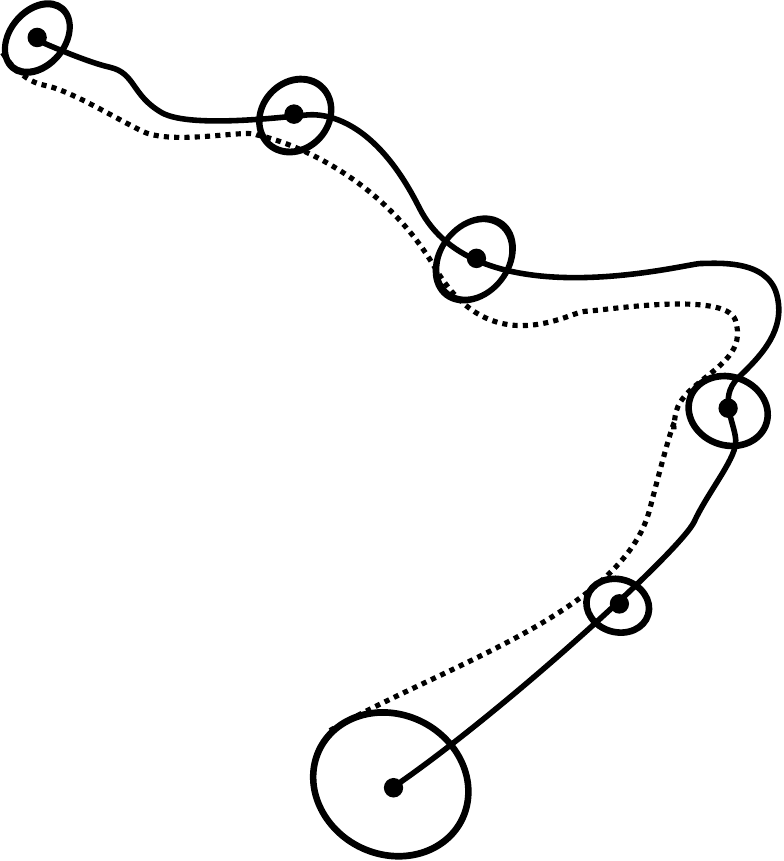}
\includegraphics[width =8cm, keepaspectratio]{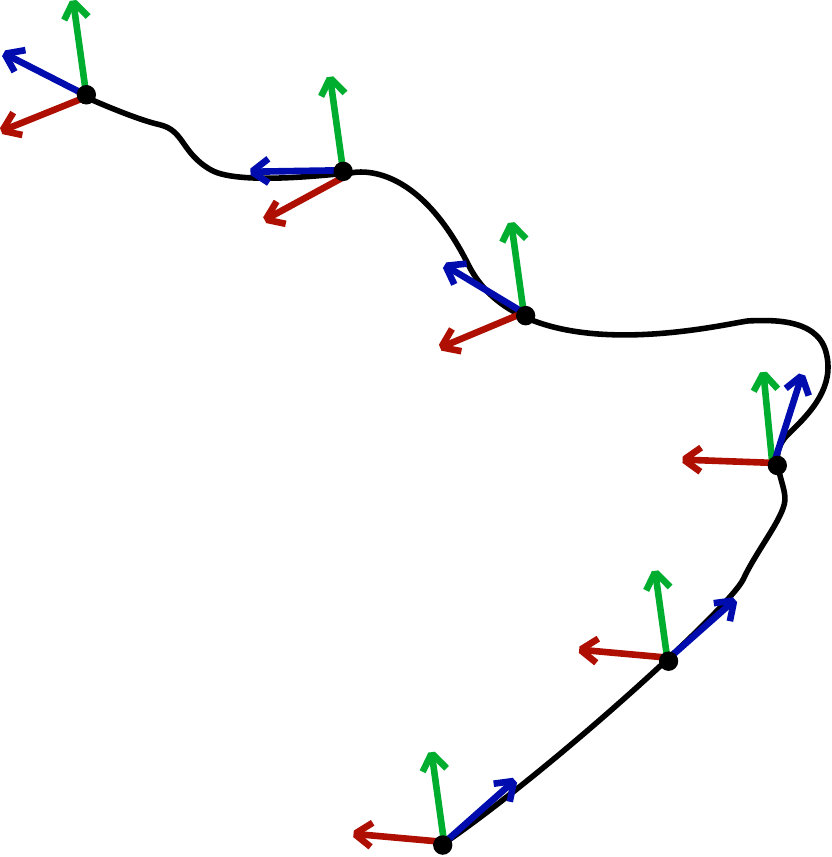}
\includegraphics[width =7cm, keepaspectratio]{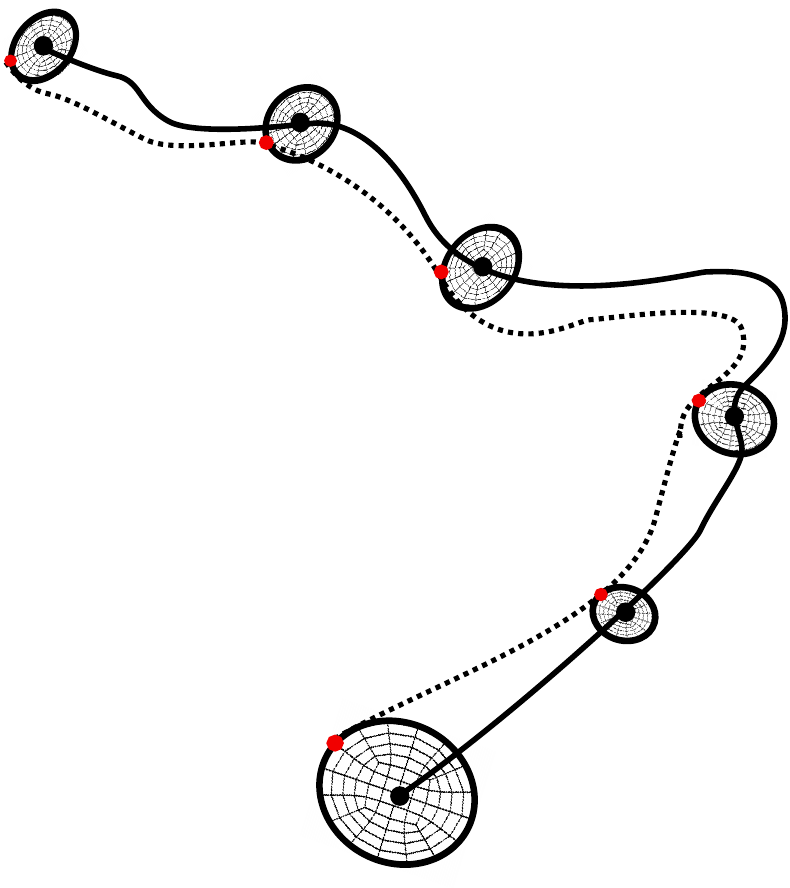}
\includegraphics[width =8cm, keepaspectratio]{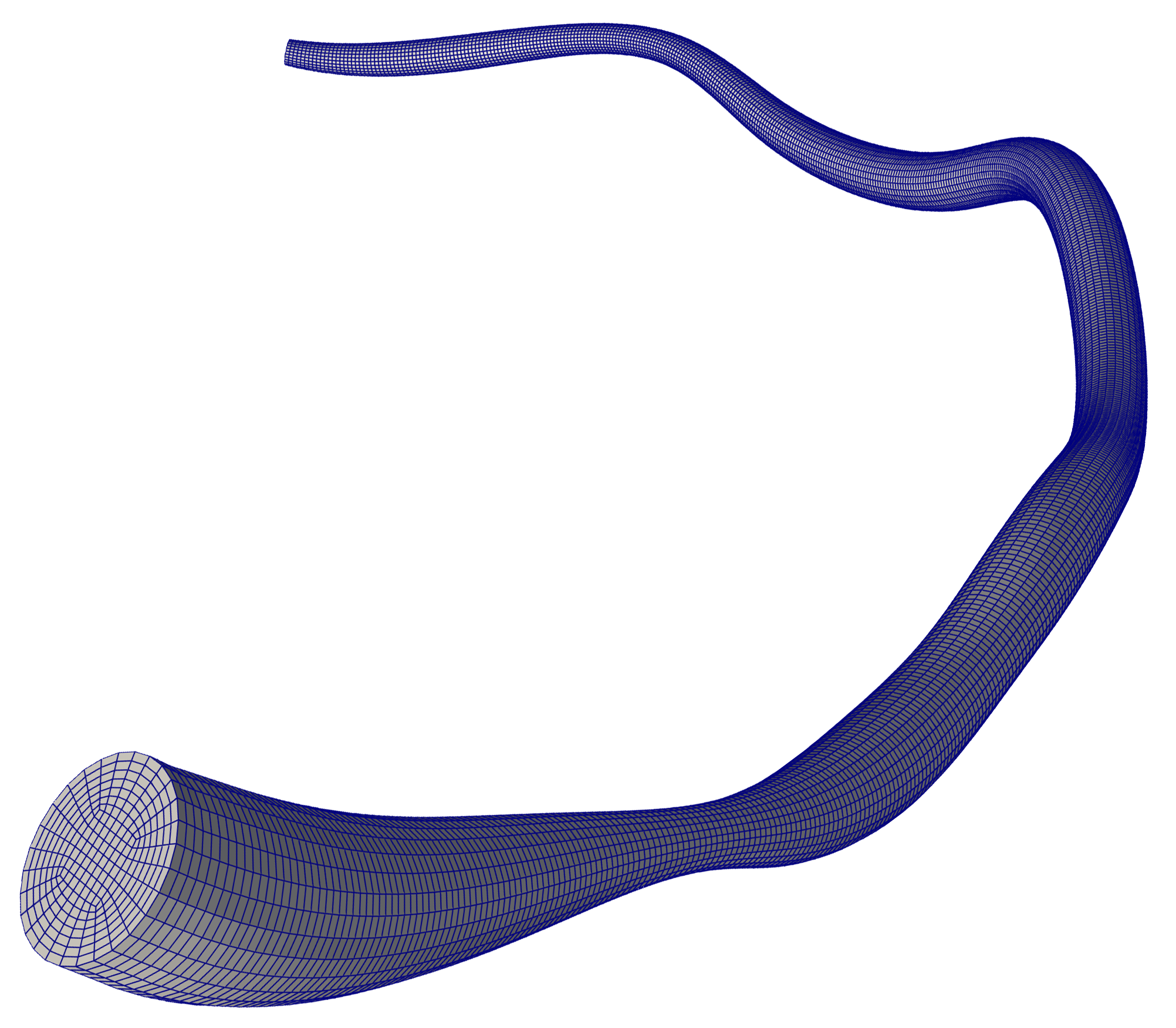}
\caption{Top-Left: Parametrization of the vessel centerline (continuous line) and radius (dashed lines). The contours with different radius are depicted in six locations.
Top-Right: The no-twist frames. The tangent, normal and binormal vectors are colored in blue, green and red, respectively.
Left-Bottom: Insertion of the grid discs along the vessel centerline following the no-twist frame. The red point shows the utility of the no-twist frame.
Right-Bottom: Generation of the volumetric mesh.}
\label{Fig: spline_parametrization}
\end{figure}

As depicted in Fig.\ \ref{Fig: spline_parametrization}, top-left, we are given a one-dimensional centerline in $\mathbb{R}^3$ as a cubic spline curve:

\begin{equation*}
\mathbf{s}(t) = \sum_i \mathbf{b}_i N_{i, 3}(t) \,,
\end{equation*}
where the $\mathbf{b}_i \in \mathbb{R}^3$ are vectorial control points, while the $\{N_{i, 3}\}_i$ span a cubic spline space with global $C^2([0, 1])$ continuity resulting from an open knot vector. 
At each $t_0 \in [0, 1]$ the radius of the vessel cross section centered at $\mathbf{s}(t_0)$ is given by $r(t_0)$, where $r: [0, 1] \rightarrow \mathbb{R}^+$ is a scalar function from the linear span of the same cubic spline space. 

At this point, we compute the unit tangent of the spline as follows:
\begin{equation*}
\mathbf{t}(t) = \frac{\mathbf{s}'(t)}{\left \| \mathbf{s}'(t) \right \|_2} \,.
\end{equation*}

With the aim of creating a volumetric mesh by inserting vessel sections (see Sec.\ \ref{sec: plane_mesh}) along the centerline, we compute a rotation minimizing frame (RMF) to minimize the torsion between two consecutive vessel sections.
Given a curve in $3$D space, a RMF can be seen as an intrinsic property of $\mathbf{s}(t)$ and it can be computed using two approaches \cite{wang2008computation}: $a)$ a discrete approach that computes the frame only at precisely the section' centroid coordinates or $b)$ a continuous approach wherein the frame along the entire curve is implicitly given as the solution of an ODE. Here, $a)$ can be regarded as a numerical approximation of $b)$ over a discrete set of data points. While sacrificing some accuracy, $a)$ should be favored since it provides clear practical advantages such as drastically improved algorithmic simplicity and computational efficiency. We have furthermore found it to be sufficiently accurate in all our test cases. 
By default, the local $z$-direction is aligned with the tangent $\mathbf{t}(t_i)$, while the local $x$- and $y$-directions are aligned with the normal $\mathbf{n}(t_i)$ and binormal $\mathbf{m}(t_i)$, respectively, which are precomputed at a select number of coordinates $\mathbf{s}(t_i)$ over a given set of abscissae $\{t_i\}_i \subset [0, 1]$ using the algorithm presented in \cite{ebrahimi2021low}. 

In Fig.\ \ref{Fig: spline_parametrization}, top-right, the no-twist frames are depicted in select points along the centerline. 
In the following section, we discuss the generation of the grid inside the section of the vessel.

\subsubsection{The circular vessel section} \label{sec: plane_mesh}
\begin{figure}[t]
\centering
\includegraphics[width =15cm, keepaspectratio]{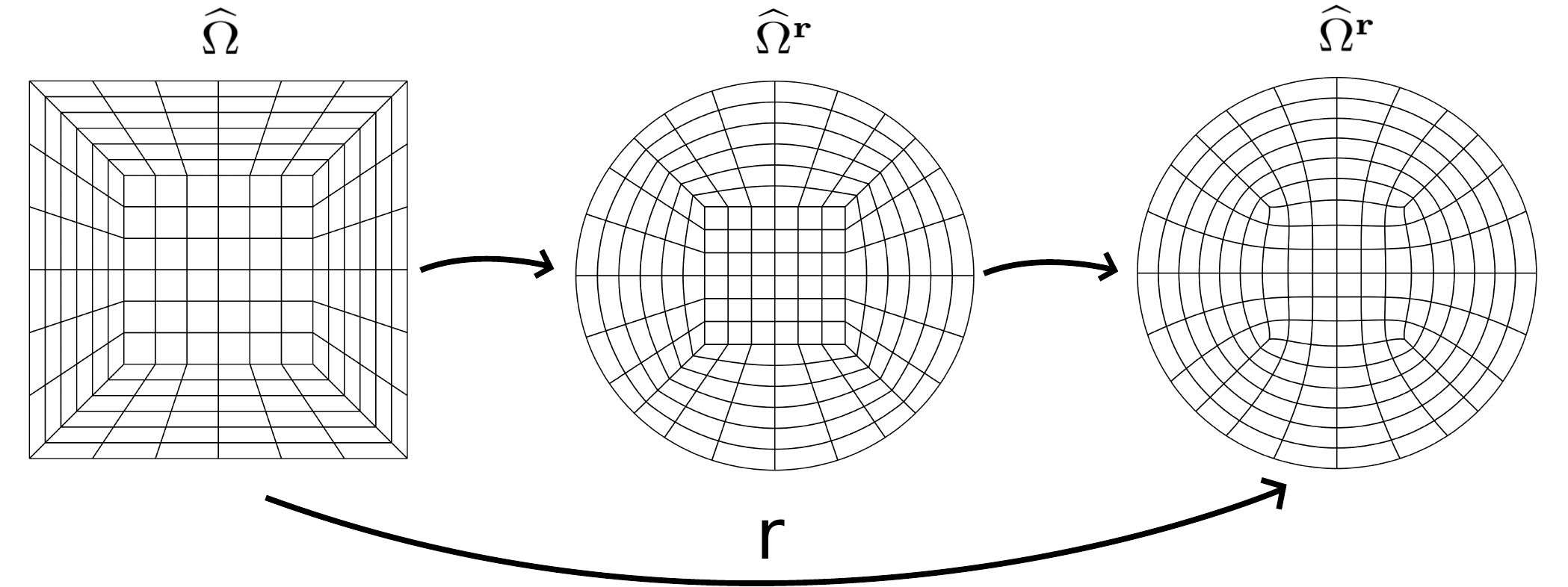}
\caption{Map from the square grid into disc grid.}
\label{Fig: square_to_smooth}
\end{figure}

\begin{figure}[t]
\centering
\includegraphics[width =15cm, keepaspectratio]{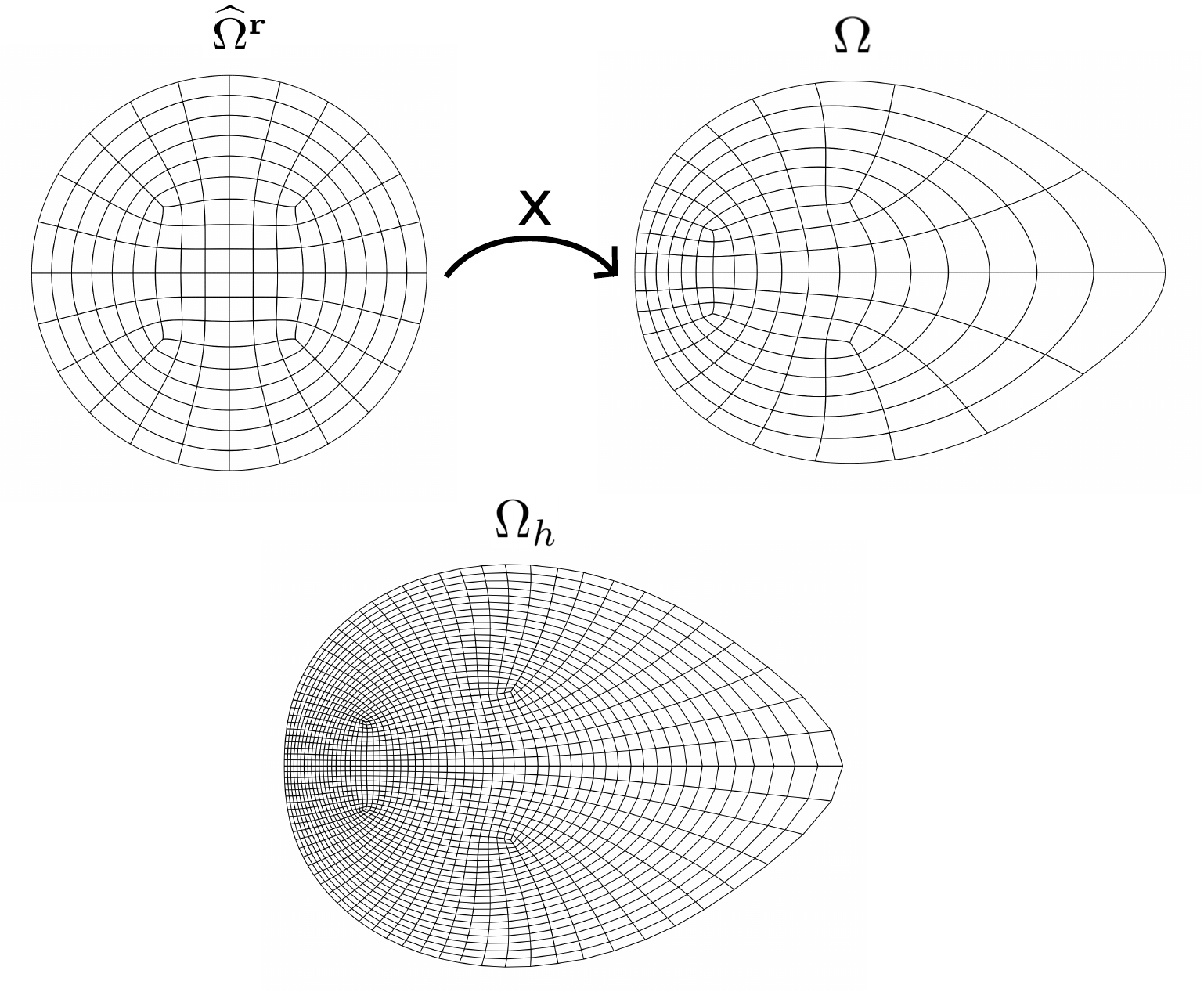}
\caption{Top: Effect of the harmonic map acting from the unit disc into a general convex section.
Bottom: A classical piece-wise linear mesh sampled from the surface spline.}
\label{Fig: harmonic_map}
\end{figure}

The generation of the $2$D circular vessel section is based on the radial expansion of a unit radius master section $\widehat{\Omega}^r$. The treatment of general convex cross-sectional shapes is discussed in Section \ref{sec: general_convex}. 

We start from a bilinear $5$-patch covering of a polygonal domain $\widehat{\Omega}$ with four edges (see Fig. \ref{Fig: square_to_smooth}, left). A finite-dimensional spline basis of polynomial order $p \geq 2$ is defined on the polygon's $5$-patch covering. The resulting finite-dimensional space satisfies $\mathcal{V}_h \subset H^1(\widehat{\Omega})$ while functions $\phi \in \mathcal{V}_h$ are assumed to be $C^{\geq 1}$ continuous on each individual patch. 
We map the $5$-patch topology defined in $\widehat{\Omega}$ onto an approximation of the unit disc using the manually constructed template map $\mathbf{r}^\prime: \widehat{\Omega} \rightarrow \widehat{\Omega}^r$, with $\mathbf{r}^\prime \in \mathcal{V}_h \times \mathcal{V}_h$, depicted in Fig. \ref{Fig: square_to_smooth}, center. For this, we apply the bilinearly-blended Coons' patch approach \cite{farin1999discrete} to each covering patch of $\widehat{\Omega}$ individually. 
While the practically motivated choice to require $\mathbf{r}^\prime \in \mathcal{V}_h \times \mathcal{V}_h$ causes $\widehat{\Omega}^r$ to no longer exactly equal the unit disc, the geometrical discrepancy is negligible and will be disregarded in the discussion that follows. 
The $5$-patch covered unit disc $\widehat{\Omega}^r$ is mapped onto itself via the use of a map $\mathbf{s}: \widehat{\Omega}^r \rightarrow \widehat{\Omega}^r$ that is sought as the discrete (over the same spline space) solution of an elliptic PDE problem that has a beneficial effect on the new map's grid quality (such as the suppression of steep interface angles between neighboring patches, see Fig. \ref{Fig: square_to_smooth}). For details, see \cite{hinz2024use}. 

We refer to the composition of the two maps as $\mathbf{r}: \widehat{\Omega} \rightarrow \widehat{\Omega}^r$. A parametrization of a disc with radius \textcolor{black}{$R > 0$} now simply follows from multiplying each control point of $\mathbf{r}: \widehat{\Omega} \rightarrow \widehat{\Omega}^r$ by \textcolor{black}{$R > 0$}. 

Given the strictly monotone sequence $\{ 0=t_0, t_1, \ldots, t_N=1 \}$ that corresponds to the displacements along the centerline curve, with above procedure, we generate a catalogue of planar geometries $\Omega_i$ with radii \textcolor{black}{$R(t_i) > 0$} and corresponding parametrization $\mathbf{x}_i^{2D}: \widehat{\Omega}^r \rightarrow \Omega_i$. In view of adopting a template grid sweeping approach, we define the lifted (into $\mathbb{R}^3$) cross sectional maps $\mathbf{x}_i^{3D}: \widehat{\Omega} \rightarrow \mathbb{R}^3$ that place the planar cross sectional parametrization along the centerline curve in the orientation dictated by the RMF. The relation is given by:
\begin{equation*}
\mathbf{x}^{3D}_{i}  = 
	\left[
	\begin{matrix}
	\mathbft{n}(t_i) & \mathbft{m}(t_i) 
	\end{matrix}
	\right] \mathbf{x}^{2D}_{i} + \mathbf{s}(t_i)
\end{equation*}
where $\mathbft{n}_i := \mathbf{n}(t_i)$ and $\mathbft{m}_i := \mathbf{m}(t_i)$ denote the normal and binormal vectors at $t = t_i$. 

As a next step, we sample dense meshes from the $\mathbf{x}_i^{3D}$. A volumetric mesh follows from connecting consecutive planar meshes via linear edges. 
Six cross-sectional disc grids are depicted in Fig.\ \ref{Fig: spline_parametrization}, bottom-left, where the red points show the orientation of the RMF. A volumetric mesh that results from this procedure is shown in Fig.\ \ref{Fig: spline_parametrization}, bottom-right.

\textcolor{black}{It is worth noting that, as show in Fig.\ \ref{Fig: spline_parametrization}, bottom-left, the use of a no-twist frame is crucial to reduce torsion along the centerline direction. 
In the figure, we highlight in red the points of the $2$D spline maps that share the same index in the $2$D connectivity of the map. The volumetric mesh is generated by connecting these points, i.e.\ , those with the same index in the $2$D connectivity but belonging to different sections.
To avoid mesh torsion, it is crucial that these points are close to each other.
This is ensured using a RMF. If a Serra-Frenet frame is used instead, the normal vector aligns with the direction of curvature and a  sudden change in curvature may lead to significant rotation between consecutive slices, potentially resulting in invalid grids.}

\subsubsection{General convex cross sections} \label{sec: general_convex}
The generation of general $2$D vessel section is performed through an approach based on harmonic maps capable of handling all convex cross-sectional shapes.
The harmonic map equation's exact solution obeys the so-called Rad\'o-Kneser-Choquet theorem which states a harmonic planar map $\mathbf{b}: \widehat{\Omega} \rightarrow \Omega$ is diffeomorphic in $\Omega$, provided $\Omega$ is convex  \cite{choquet1945type}. 
As a consequence, also a fine enough numerical approximation is folding free. Since each planar section is folding-free, a sufficiently tight stacking of the vessel sections results in a folding-free volumetric map (provided the input data satisfies reasonable regularity requirements). 

Given a boundary correspondence $\mathbf{F}: \partial \widehat{\Omega}^r \rightarrow \partial \Omega$, where $\Omega \subset \mathbb{R}^2$ denotes the convex target domain, the harmonic map problem for $\mathbf{x}: \widehat{\Omega}^r \rightarrow \Omega$ reads:
\begin{equation}
	i \in \{1, 2\} \quad \Delta \mathbf{x}_i = 0 \quad \text{in} \quad \widehat{\Omega}^{\mathbf{r}}, \quad \text{s.t.} \quad \mathbf{x} \vert_{\partial \widehat{\Omega}^{\mathbf{r}}} = P(\mathbf{F}),
    \label{eq: projection}
\end{equation}
where we denote the push-forward of $\mathcal{V}_h$ to $\widehat{\Omega}^r$ via $\mathbf{r}: \widehat{\Omega} \rightarrow \widehat{\Omega}^r$ by $\mathcal{V}_h^r$. Here, P($\cdot$) denotes a suitable least squares approximation operator that fits the boundary data $\mathbf{F}: \partial \widehat{\Omega}^\mathbf{r} \rightarrow \partial \Omega$ (which may alternatively come in the form of scattered point data) to the spline DOFs of $\mathcal{V}_h^{\mathbf{r}}$ that are nonvanishing on $\partial \widehat{\Omega}^\mathbf{r}$. 
The advantage of a spline-based approach is further underlined in this context: an accurate spline-based representation of the correspondence $\mathbf{F}: \partial \widehat{\Omega}^r \rightarrow \partial \Omega$ is achieved with drastically reduced dimensionality of $\mathcal{V}_h$ over an analogous piecewise-linear approach for virtually all clinically relevant input contours. 

The discretization over $\mathcal{V}_h$ results in a decoupled linear problem with one linear system for both entries of $\mathbf{x}: \widehat{\Omega}^{\mathbf{r}} \rightarrow \Omega$. In fact, only the right-hand sides of the equations differ. The problem can be cast into the form:
\begin{equation}
\label{eq:harmonic_map_discretized}
	i \in \{1, 2\} \quad A \mathbf{c}^i = \mathbf{f}^i,
\end{equation}
where $A$ denotes the stiffness matrix over the canonical basis of $\mathcal{V}_h^{\mathbf{r}}$, while $\mathbf{f}^i$ follows from the Dirichlet data. The $\mathbf{c}^i$ then represent the weights of the $\mathbf{x}_i \in \{\mathbf{x}_1, \mathbf{x}_2 \}$ with respect to the canonical basis of $\mathcal{V}_h^{\mathbf{r}} \cap H^1_0(\widehat{\Omega}^{\mathbf{r}})$. All remaining weights follow from the Dirichlet data. We compute the Cholesky decomposition of the linear SPD operator $A$ which is then employed in all planar parameterisation problems (usually in the hundreds). The only ingredients that need to be reassembled are the right-hand-sides $\mathbf{f}^i$, which require no more than a number of sparse boundary integrals and are therefore computationally inexpensive.

\begin{figure}[t]
\centering
\includegraphics[width =7cm, keepaspectratio]{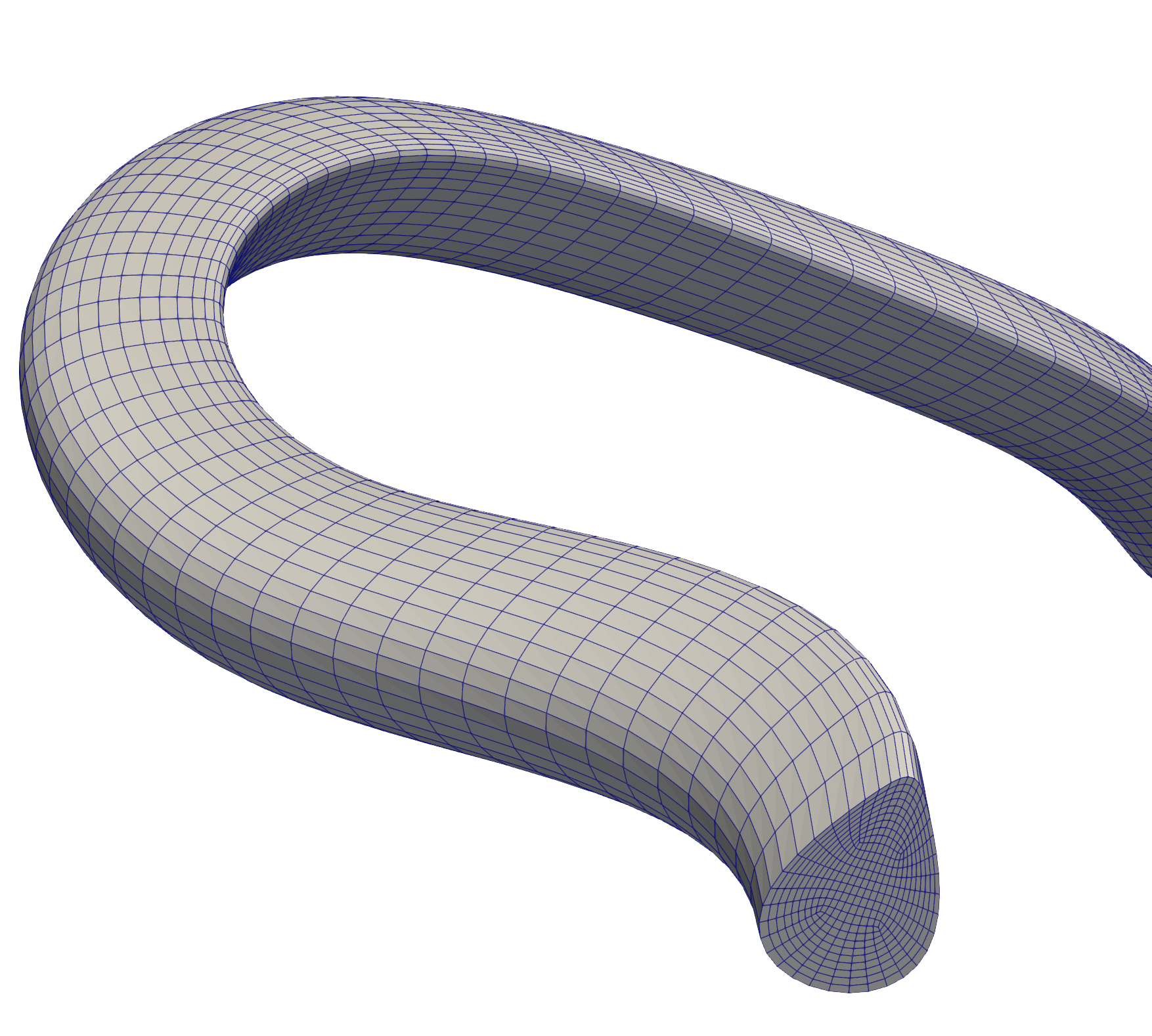} \quad\quad\quad
\includegraphics[width =7cm, keepaspectratio]{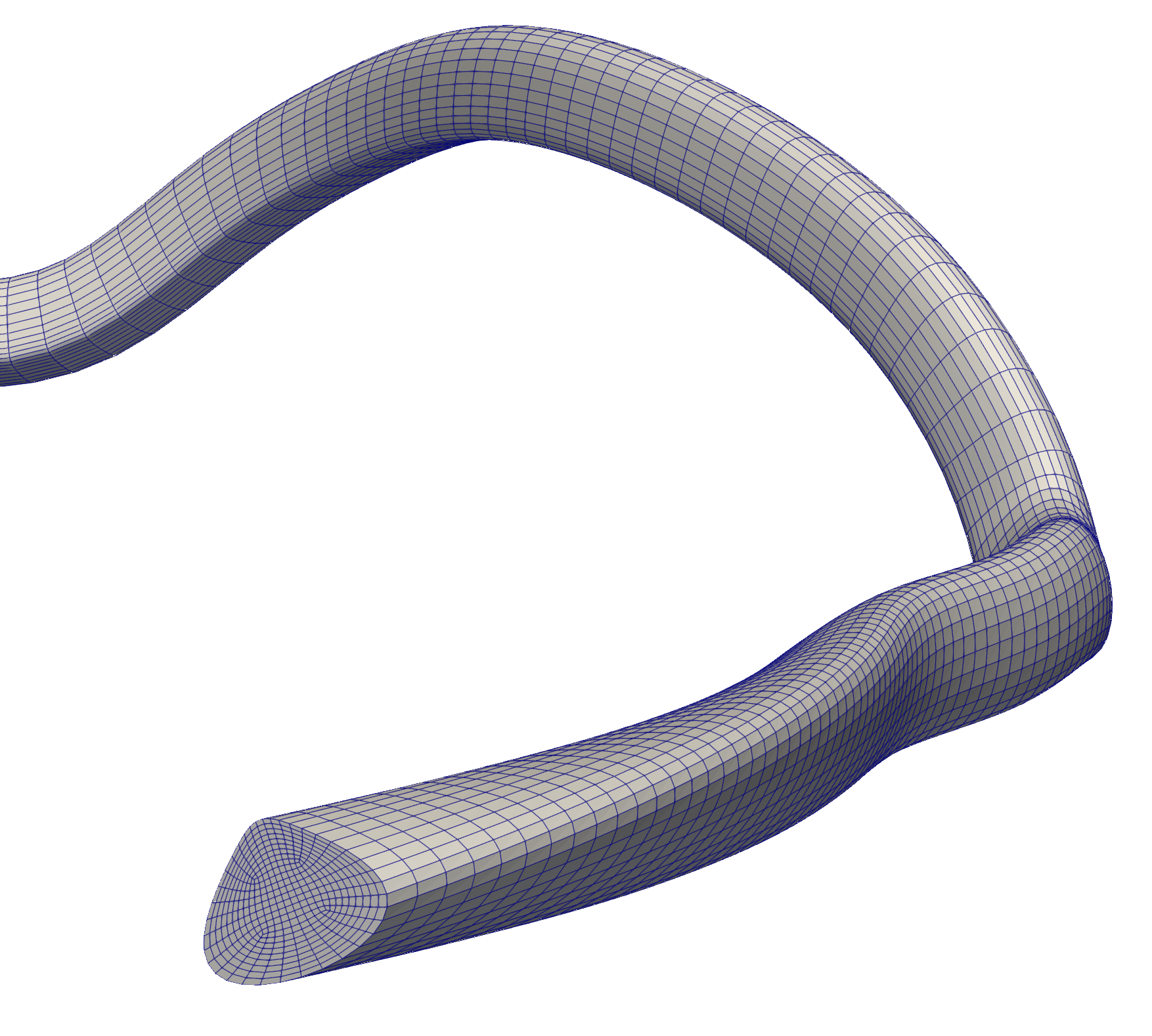}
\caption{A patient specific coronary artery with a general convex cross sections. For the sake of the representation, only the initial and end parts are depicted.}
\label{Fig: 3D_general_convex}
\end{figure}

Fig.\ \ref{Fig: harmonic_map} shows the spline map that results from mapping harmonically into non-circular convex domains along with a dense classical mesh sampled from the surface spline. Fig.\ \ref{Fig: 3D_general_convex} shows the volumetric mesh that results from a large number of non-circular cross-sectional meshes. 
\textcolor{black}{
The centerline curve will at every point $t \in [0, 1]$ of the parametric domain coincide with the center of mass of the cross-sectional disc placed along the vessel at the same value of $t$.}

In the following, we summarize all potential advantages of the proposed methodology:
\begin{itemize}
\item The ability to determine the quantity of vessel sections to be inserted along the centerline allows direct control over the number of vertices. 
Additionally, considering the generation of multiple grids, it becomes feasible to maintain the same topology across all the meshes in view of integrating a reduced basis framework.



\item The mesh is built using a spline approach in both the radial and longitudinal directions, leaving open the possibility of performing isogeometric analysis (IGA) simulations. Moreover, the use of spline control points results in an efficient and light representation of the mesh.

\item The implementation of flow-oriented cells improves the computational efficiency, see \cite{de2010patient,ghaffari2017large}; 


\item The use of a harmonic map allows for the parametrization of any convex cross sectional shape without the need to adjust the methodology based on the type of cross section.
\end{itemize}

We note that a generalisation of harmonic maps to nonconvex shapes exists but leads to a nonlinear problem \cite{azarenok20112d}. 
\textcolor{black}{
However, as reported in \cite{galassi20183d}, for coronary arteries, healthy tracts and uniform stenoses are accurately modeled using nearly circular cross-sections. Meanwhile, if the distribution of the stenosis is less uniform, shapes could be semi-circular, oval or even triangular, which are all convex shapes, i.e.\ they can be modeled with our approach.}

\subsection{Non-planar bifurcations} \label{sec: bifurcation}

While a number of studies on the generation of structured grids inside a vessel bifurcation have appeared in the literature, the proposed methods require a concrete reconstruction of the vessel surface \cite{zhang2007patient,de2010patient,de2011patient,decroocq2023modeling} or iterative post processing (usually, grid smoothing) to finalize the grid \cite{ghaffari2017large, decroocq2023modeling}.

In the present section, starting from the bifurcation geometrical models introduced by \cite{zakaria2008parametric}, and the extended versions of \cite{ghaffari2017large, decroocq2023modeling}, we propose a technique for generating the bifurcation model and its numerical grid exclusively based on Hermite curves. Our approach does not need the computation of the vessel surfaces or their intersections (a step that constitutes a critical robustness bottleneck) and other smoothing post-processing steps. 
In particular, we propose a method for meshing bifurcation geometries from no more than the three end-vessel sections' positional and contour information. 

For the sake of the exposition, we restrict our attention to circular contours. However, the generalization to generic convex cross-sectional shapes, using the techniques from Sec.\ \ref{sec: general_convex}, is straightforward.

\subsubsection{Geometrical modeling}

\begin{figure}[t!]
\centering
\includegraphics[width = 6.7cm, keepaspectratio]{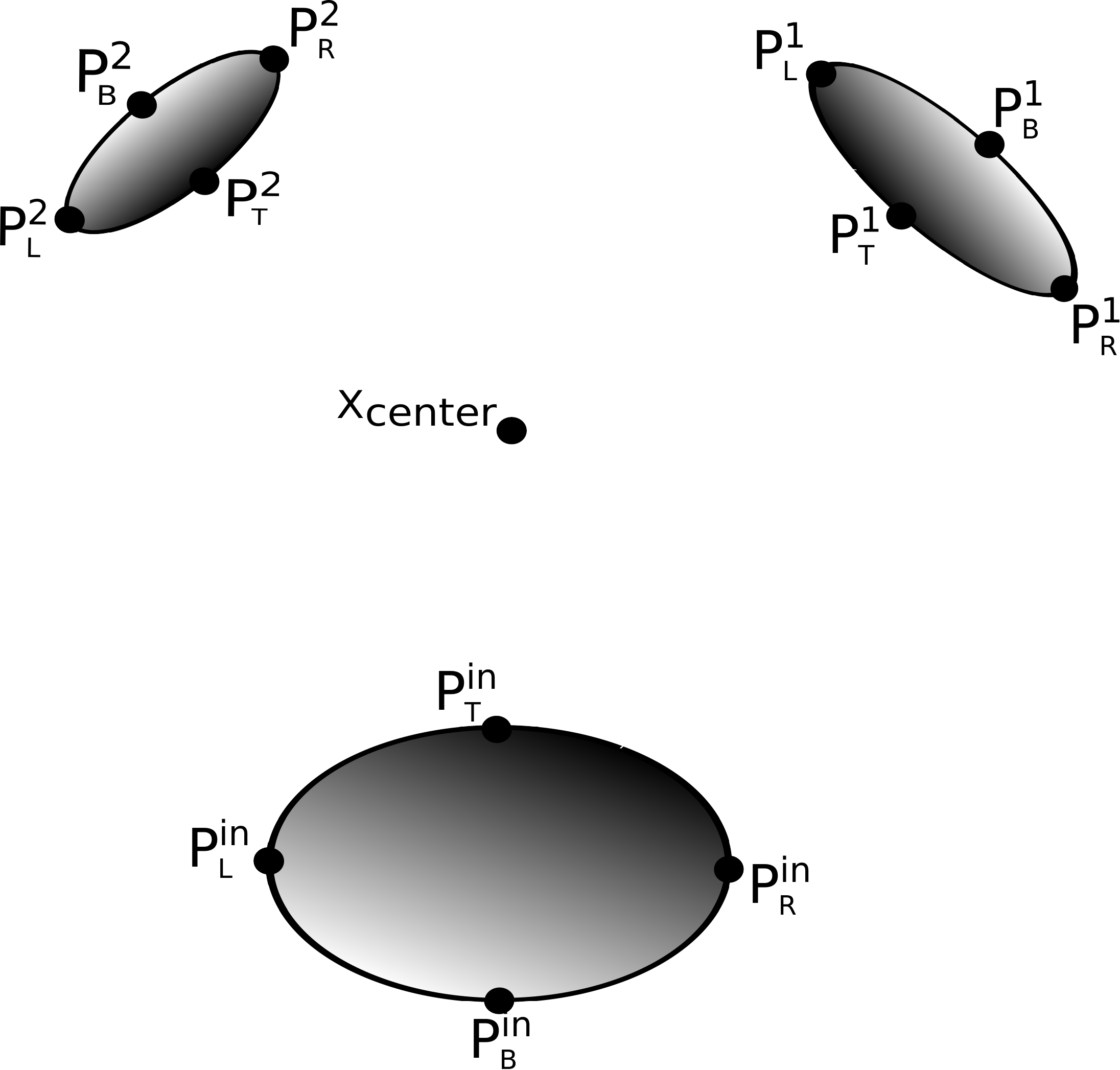}\quad\quad\quad\quad
\includegraphics[width = 6.7cm, keepaspectratio]{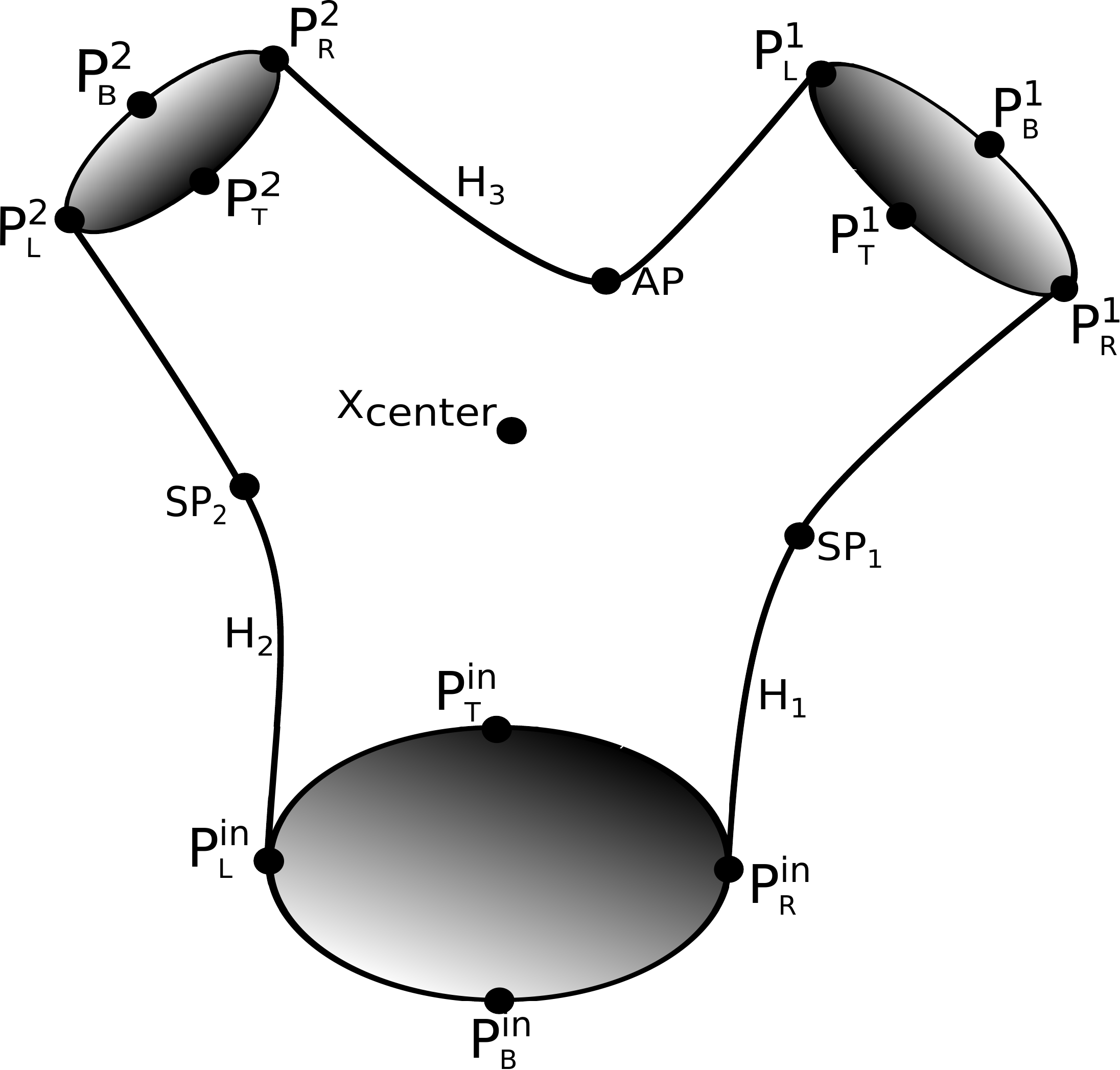}\vspace{1cm}
\includegraphics[width = 7.5cm, keepaspectratio]{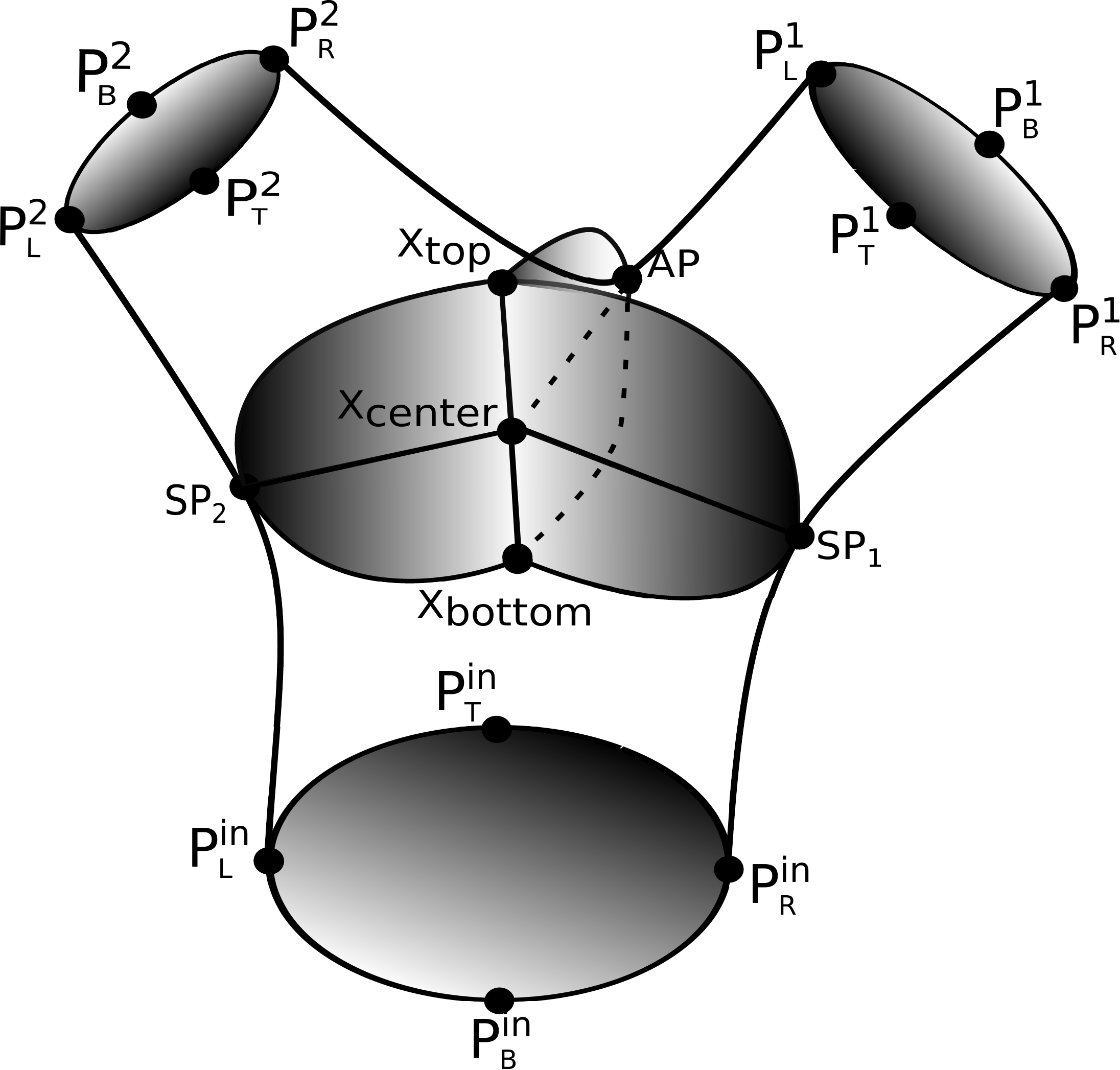}
\caption{Top-Left: The three contours are positioned without twisting with respect to the others. The anchor points (left, right, bottom, top) are highlighted in black.
The center points are not reported for the sake of the representation.
Top-Right: Construction of the Hermite curves ($\operatorname{H_{1}}$, $\operatorname{H_{2}}$, $\operatorname{H_{3}}$) to compute $\operatorname{SP_1}$, $\operatorname{SP_2}$, $\operatorname{AP}$.
Bottom: The contours of the three petals of the butterfly structure is depicted. For the sake of the representation, the Hermite curves to compute
$\operatorname{X_{top}}$, $\operatorname{X_{center}}$ and $\operatorname{X_{bottom}}$ are not reported.}
\label{Fig: bifurcation}
\end{figure}

\noindent We are given the circular inlet section $\Omega_{\text{in}} \subset \mathbb{R}^3$ (see Fig.\ \ref{Fig: bifurcation}, top-left) whose boundary is parameterized by a given correspondence $\mathbf{F}_{\text{in}}: \partial \widehat{\Omega}^r \rightarrow \partial \Omega_{\text{in}} \subset \mathbb{R}^3$.
Referring again to Fig.\ \ref{Fig: bifurcation}, top-left, we associate this contour with four anchor points, $P_R^{\text{in}}$, $P_T^{\text{in}}$, $P_L^{\text{in}}$, $P_B^{\text{in}}$ (right, top, left, bottom) which are, in order, assumed at $\mathbf{F}(\theta = \theta_i)$ with $\theta_i = \frac{i \pi}{2}, \, n \in \{1, 2, 3, 4\}$, where $\theta$ denotes the azimuth angle on $\partial \Omega^r$. Additionally, we introduce the cross section's center of mass $P_C^{\text{in}}$.
Similarly, we are given two outlet cross sections, $\Omega_1$ and $\Omega_2$, with associated correspondences $\mathbf{F}_{i}: \partial \widehat{\Omega}^r \rightarrow \partial \Omega_i \subset \mathbb{R}^3$ and centers of mass $P_C^i$, $i \in \{1, 2\}$ (for the sake of the representation, $P_C^i$ is not reported in Fig.\ \ref{Fig: bifurcation}, top-left). 

The cross sections $\Omega_{\text{in}}, \Omega_1, \Omega_2$ are associated with the inward-oriented unit tangent vectors denoted by $\mathbf{t}_{\text{in}}$, $\mathbf{t}_1$ and $\mathbf{t}_2$ which represent the normals to the cross sections' embedding planes. In the following, we reasonably assume $\mathbf{t}_1$ and $\mathbf{t}_2$ to be linearly independent.
To create a torsion-minimizing parameterization, we require the outlet cross sections' ($\Omega_1$ and $\Omega_2$) major axes to be oriented in a twist-minimizing fashion in relation to those of $\Omega_{\text{in}}$, see Fig. \ref{Fig: bifurcation}, top-left. For this, we compute the cross product $\mathbf{t_{1}} \times \mathbf{t_{2}}$ and project the result onto the inlet cross section's embedding plane. The resulting vector is denoted by $\mathbf{n}_{\text{in}}$. The binormal vector $\mathbf{b}_{\text{in}}$ is computed as the normalized cross product of $\mathbf{t}_{\text{in}}$ and $\mathbf{n}_{\text{in}}$.
In the unusual case where $\mathbf{t}_1$ and $\mathbf{t}_2$ are parallel, the normal vector $\mathbf{n}_{\text{in}}$ can be computed as the projection of the perpendicular vector to the plane spanned by the three anchor points $P_C^{\text{in}}$, $P_C^{\text{1}}$ and $P_C^{\text{2}}$ onto the inlet cross section's embedding plane. The vector triplet $(\mathbf{t}_{\text{in}}, \mathbf{n}_{\text{in}}, \mathbf{b}_{\text{in}})$ creates the orthonormal frame $R_{\text{in}}$ for the inlet's embedding plane. The associated outlet frames $R_1$ and $R_2$ follow from computing a discrete twist minimization (see Section \ref{sec: single_branch}) on the point sequence $(P^{\text{in}}_C, P_C^1, P_C^2)$ using $R_{\text{in}}$ as initial configuration on $P_C^{\text{in}}$. 
In this manner, as depicted in Fig. \ref{Fig: bifurcation}, top-left, the grids' coordinate frames exhibit minimal mutual torsion, as shown through the consistent orientation of the anchor points. 
In what follows, we denote the anchor points on the $\Omega_i$, as dictated by the RMF, by $(P_L^{i}$, $P_R^{i}$, $P_B^{i}$, $P_T^{i})$. The frames induced at $(\Omega_{\text{in}}, \Omega_1, \Omega_2)$ serve as initial conditions for the RMF of each individual branch.


We introduce the Hermite curve as $C: [0, 1] \rightarrow \mathbb{R}^3$ between points $p \in \mathbb{R}^3$ and $q \in \mathbb{R}^3$ as follows:
\begin{align}
C(t) = \operatorname{H_{0, 0}}(t) p + \operatorname{H_{0, 1}}(t) \|p - q\| \mathbf{m}_1 + \operatorname{H_{1, 0}}(t) q + \operatorname{H_{1, 1}}(t) \|p - q\| \mathbf{m}_2,
\label{eq: curve_hermite}
\end{align}
where $\{ \mathbf{m}_1, \mathbf{m}_2 \}$ denote the points' normalized tangent vectors. Expressions for the $\operatorname{H_{i, j}}(t)$ can be found in \cite{farin2002curves}.
As depicted in Fig. \ref{Fig: bifurcation}, top-right, we construct a total of three Hermite curves: $\operatorname{H_{1}}$, $\operatorname{H_{2}}$ and $\operatorname{H_{3}}$, connecting the point pairs $(P_R^{\text{in}}, P_R^1)$, $(P_L^\text{in}, P_L^2)$ and $(P_R^2, P_L^1)$, respectively. As before, we impose that the curve assume the embedding planes' tangents scaled by the Euclidean distance between the boundary points. For example, for $\operatorname{H_{1}}$, we impose $\mathbf{m}_1 = \| P_R^{\text{in}} - P_R^1 \| \mathbf{t}_{\text{in}}$ and similarly for $\mathbf{m}_2$.
We evaluate the Hermite curves in $t = 0.5$, and following the notation proposed in \cite{decroocq2023modeling}, points $\operatorname{SP_1}$, $\operatorname{SP_2}$, $\operatorname{AP}$ are directly computed as follows:

\begin{equation*}
\operatorname{SP_1} = \operatorname{H_{1}}(0.5), \quad \operatorname{SP_2}= \operatorname{H_{2}}(0.5), \quad \operatorname{AP} = \operatorname{H_{3}}(0.5).
\end{equation*}

The use of Hermite interpolation allows us to avoid the apex smoothing procedure, presented in \cite{decroocq2023modeling}, where the apical point ($\operatorname{AP}$) of the bifurcation is computed by means of the intersection of the two vessel surfaces.

A similar procedure is repeated for the top, center and bottom anchor points of the three grids.
In particular, we create Hermite curves between $P_T^\text{in}$ and $P_T^1$ ($\operatorname{H_{4}}$), $P_T^\text{in}$ and $P_T^2$ ($\operatorname{H_{5}}$) and between $P_T^2$ and $P_T^1$ ($\operatorname{H_{6}}$).
The same procedure is repeated for the center anchor points creating 
$\operatorname{H_{7}}$, $\operatorname{H_{8}}$, $\operatorname{H_{9}}$, and for the bottom anchor points leading to $\operatorname{H_{10}}$, $\operatorname{H_{11}}$, $\operatorname{H_{12}}$.
A total of 12 Hermite curves defining the bifurcation are obtained. 
Using the generated Hermite curves, we define three new points:

\begin{align*}
\operatorname{X_{top}} = \displaystyle\frac{\operatorname{H_{4}}(0.5) + \operatorname{H_{5}}(0.5) + \operatorname{H_{6}}(0.5)}{3},\\
\operatorname{X_{center}} = \displaystyle\frac{\operatorname{H_{7}}(0.5) + \operatorname{H_{8}}(0.5) + \operatorname{H_{9}}(0.5)}{3},\\
\operatorname{X_{bottom}} = \displaystyle\frac{\operatorname{H_{10}}(0.5) + \operatorname{H_{11}}(0.5) + \operatorname{H_{12}}(0.5)}{3}.
\end{align*}

We note that the points $\operatorname{X_{top}}$, $\operatorname{X_{center}}$ and $\operatorname{X_{bottom}}$ are computed using the same strategy and they lie on the same line by construction. 
As reported in Fig.\ \ref{Fig: bifurcation}, bottom, the point $\operatorname{SP_1}$ is connected to $\operatorname{X_{top}}$ via the bounding curve of one ellipse quadrant.
This quadrant is defined by the points $\operatorname{SP_1}$, $\operatorname{X_{top}}$ and $\operatorname{X_{center}}$.
Following this geometric framework, we extend the connections to include all relevant points: $\operatorname{AP}$, $\operatorname{SP_1}$, $\operatorname{SP_2}$, $\operatorname{X_{top}}$, $\operatorname{X_{center}}$, and $\operatorname{X_{bottom}}$. This process allows us to construct the contours of an ad-hoc \textit{butterfly structure}, completing the skeleton of the bifurcation, see Fig.\ \ref{Fig: bifurcation}, bottom.

\textcolor{black}{The choices of evaluating the Hermite curves at $t = 0.5$, and the parameter used for scaling the tangent magnitudes are arbitrary.
The Hermite curves’ tangents are scaled with case-specific parameters, such as the distances between the
points or the diameter of the cross section; this avoids distortion or wiggles, since the curve automatically adapts to the specific case.
However, the user may always tune both the magnitude of the tangents and the choice of the sampling point of the Hermite curves, to differently shape the bifurcation. 
}

\subsubsection{Meshing}

\begin{figure}[t!]
\centering
\includegraphics[width = 16cm, keepaspectratio]{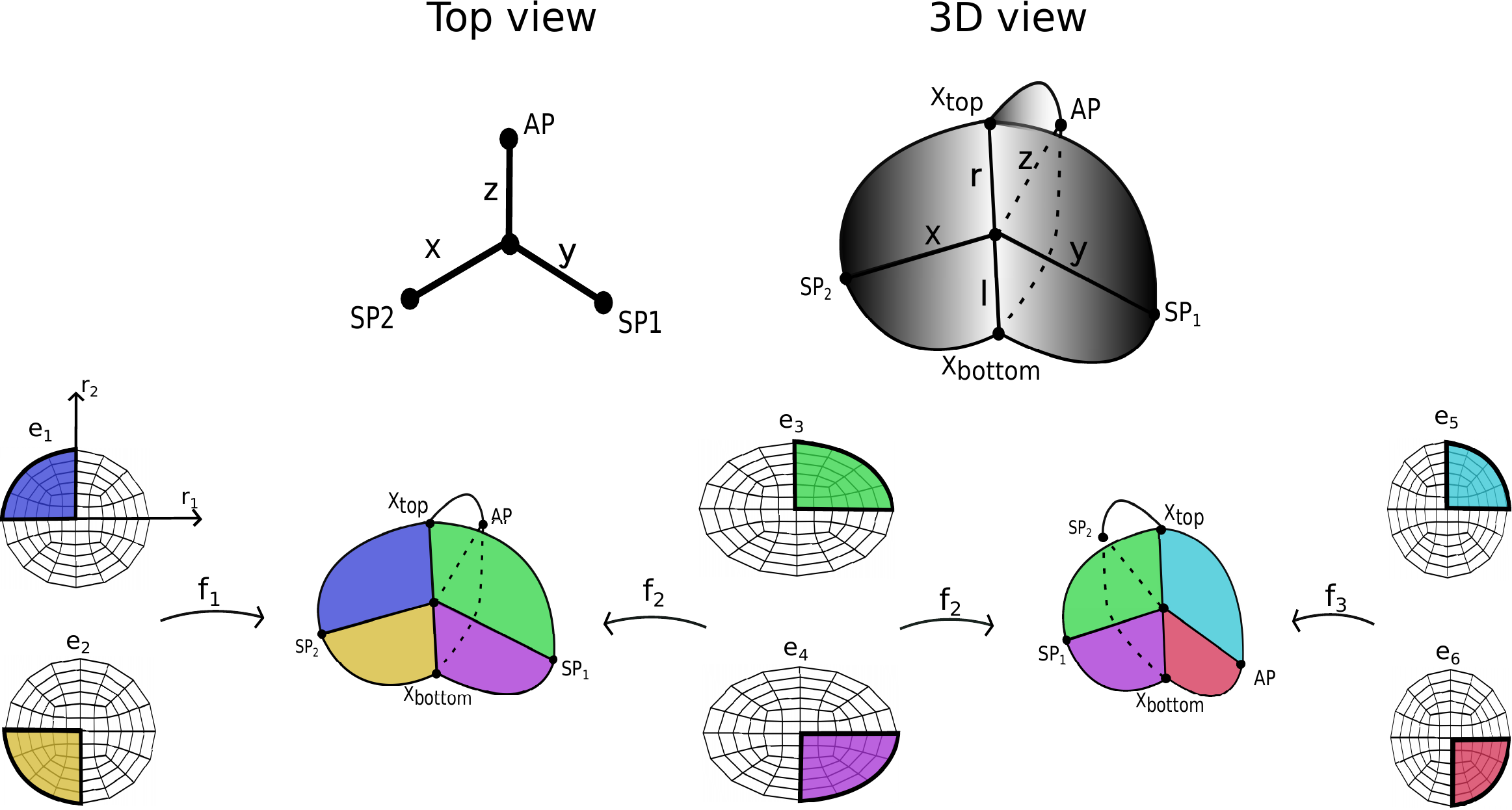}
\caption{Construction of the spline map inside the butterfly structure.}
\label{Fig: butterfly_construction}
\end{figure}

Once the creation of the bifurcation's skeleton has been finalized, in view of obtaining a conformal mesh suited for both FEM and IGA, we assume that the three end-vessel section patches have the same knot vectors. This condition is not restrictive, as we have the flexibility to decide how to generate the meshes, and it is not imposed by any external party.
While the grid generation inside the vessels has been already discussed in Sec.\ \ref{sec: single_branch}, the grid inside the butterfly structure requires more attention.

Following Fig.\ \ref{Fig: butterfly_construction}, top, we define the following vectors: 
\begin{align*}
& \mathbf{x} := \operatorname{SP_2} - \operatorname{X_{center}}, \quad \mathbf{y} := \operatorname{SP_1} - \operatorname{X_{center}}, \quad \mathbf{z} := \operatorname{AP}
 - \operatorname{X_{center}}, \\
& \mathbf{r} := \operatorname{X_{top}} - \operatorname{X_{center}}, \quad \mathbf{l} := \operatorname{X_{bottom}} - \operatorname{X_{center}}.
\label{eq: vector_butterfly1}
\end{align*}



By $\operatorname{e}(a, b) \subset \mathbb{R}^2$, we denote a spline map defined inside a planar ellipse with semi axis lengths of $a$ and $b$ aligned with the $r_1$ and $r_2$ axis, respectively.
By $S = \{\operatorname{e}(a, b)| \, r_1\leq0,r_2\leq0\}$, we denote the subset contained in the ellipse's third quadrant where both coordinates are negative and similar for the other quadrants.
Using the spline-based approach presented in Sec.\ \ref{sec: plane_mesh} and denoting $\norm{\, \cdot \,}$ as the Euclidean norm of a vector, we compute the following spline maps (see Fig.\ \ref{Fig: butterfly_construction}, bottom):

\begin{align*}
& \operatorname{e_1} = \{\operatorname{e}(\norm{\mathbf{x}},\norm{\mathbf{r}})| \, r_1\leq0,r_2\ge0) \}, 
\operatorname{e_3} = \{\operatorname{e}(\norm{\mathbf{y}},\norm{\mathbf{r}})| \, r_1\ge0,r_2\ge0)\ \}, 
\operatorname{e_5} = \{\operatorname{e}(\norm{\mathbf{z}},\norm{\mathbf{r}})| \, r_1\ge0,r_2\ge0)\ \}, \\ \vspace{0.5cm}
& \operatorname{e_2} = \{\operatorname{e}(\norm{\mathbf{x}},\norm{\mathbf{l}})| \, r_1\leq0,r_2\leq0)\ \}, 
\operatorname{e_4} = \{\operatorname{e}(\norm{\mathbf{y}},\norm{\mathbf{l}})| \, r_1\ge0,r_2\leq0)\ \}, 
\operatorname{e_6} = \{\operatorname{e}(\norm{\mathbf{z}},\norm{\mathbf{l}})| \, r_1\ge0,r_2\leq0)\ \} .
\end{align*}

We introduce three rotation matrices:

\begin{equation*} 
\operatorname{R_1} := [\mathbf{x}, \mathbf{r}], \quad 
\operatorname{R_2} := [\mathbf{y}, \mathbf{r}], \quad
\operatorname{R_3} := [\mathbf{z}, \mathbf{r}]
\label{eq: rotation_matrix}
\end{equation*}

acting on the 2D ellipse quadrant in a roto-translational fashion, $\operatorname{f} : \mathbb{R}^2 \rightarrow \mathbb{R}^3 $, as follows:
\begin{align}
    \operatorname{f_1} = \operatorname{R_1}  \operatorname{e_j} + \operatorname{X_{center}}, \quad \operatorname{j}\, =\, 1,2,  \\
    \operatorname{f_2} = \operatorname{R_2} \operatorname{e_j} + \operatorname{X_{center}}, \quad \operatorname{j}\, =\, 3,4 , \\
    \operatorname{f_3} = \operatorname{R_3} \operatorname{e_j} + \operatorname{X_{center}}, \quad \operatorname{j}\, =\, 5,6. 
\end{align}

In particular, continuing to follow Fig.\ \ref{Fig: butterfly_construction}, bottom,  we roto-translate the top-left quadrant $\operatorname{e_1}$ and the bottom-left quadrant $\operatorname{e_2}$ applying $\operatorname{f_1}$.
Note that the anchor point "top" of $\operatorname{e_1}$ coincides with $\operatorname{X_{top}}$, the "center" anchor points of both $\operatorname{e_1}$ and $\operatorname{e_2}$ coincide with $\operatorname{X_{center}}$, the "bottom" anchor point of $\operatorname{e_2}$ coincides with $\operatorname{X_{bottom}}$ and, finally, the "left" anchor points of both $\operatorname{e_1}$ and $\operatorname{e_2}$ coincide with $\operatorname{SP_2}$.
Using the same procedure, we translate the top-right quadrant $\operatorname{e_3}$ and the bottom-right quadrant $\operatorname{e_4}$ using $\operatorname{f_2}$.
Note that, as before, the anchor point "top" of $\operatorname{e_3}$ coincides with $\operatorname{X_{top}}$, both of the ellipses' "center" points coincide with $\operatorname{X_{center}}$, the "bottom" of $\operatorname{e_4}$ coincides with $\operatorname{X_{bottom}}$ and, additionally, the "right" anchor points coincide with $\operatorname{SP_1}$.
Finally, using $\operatorname{f_3}$, we position the spline map of the third petal of the butterfly structure.

Once the generation of the spline maps inside the butterfly structure has been finalized, the volumetric Hermite interpolation between the vessel sections and the butterfly structure requires the definition of a space-dependent tangent vector, $\mathbf{t_\text{BF}}(\mathbf{x})$, assigned to the butterfly structure. 
We denote $\Omega^\text{in}_\text{BF}$ as the domain of the butterfly structure containing the two petals pointing towards $\Omega_\text{in}$ (the ones composed by $e_1, e_2, e_3$ and $e_4$, see Fig.\ \ref{Fig: butterfly_construction}).
For example, we consider the interpolation between $\Omega_\text{in}$ and $\Omega^\text{in}_\text{BF}$. The tangent for all the points contained in $\Omega_\text{in}$ is a constant vector, $\mathbf{t_{\text{in}}}$. 
While in the case of $\Omega^\text{in}_\text{BF}$, we define a space-dependent tangent, i.e., $\mathbf{t}_{\text{BF}} = \mathbf{t}_{\text{BF}}(r_1, r_2)$ in the parametric disc $\widehat{\Omega}^r$. It reads:
\begin{align}
\mathbf{t}_{\text{BF}}(r_1, r_2) = P^3(f_0(r_1) \operatorname{H_1^\prime}(0.5) + f_1(r_1) (X_{\text{center}} - P_c^{\text{in}}) + f_2(r_1) \operatorname{H_2^\prime}(0.5))
\label{eq: tangent_butterfly}
\end{align}

\begin{figure}[t!]
\centering
\includegraphics[width = 5.8cm, keepaspectratio]{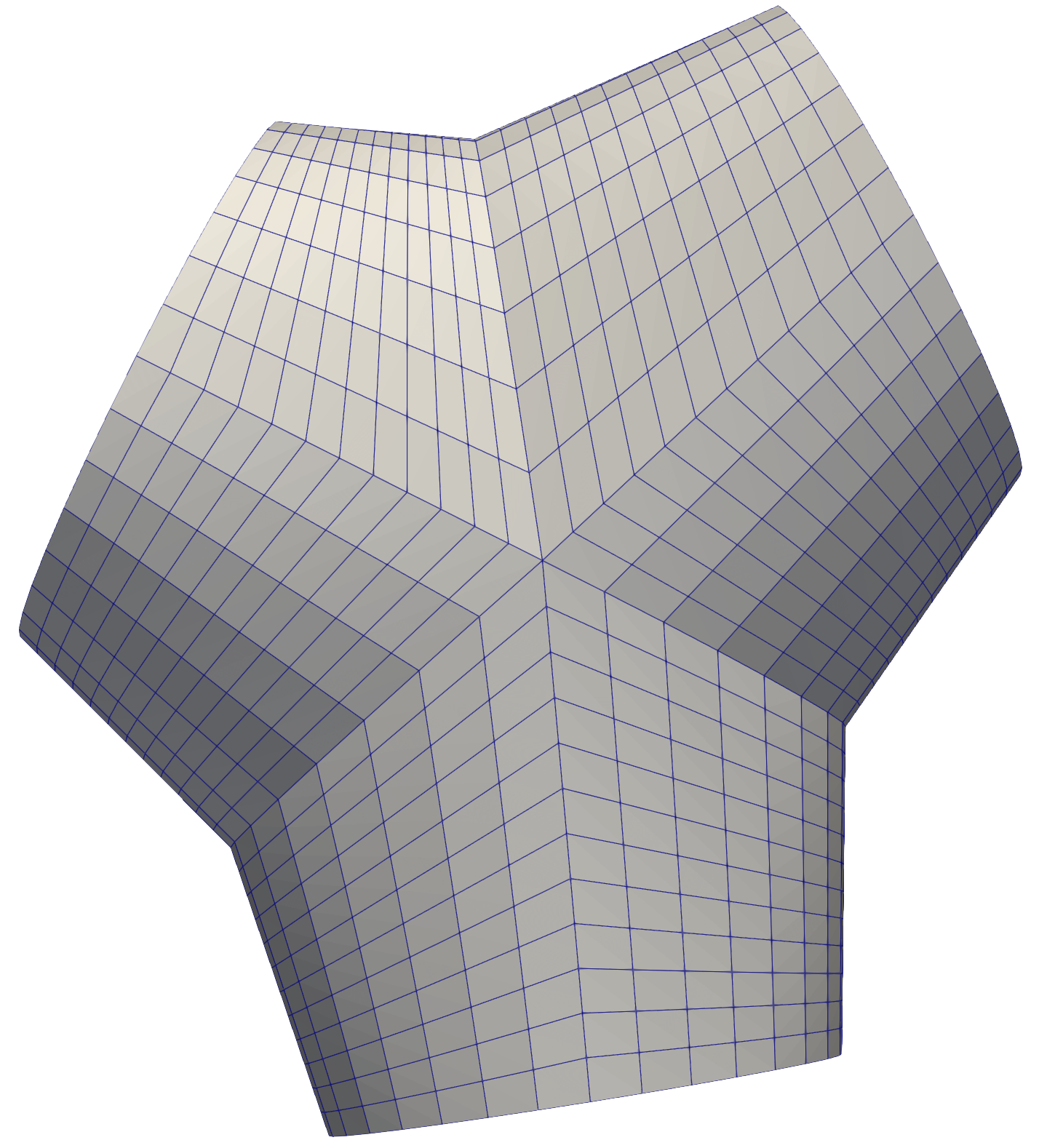}
\includegraphics[width = 7cm, keepaspectratio]{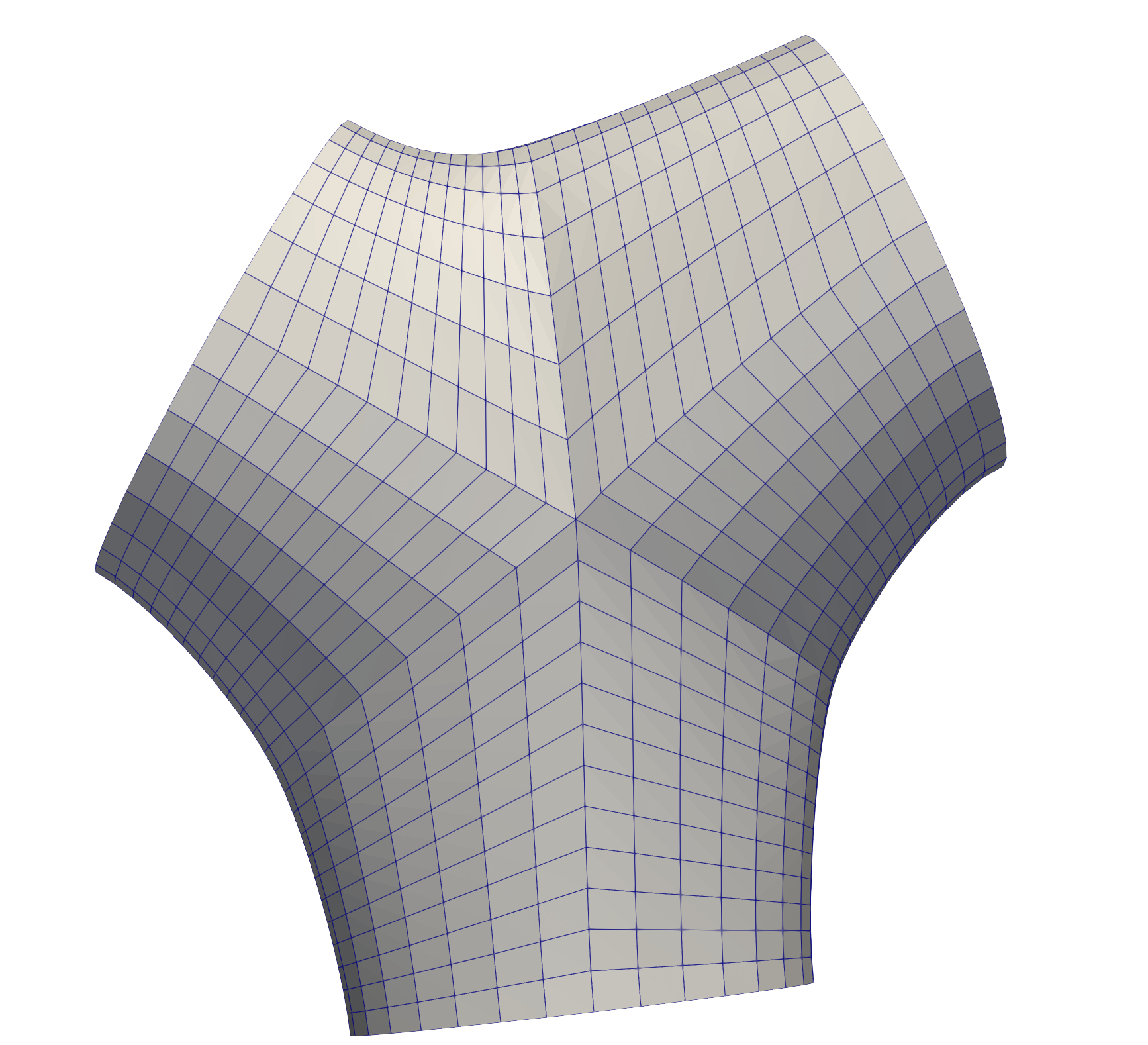}
\caption{Left: Linear interpolation. Right: Hermite interpolation.}
\label{Fig: linear_vs_Hermite}
\end{figure}

where $f_0, f_1, f_2$ are piecewise linear hat functions satisfying $f_0(-1) = 1, f_1(0) = 1, f_2(1) = 1$ (and zero in the other two), while $P^3(\cdot)$ denotes the spline projection onto $(\mathcal{V}_h^{\mathbf{r}})^3$, c.f. Eq.~\eqref{eq: projection}. 
Note that $\operatorname{H^\prime}$ denotes the derivative of the Hermite curve, i.e., the tangent. Eq.\ \ref{eq: tangent_butterfly} results in a continuous transition along the separation line (containing $\operatorname{X_{top}}$,$\operatorname{X_{center}}$, $\operatorname{X_{bottom}}$).

The domain $\Omega_\text{{in}}$ is glued to $\Omega^\text{in}_\text{BF}$ through the following volumetric Hermite interpolation, $V: \widehat{\Omega}^r \times [0, 1] \rightarrow \mathbb{R}^3$:
\begin{align}
V(r_1, r_2, t) = \operatorname{H_{0, 0}}(t) p(r_1,r_2) + \operatorname{H_{0, 1}}(t)\mathbf{t}_{in} + \operatorname{H_{1, 0}}(t) q(r_1,r_2) + \operatorname{H_{1, 1}}(t)\mathbf{t}_{BF}(r_1, r_2), 
\label{eq: vol_interpolation}
\end{align}

where $p \in \Omega_\text{in}$ and $q \in \Omega^\text{in}_\text{BF}$. 
The map $V: \widehat{\Omega}^r \times [0, 1] \rightarrow \mathbb{R}^3$ is a volume (result of linear combinations of functions described by splines) representing a continuous extrusion from $\Omega_\text{{in}}$ to $\Omega^\text{in}_\text{BF}$. For given $t = t_0$, it corresponds to a unique vessel section between $\Omega_\text{{in}}$ and $\Omega^\text{in}_\text{BF}$. In particular, $V(r_1,r_2,t = 0)$ corresponds to $\Omega_\text{in}$, while $V(r_1,r_2,t = 1)$ coincides with $\Omega^\text{in}_\text{BF}$.
We note that the spline-based description resulting from the $V(r_1, r_2, t)$ can be converted to an arbitrarily dense mesh via sampling.
A possible linear and Hermite-based meshes are depicted in Fig. \ref{Fig: linear_vs_Hermite}.

In the following, we summarize the main steps:
\begin{enumerate}
\item Position of the three vessel sections while minimizing their mutual torsion. 
\item Creation of the Hermite curves between the right and left anchor points of the sections, evaluate the curves in the center of the parametric domain to find the $\operatorname{AP}$, $\operatorname{SP_1}$ and $\operatorname{SP_2}$ points.
\item Creation of the Hermite curves between the top, center and bottom anchor points of the sections, then average the points that result from evaluating these curves (always in in the center of the parametric domain) to find $\operatorname{X_{top}}$, $\operatorname{X_{center}}$ and $\operatorname{X_{bottom}}$, respectively.
\item  Position the butterfly structure within the three sections.
\item Perform a linear/Hermite interpolation between the three sections and the butterfly structure to obtain a mesh ready for numerical simulation.
\end{enumerate}

We note that this procedure requires no more than the positioning of three grids in $3$D space while no explicit vessel surfaces or centerline information is needed.

\begin{figure}[t]
\centering
\includegraphics[width = 6.8cm, keepaspectratio]{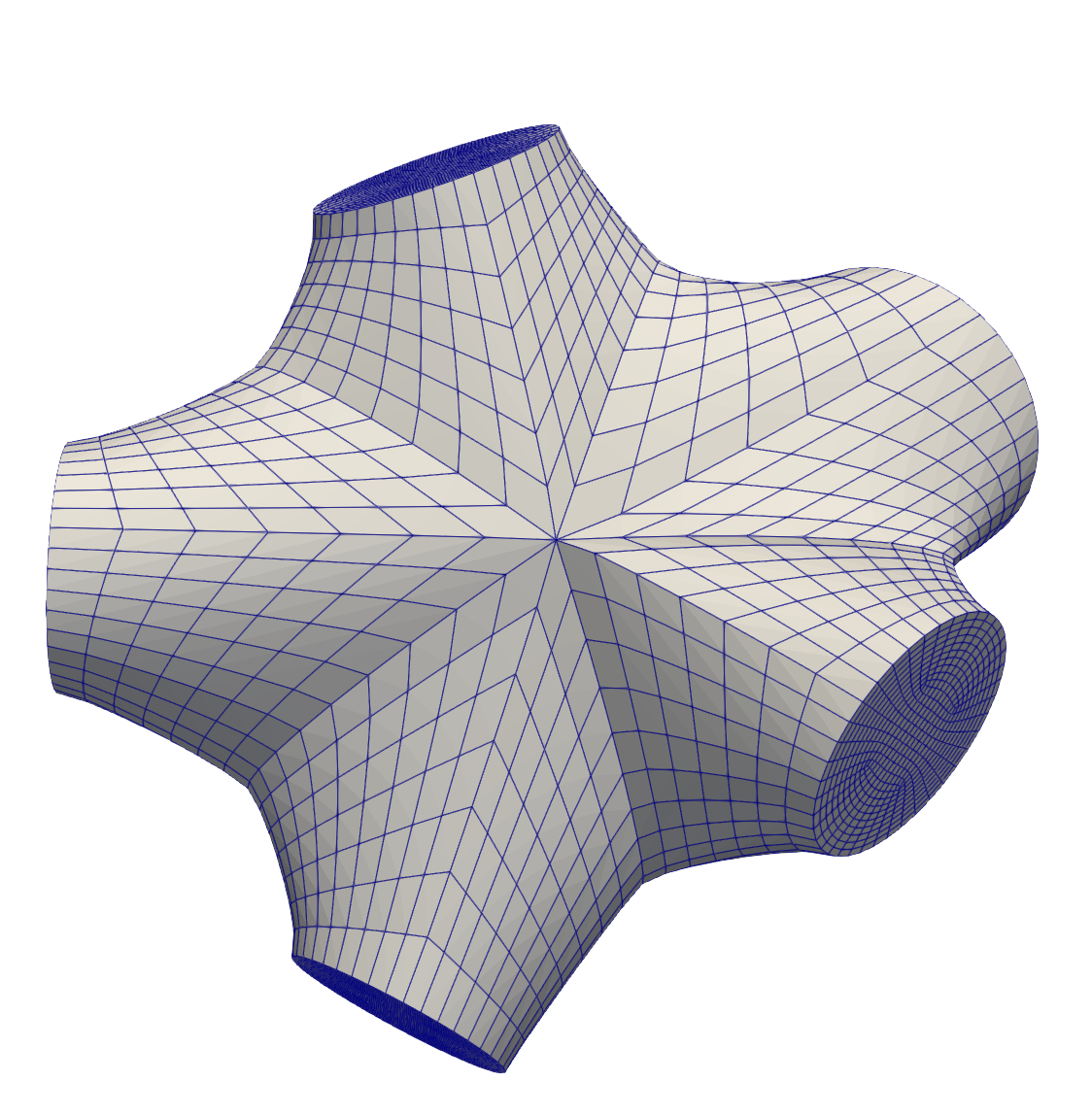} \quad
\includegraphics[width = 8cm, keepaspectratio]{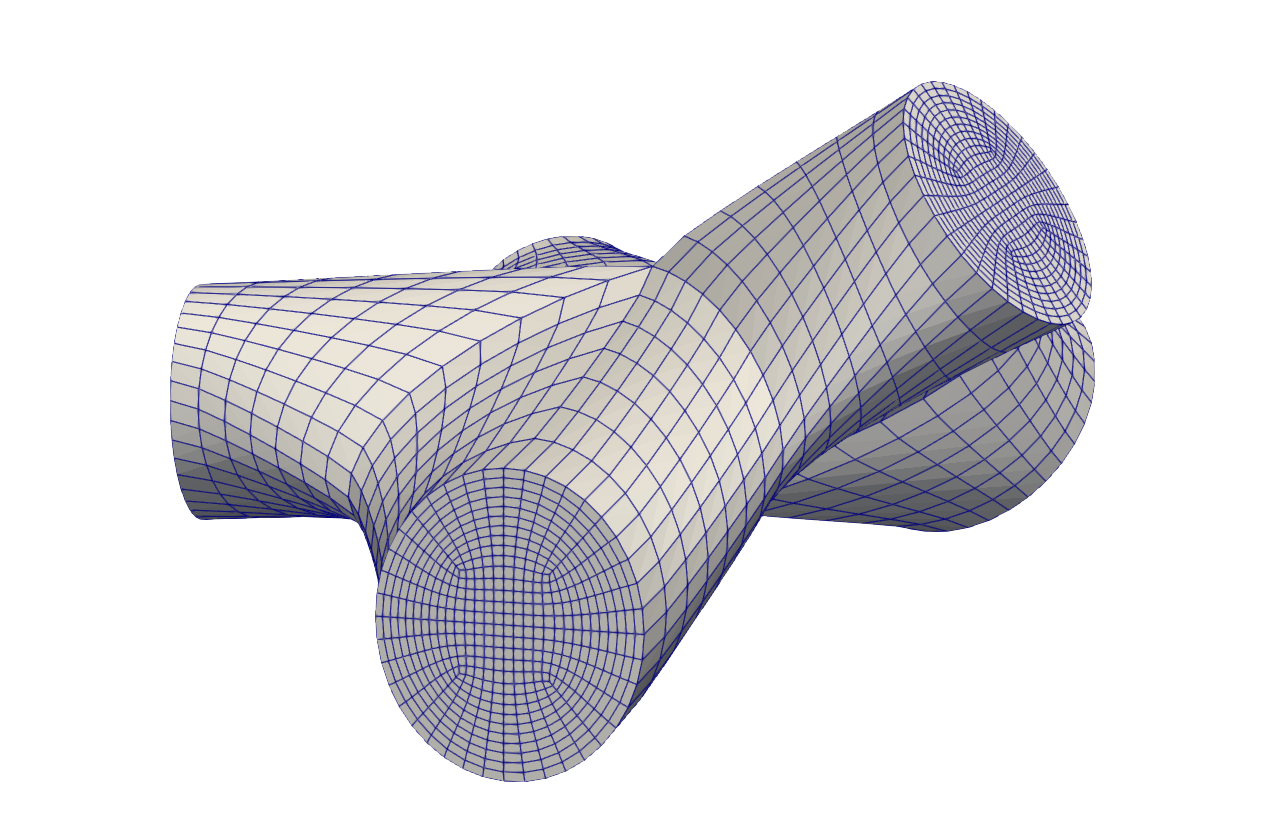} 
\caption{Left: Top view of 5 non-planar intersecting branches. Right: Lateral view of non-planar 5 intersecting branches. Note that the vessel sections can have different radius and tangents.}
\label{Fig: n_furcations}
\end{figure}

\subsection{Non-planar n-furcations} \label{sec: n_furcations}
Taking advantage of the intrinsic properties of Hermite curves, the bifurcation algorithm generalizes straightforwardly to any number of intersecting branches.
As an extension of the approaches from \cite{ zhang2007patient,ghaffari2017large,decroocq2023modeling}, which require the ad-hoc introduction of a master branch, are incompatible with $n$-furcations or require them to be planar, our approach imposes no such limitations. 
Depending on the number of intersecting branches ($\operatorname{n_{b}}$), the required number of Hermite curves ($\operatorname{n_H}$) and the number of ellipse quadrants $\operatorname{n_e}$ is given by:

\begin{align*}
    \operatorname{n_{H}} = 4 * \operatorname{n_{b}}, \\
    \operatorname{n_{e}} = 2* \operatorname{n_{b}}.
\end{align*}

Considering the case of $5$ non-planar intersecting branches, for instance, building the butterfly structure requires computing a total of $20$ Hermite curves and $10$ ellipse quadrants, following the same procedure presented in Sec.\ \ref{sec: bifurcation} (which is not repeated here for the sake of brevity). Note that, in this case, the butterfly structure has $5$ petals, and moreover, 
\textcolor{black}{the $C^1$ continuity is lost at the intersection point.}
In Fig. \ref{Fig: n_furcations}, we depict the top and lateral view of the mesh with $5$ non-planar intersecting \textcolor{black}{branches obtained with Hermite interpolation}. 

\subsection{Boundary Layers}\label{sec: BL}



In fluid dynamics simulations, capturing the regions characterized by strong velocity gradients is crucial to obtain meaningful results. Typically, in h{\ae}modynamics, these regions are situated close to the vessel surfaces. 

In the context of this paper's methodology based on splines, a convenient method for creating a boundary layer is directly introducing it in the multipatch covering of the reference domain $\widehat{\Omega}^r$. In addition, the boundary layer will, at least in part, carry over to the harmonic map between $\widehat{\Omega}^r$ and $\Omega$ in the case of non-circular vessel sections (see Sec. \ref{sec: general_convex}). 

\begin{figure}[t]
\centering
\includegraphics[width = 15cm, keepaspectratio]{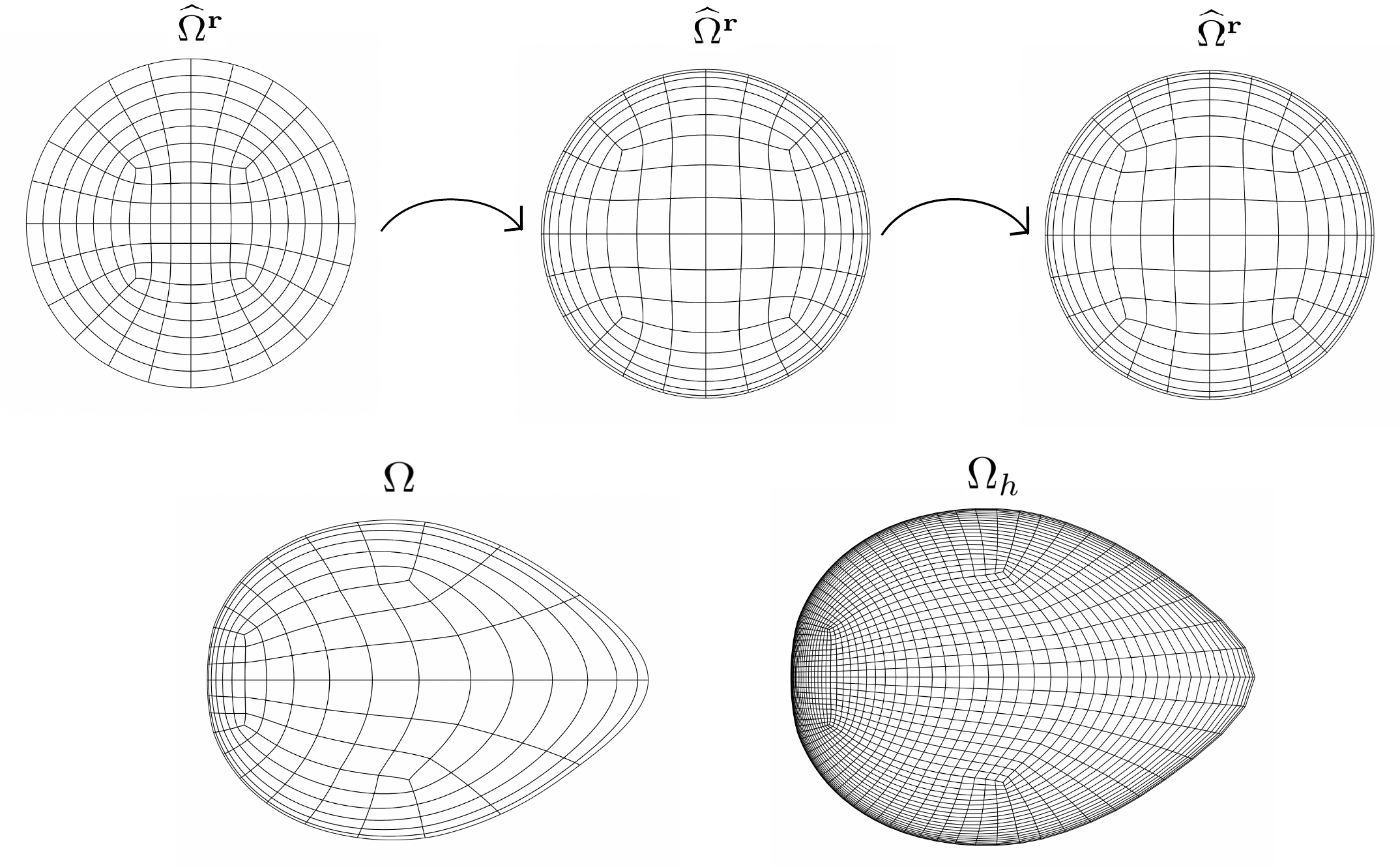}
\caption{Figure depicting, in order, $1)$ the original controldomain $\widehat{\Omega}^r$'s covering, $2)$ the boundary layer covering under the control function $f_\alpha(\mathbf{r})$ with $\alpha = 0.95$, $3)$ the covering whose four outer patch isolines have been projected back onto straight line segments, $4)$ the boundary layer parameterization of a convex non-circular cross-section and $5)$ a mesh sampled from the spline parameterization whose number of mesh points exceeds the dimension of $\mathcal{V}_h$ by one order of magnitude. }
\label{Fig: boundary_layer_sampling}
\end{figure}

A boundary layer is created in a two-step procedure. The first step introduces a global map from the default covering of $\widehat{\Omega}^r$ (see Figure \ref{Fig: boundary_layer_sampling}, left) onto itself. This map aims at accumulating control points close to $\partial \widehat{\Omega}^r$ (while leaving boundary control points intact). Denoting the free Cartesian coordinate functions in $\widehat{\Omega}^r$ by $(r_1, r_2)^T$, we achieve this by the introduction of a map
\begin{align}
    \mathbf{r}^\prime(\mathbf{r}) = f(r) \hat{\mathbf{r}}, \quad \text{with } r = \| \mathbf{r} \| \text{ and } \hat{\mathbf{r}} = \frac{\mathbf{r}}{r}.
\end{align}
Here, $f(r)$ satisfies $f(0) = 0$, $f(1) = 1$ and $f^\prime(r) > 0$. A boundary layer is the result of $f^\prime(r)$ being small in a neighbourhood of $r = 1$, whereby smaller values create a stronger layer. A possible choice of $f(r)$ is
\begin{align}
    f_\alpha(r) = r + \alpha \left(1 - r\right) r,
    \label{eq: bl}
\end{align}
where the parameter $0 < \alpha < 1$ tunes the boundary layer density. Note that $f_\alpha(0) = 0$, $f_\alpha(1) = 1$ and $f^\prime(r) > 0$ as required, while $f_\alpha^\prime(1) = 1 - \alpha > 0$. For an overview of suitable accumulation functions $f(\, \cdot \,)$ with varying properties, we refer to \cite[Chapter~4]{liseikin1999grid}. 

Since the reference controlmap's outer patch isolines do not align with straight rays drawn radially outward from $\widehat{\Omega}^r$'s center of mass, a strong boundary layer (i.e., $\alpha \approx 1$) may introduce a considerable curvature in the adjusted controlmap's isolines near the boundary $\partial \widehat{\Omega}^r$. We deem this undesirable from a parameterisation quality perspective. As a remedy, we project the outer patch's isolines back onto straight line segments between the interior patch's interface and the endpoint on $\partial \widehat{\Omega}^r$. This is most conveniently achieved by directly moving the controlmap's controlpoints onto the associated ray via an orthogonal projection. Hereby, the interior patch is not changed. 

Fig.\ \ref{Fig: boundary_layer_sampling} shows an example of boundary layer creation based on this strategy. The boundary layer created in $\widehat{\Omega}^r$ carries over to the convex, non-circular target cross section via the harmonic map. We mention that changing the boundary layer is associated with a non-negligible computational burden because the harmonic map's associated PDE problem has to be solved over $\widehat{\Omega}^r$ in the new coordinate system and therefore requires re-assembly of the associated linear operator. As such, once introduced, the boundary layer is not changed such that the same linear operator can be used for the computation of all cross sections.

\clearpage

\section{Numerical results}\label{sec: results}
The section of the results is divided into two main parts:
\begin{itemize}
\item Meshing results (Sec.\ \ref{sec: meshing_results}). In Sec.\ \ref{sec: grid_comparison}, we analyse our structured grids together with unstructured grids created using VMTK and Gmsh for both the single branch and bifurcation cases.
In Sec.\ \ref{sec: IGA}, we compare our approach against the state of the art method for the generation of structured grid in the h{\ae}modynamic context.
Sec.\ \ref{sec: topology} discusses the flexibility of the mesh generator in changing the vessel geometry, i.e.\ 
smoothly creating pathological cases (stenoses and aneurysms) without affecting the vessel mesh topology. \textcolor{black}{Finally, in Sec.\ \ref{sec: complex_tree}, we apply our method to a complex coronary tree to show the robustness of the proposed approach.}
\item Fluid results (Sec.\ \ref{sec: fluid_results}). In Sec.\ \ref{sec: fluidynamics}, the fluid dynamics and the lumped parameter models used to impose outflow physiological boundary conditions are described, while Sec.\ \ref{sec: numerical_results} showcases physiological numerical results for a patient-specific coronary geometry along with a mesh convergence test. Finally, in Sec.\ \ref{sec: FFR_vs_CFR}, we compute h{\ae}modynamics indices relevant in the prevision of myocardial infarction (MI) of coronary arteries.
\end{itemize}

\subsection{Meshing Results}\label{sec: meshing_results}

\subsubsection{Structured vs Unstructured approaches (VMTK, Gmsh)}\label{sec: grid_comparison}

In Finite Element (FE) simulations, CPU time and obtaining meaningful results are often among the primary bottlenecks. Mesh quality significantly impacts both aspects. However, there are no universally accepted, minimally required standards for a numerical mesh, as mesh quality metrics strongly depend on the type of application and the governing equations \cite{yang2018isoparametric}.
Numerous metrics for assessing numerical grid quality are available in the literature. In h{\ae}modynamics, the scaled Jacobian (SJ) and normalized equiangular skewness (NES) grid quality metrics constitute the default choices, and they can be applied both to structured and unstructured grids \cite{ghaffari2017large,decroocq2023modeling,bovsnjak2023higher}. In what follows, we briefly summarize their significance \cite{yang2018isoparametric}:

\begin{itemize}
\item Scaled Jacobian (SJ). This provides a scalar measure of the deviation from the reference element. It is defined as the ratio between the minimum and the maximum value of the element's Jacobian determinant. For a valid FE element, the values range between $0$ and $1$, with $1$ considered optimal and $0$ indicating a collapsed side. The minimum value should preferably exceed $0.6$. Negative values denote invalid elements. 

\item Normalized equiangular skewness (NES). A measure of the quality of an element's shape. 
The skewness varies between $0$ (good) and $1$ (bad); a value close to $1$ indicates a degenerate cell where nodes are almost coplanar. It is defined as $\displaystyle\operatorname{max}\left[\frac{\theta_{max}- \theta_{e}}{180 - \theta_{e}},\frac{\theta_{e}- \theta_{min}}{\theta_{e}}\right]$ where $\theta_{\text{max}}, \theta_{\text{min}}$ are the largest / smallest angles in the cell/face while $\theta_{e}$ is the optimal angle, which creates an equiangular face/cell ($\theta_e = 60$ for triangles, $\theta_e = 90$ for a quadrilateral). Preferably, the skewness should not exceed $0.5$.
\end{itemize}

To demonstrate the proficiency of the proposed mesh generator, we conduct a comparative grid quality analysis. We base this comparison on aforementioned metrics applied to numerical grids generated by our tool and those produced by the prevalent general purpose meshing software suites VMTK and Gmsh. The unstructured grids are generated starting from a surface input extracted from the external structured surface mesh generated with our approach, thus eliminating potential mesh quality biases stemming from discrepancies in the quality of the input surface mesh. It is worth noting that the quality analysis on the meshes has been performed through the "mesh quality" filter of ParaView \cite{paraviewuserguide} on the mesh VTK \cite{vtkBook} format. 

We assess the performance of the mesh generators in two distinct scenarios: 1) a single arterial branch with a severe stenosis and high curvature and 2) two intersecting branches forming a highly non-planar bifurcation. Note that all meshes are reconstructed from patient-specific invasive coronary angiographies.

\begin{figure}[th!]
\centering
\includegraphics[width = 1.5cm, keepaspectratio]{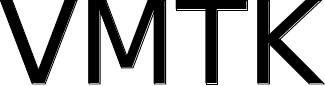}\hspace{4cm}
\includegraphics[width = 1.5cm, keepaspectratio]{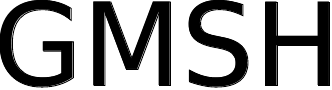}\hspace{4cm}
\includegraphics[width = 1.5cm, keepaspectratio]{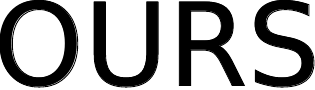}\\\vspace{0.5cm}
\includegraphics[width = 5cm, keepaspectratio]{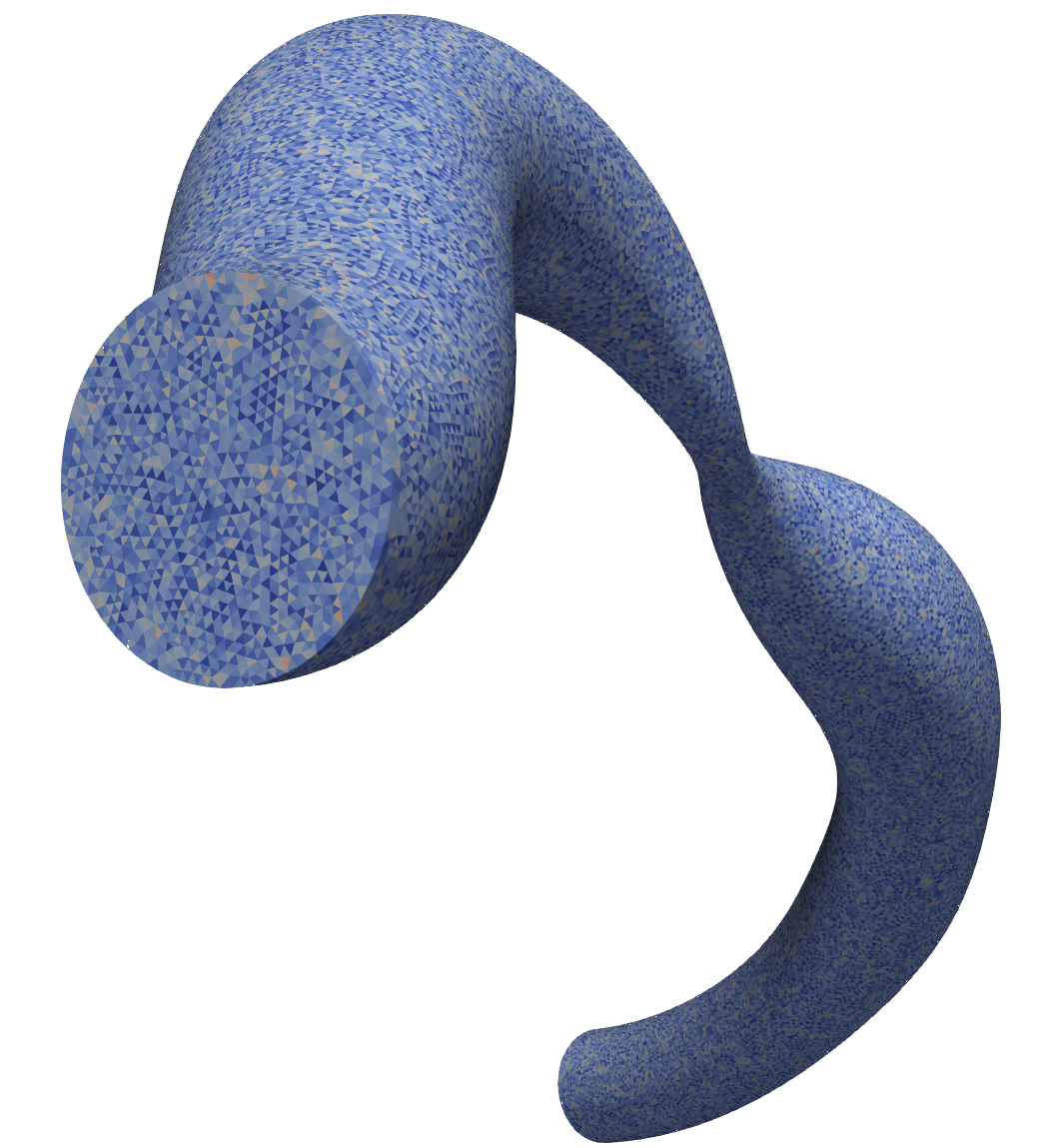}\quad
\includegraphics[width = 5cm, keepaspectratio]{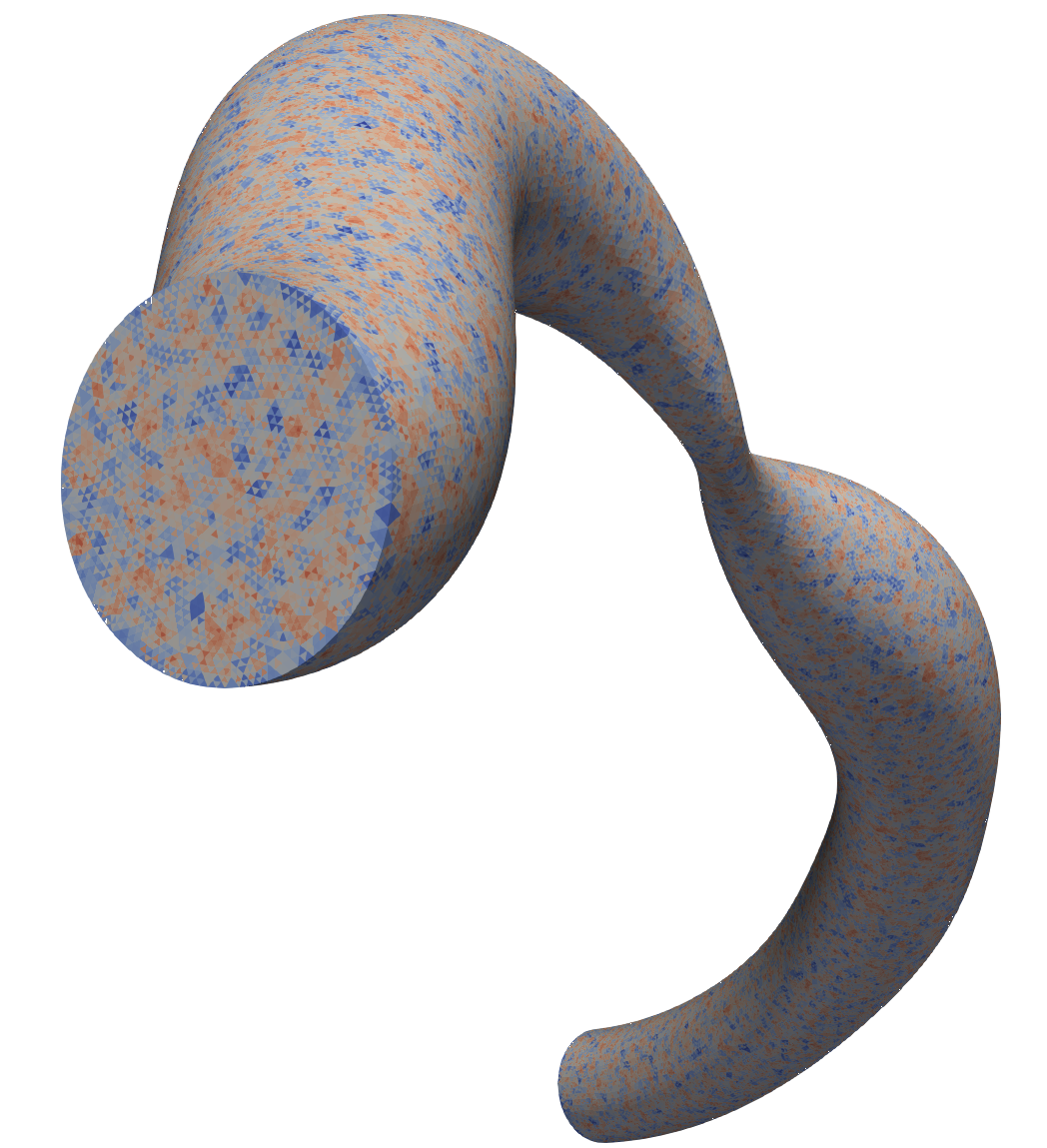}\quad
\includegraphics[width = 5cm, keepaspectratio]{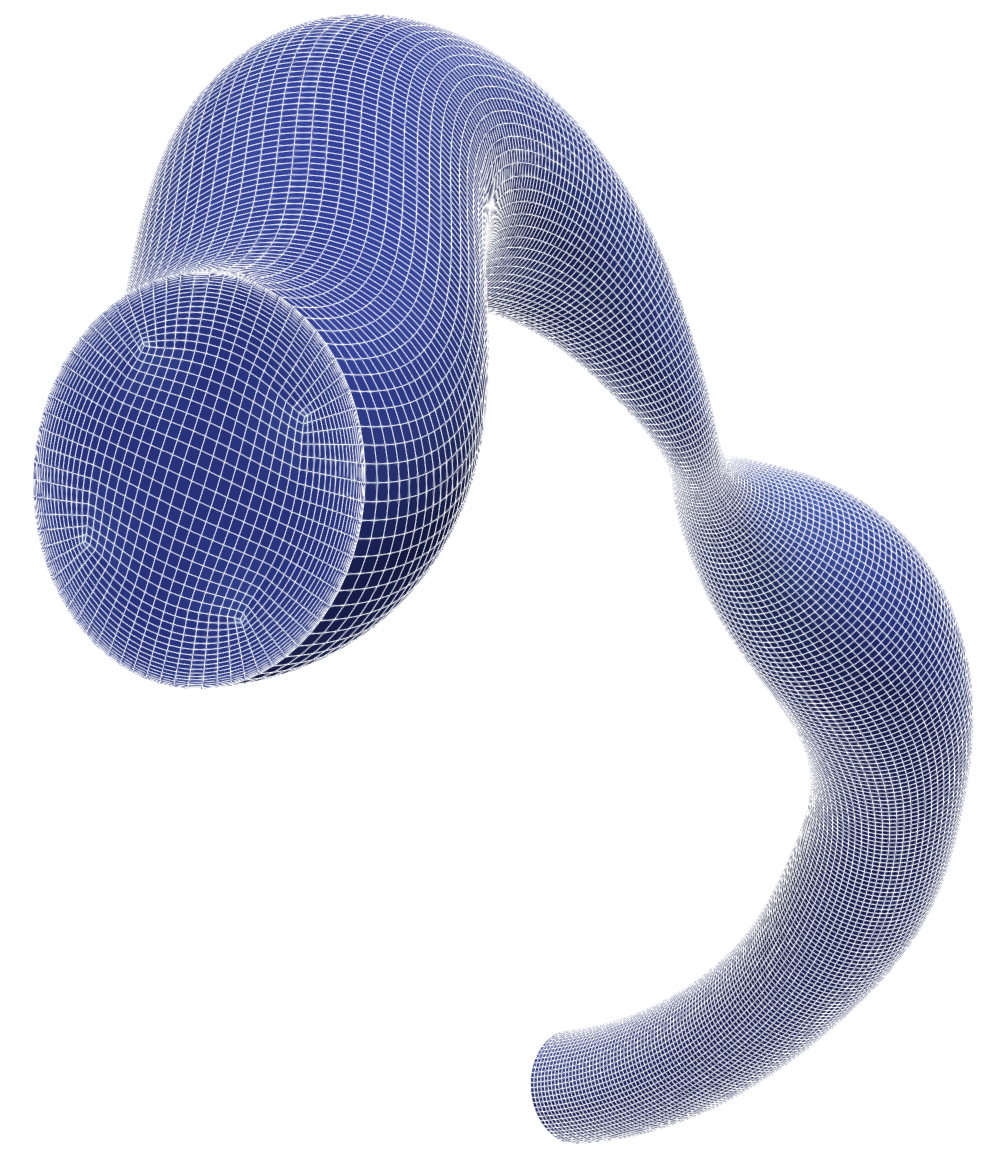} \\\vspace{1cm}
\includegraphics[width = 9cm, keepaspectratio]{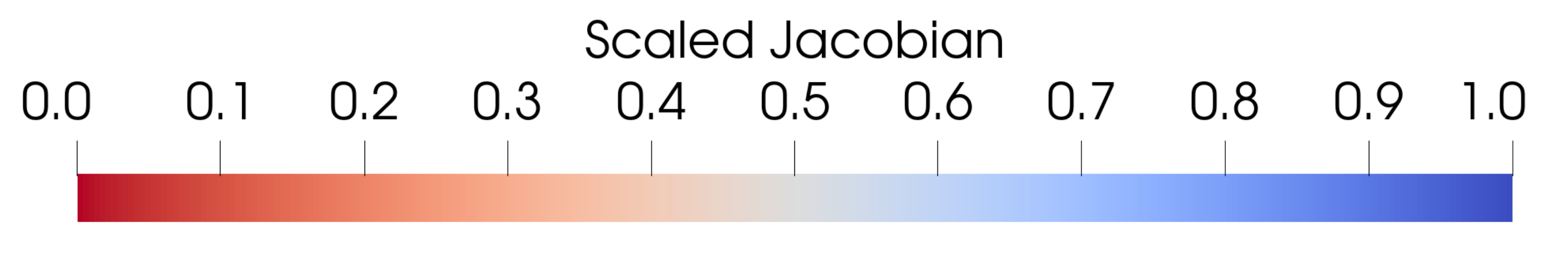}\\ \vspace{1cm}
\includegraphics[width = 1.5cm, keepaspectratio]{images/results/first_section/vmtk-eps-converted-to.pdf}\hspace{4cm}
\includegraphics[width = 1.5cm, keepaspectratio]{images/results/first_section/gmsh_logo-eps-converted-to.pdf}\hspace{4cm}
\includegraphics[width = 1.5cm, keepaspectratio]{images/results/first_section/our_logo-eps-converted-to.pdf}\\\vspace{0.5cm}
\includegraphics[width = 5.2cm, height = 3.5cm]{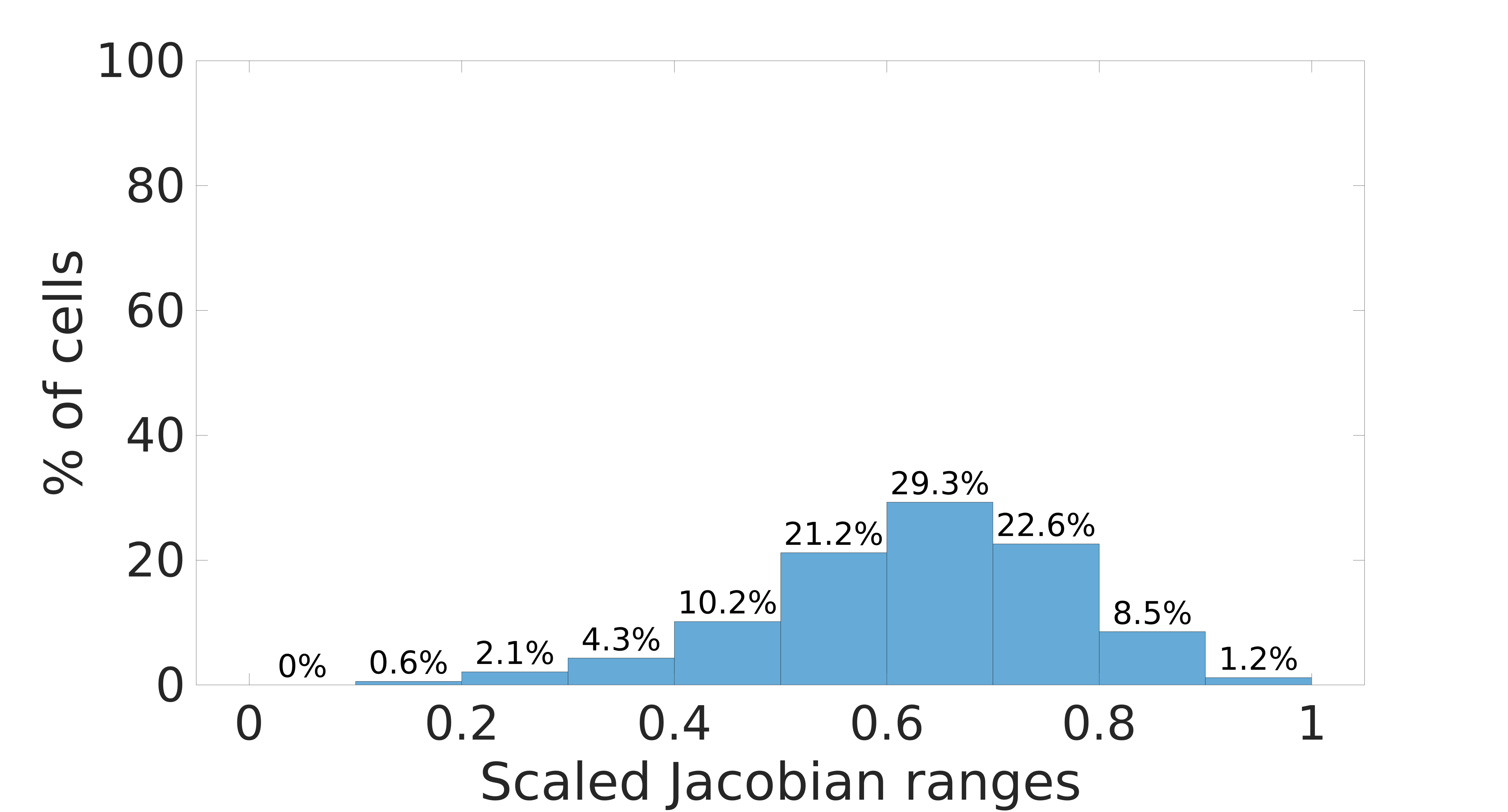}
\includegraphics[width = 5.2cm, height = 3.5cm]{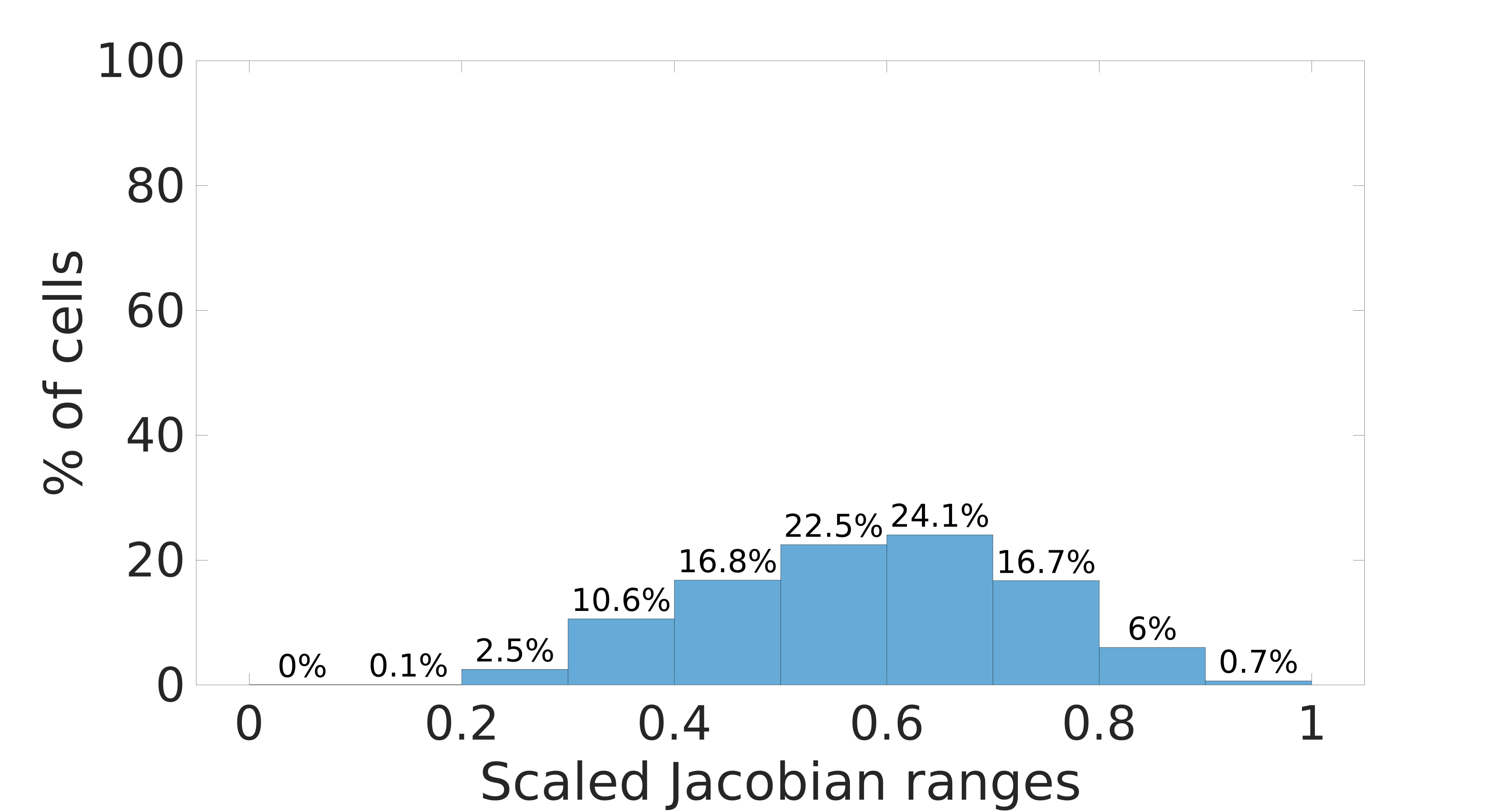}
\includegraphics[width = 5.2cm, height = 3.5cm]{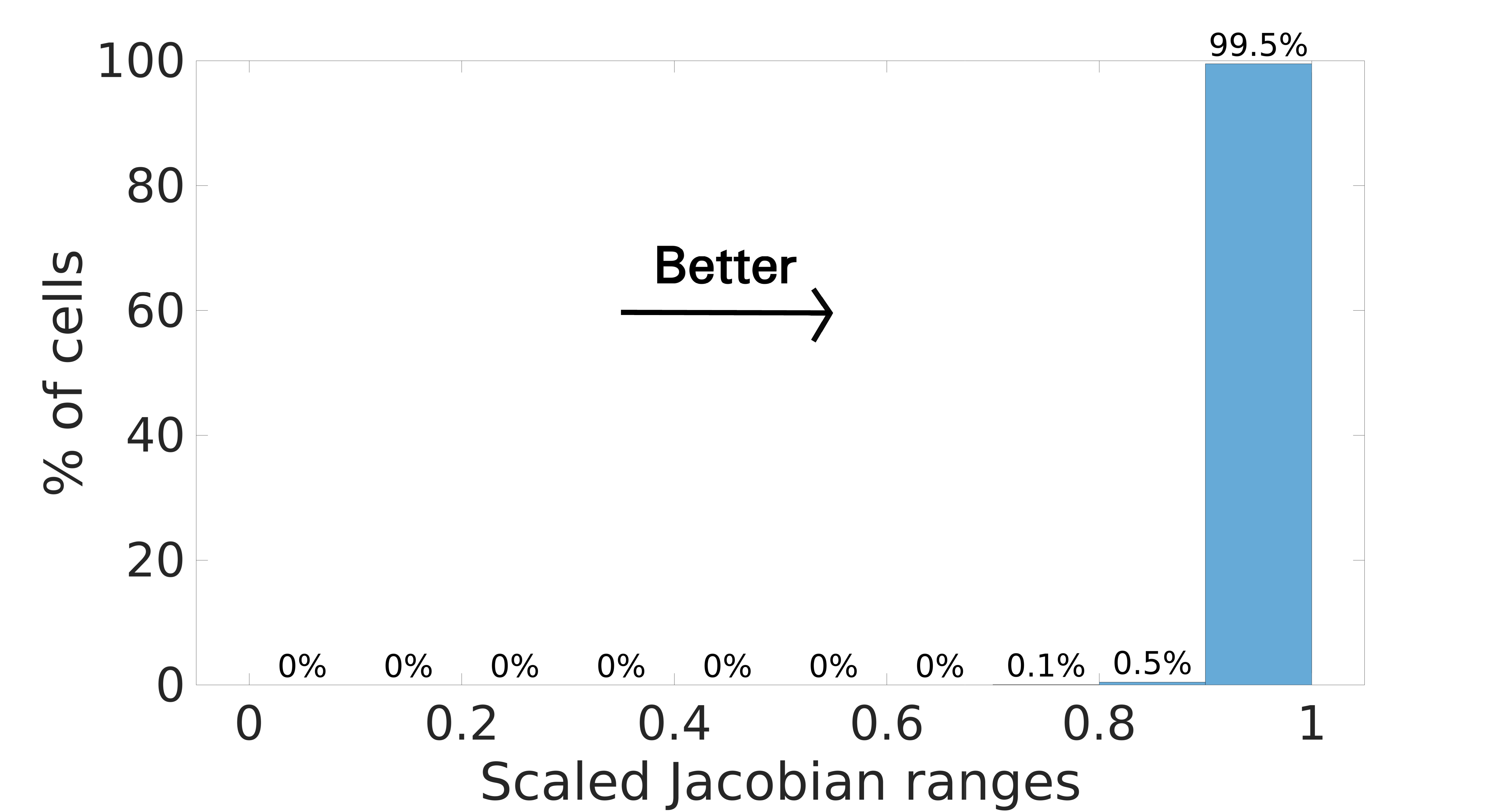}
\includegraphics[width = 5.2cm, height = 3.5cm]{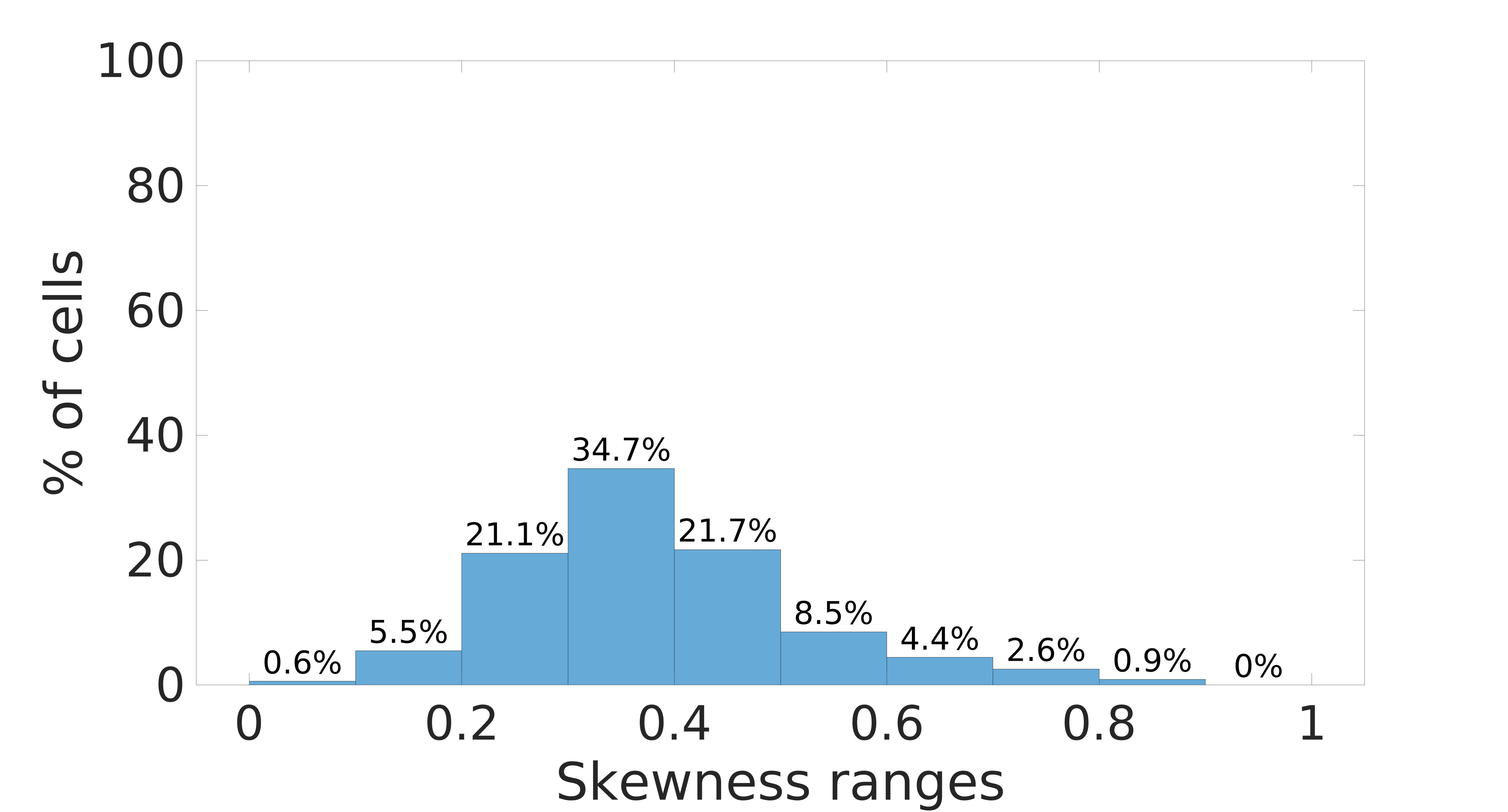}
\includegraphics[width = 5.2cm, height = 3.5cm]{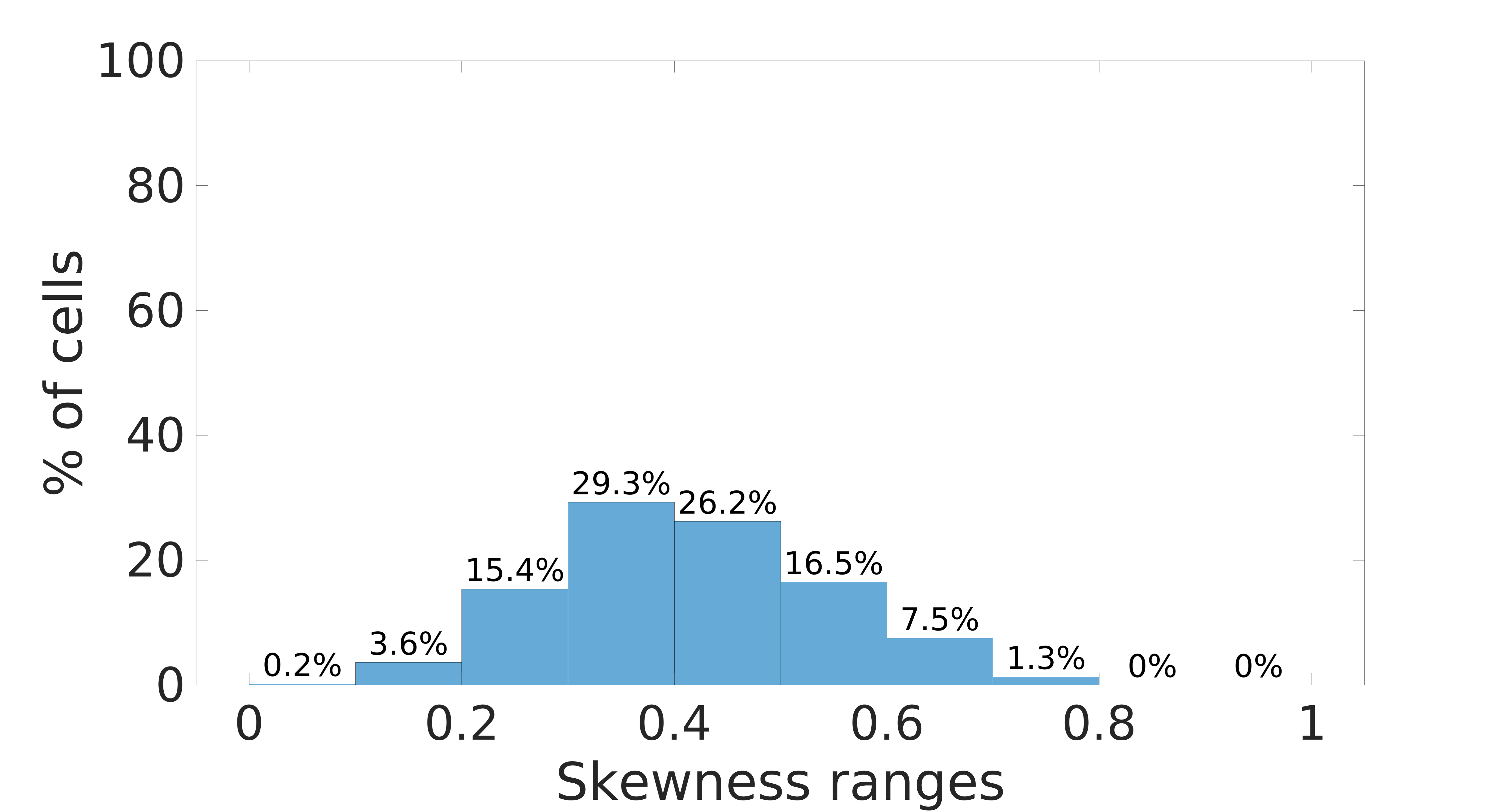}
\includegraphics[width = 5.2cm, height = 3.5cm]{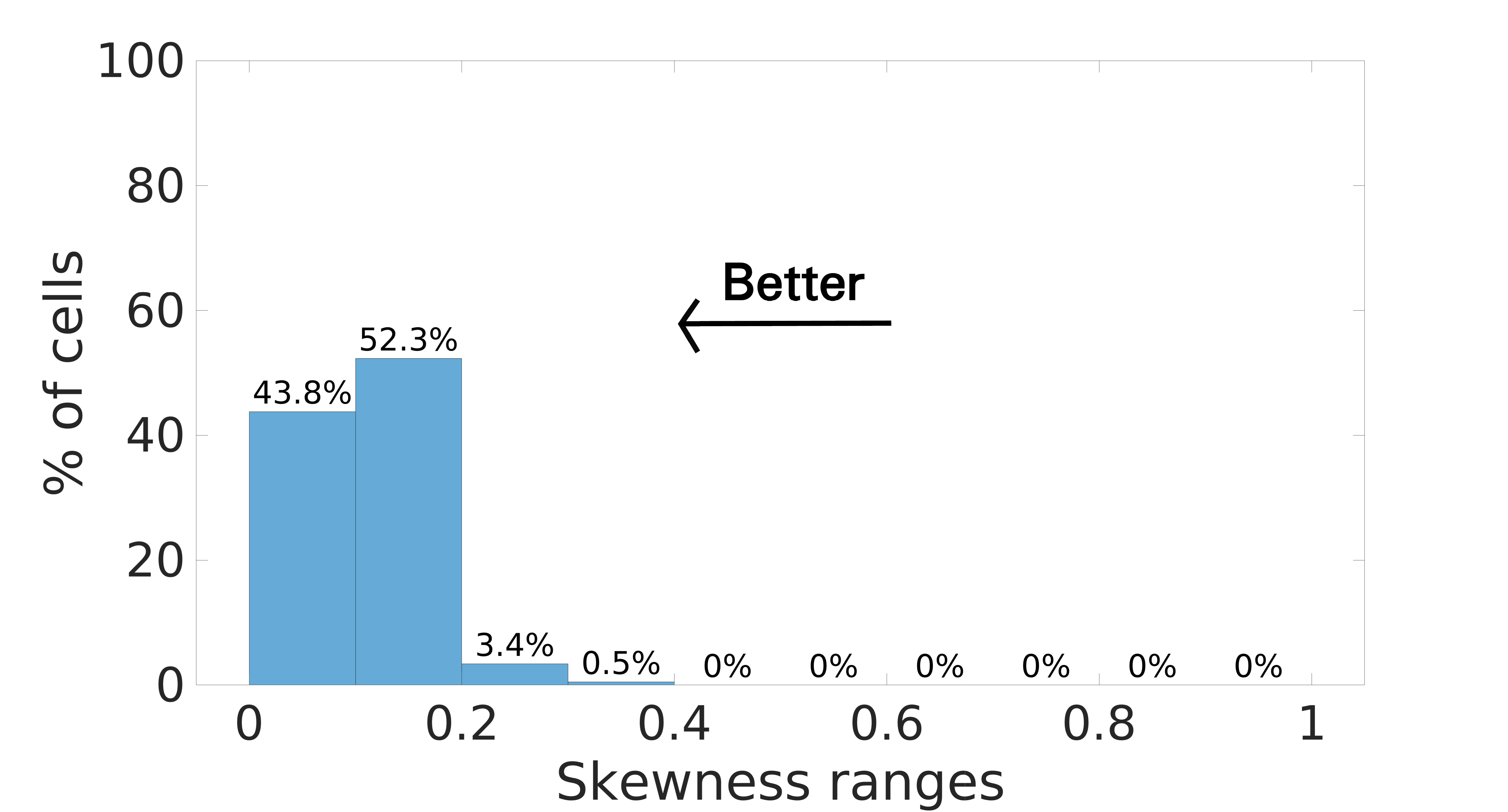}
\caption{Comparison between structured and unstructured grids for the single branch case. The first row depicts the values of the elementwise SJ measure. Our grid is highlighted in white.
The second row shows histograms illustrating the percentage of cells that fall within specified ranges of the SJ. The third row shows the percentage of cells within a certain range of the NES measure.}
\label{Fig: comparison_single_branch}
\end{figure}

Figure \ref{Fig: comparison_single_branch}, top-row, depicts a $3$D visualization of the scaled Jacobian per element for the single branch geometry; from the left to the right column, VMTK, Gmsh and our mesh generator, respectively, are reported.
The second row presents histograms illustrating the percentage of cells falling within specified ranges of the SJ. Meanwhile, the last row displays the percentage of cells within various ranges of the NES. Starting from the SJ measure, the structured mesh contains no elements that fall below the value of $0.78$ and $99.5$ \% of the cells have a value between $0.9$ and 1.
VMTK, which is specifically designed for meshing vessel geometries, shows better results than Gmsh.
The same scenario repeats itself for the NES measure. More than $80 \%$ of all structured cells exhibit minuscule skewness which is indicative of excellent grid quality. Meanwhile, both VMTK and Gmsh create meshes wherein more than $30 \%$ of the cells exceed the preferred skewness limit of $0.5$. 
Moreover, it is worth noting that our results are comparable to the ones obtained in \cite{ghaffari2017large, decroocq2023modeling}.  
In Tab.\ \ref{tab: unstructured_grid_comparison_single_branch}, the minimum, average and maximum SJ ($\text{SJ}_{\text{min}}$, $\text{SJ}_{\text{mean}}$ and $\text{SJ}_{\text{max}}$, respectively) and the same for the NES ($\text{EJ}_{\text{min}}$, $\text{EJ}_{\text{mean}}$ and $\text{EJ}_{\text{max}}$) are reported for all grids.

\begin{figure}[th!]
\centering
\includegraphics[width = 1.5cm, keepaspectratio]{images/results/first_section/vmtk-eps-converted-to.pdf}\hspace{4cm}
\includegraphics[width = 1.5cm, keepaspectratio]{images/results/first_section/gmsh_logo-eps-converted-to.pdf}\hspace{4cm}
\includegraphics[width = 1.5cm, keepaspectratio]{images/results/first_section/our_logo-eps-converted-to.pdf}\\
\includegraphics[width = 5.3cm, keepaspectratio]{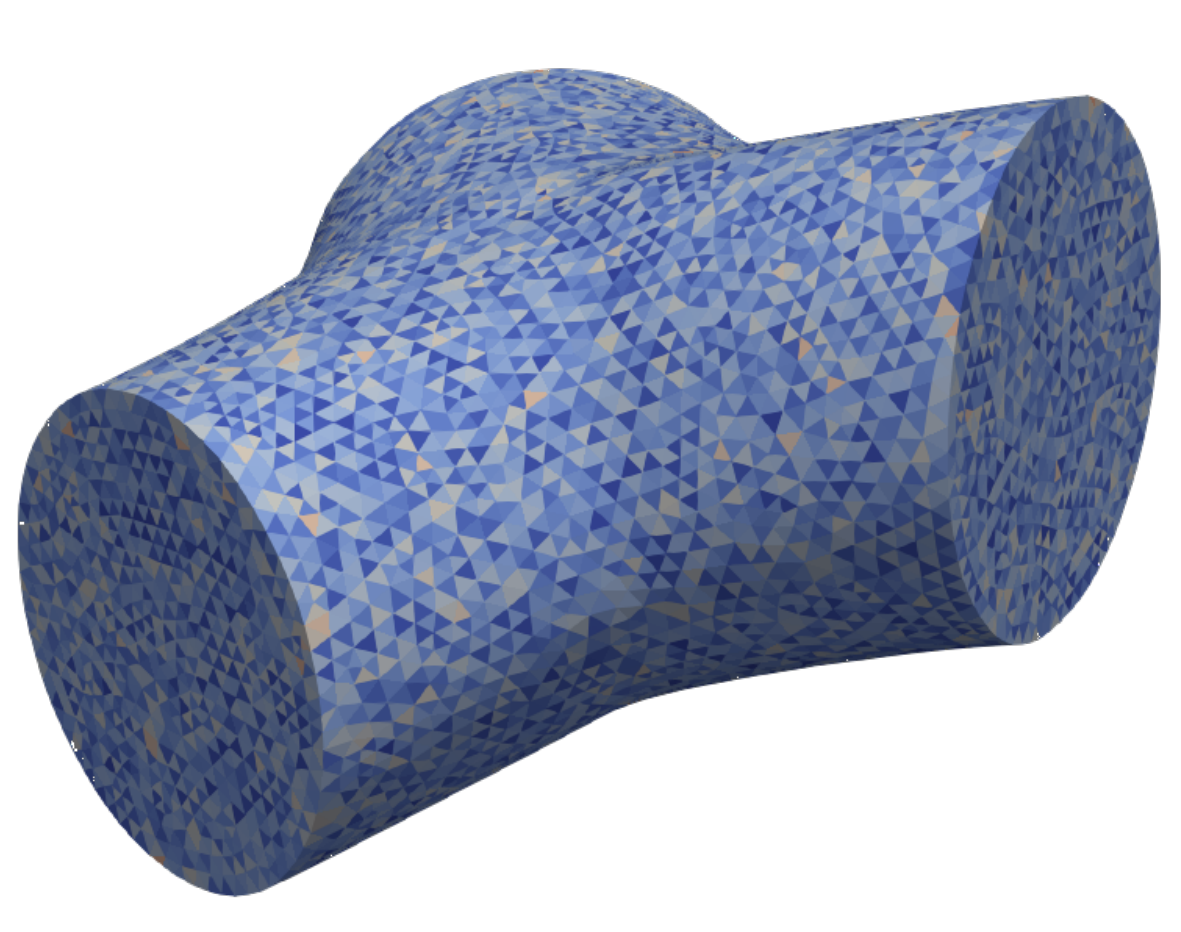}
\includegraphics[width = 5.3cm, keepaspectratio]{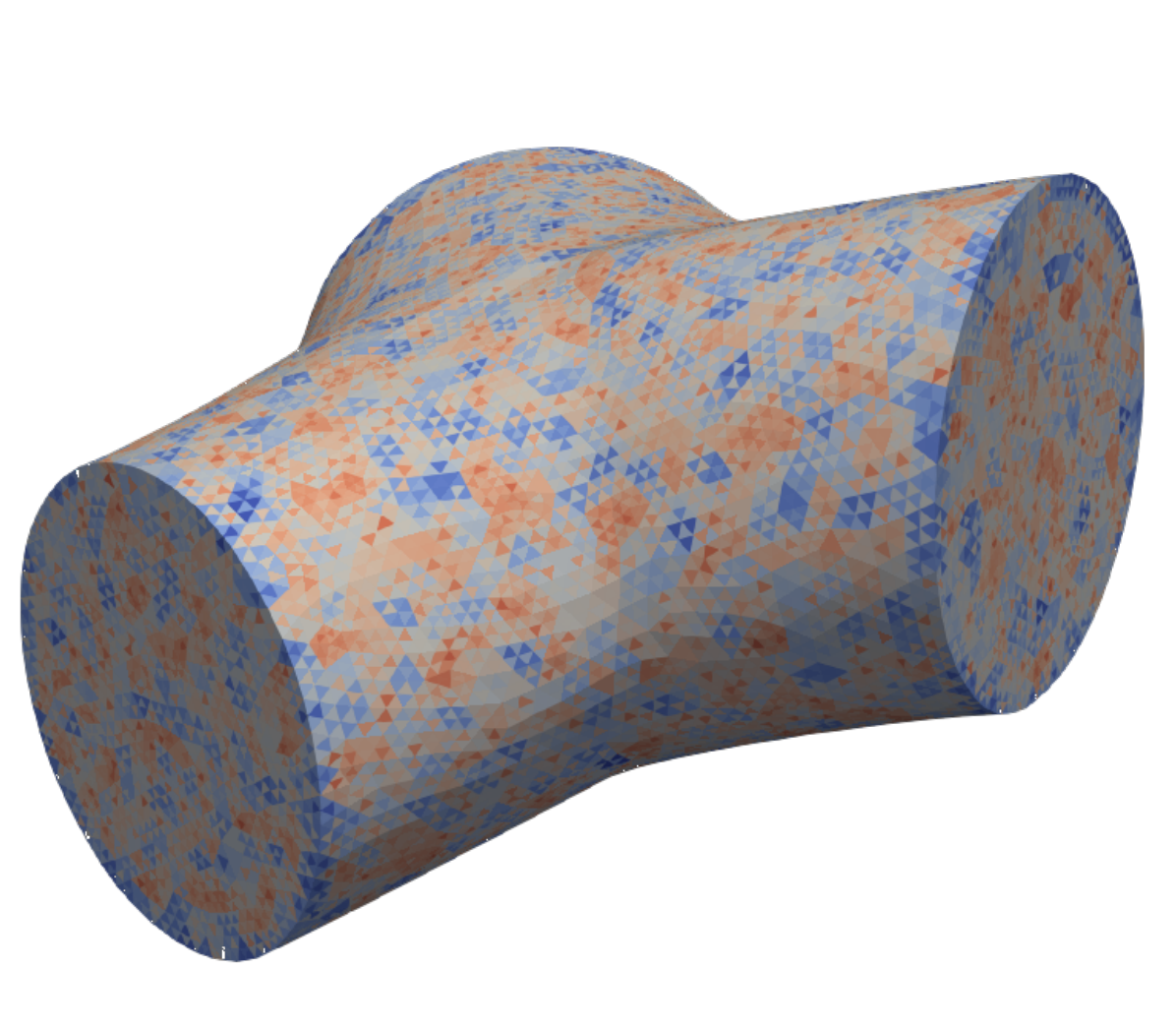}
\includegraphics[width = 5.6cm, keepaspectratio]{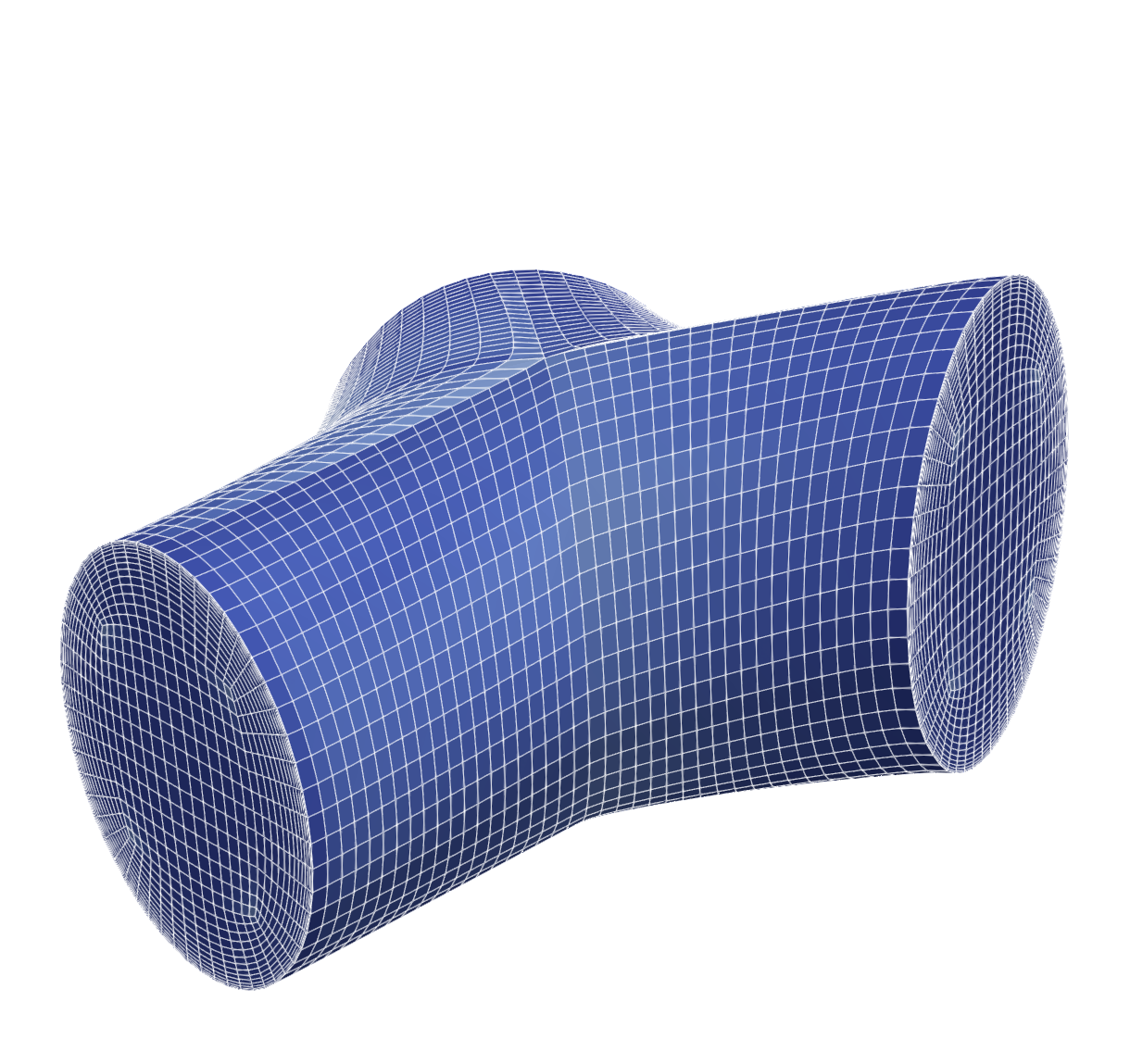}\vspace{1cm}
\includegraphics[width = 8cm, keepaspectratio]{images/results/first_section/scale_jacobian.pdf}\\ \vspace{1cm}
\includegraphics[width = 1.5cm, keepaspectratio]{images/results/first_section/vmtk-eps-converted-to.pdf}\hspace{4cm}
\includegraphics[width = 1.5cm, keepaspectratio]{images/results/first_section/gmsh_logo-eps-converted-to.pdf}\hspace{4cm}
\includegraphics[width = 1.5cm, keepaspectratio]{images/results/first_section/our_logo-eps-converted-to.pdf}\\\vspace{0.5cm}
\includegraphics[width = 5.2cm, height = 3.5cm]{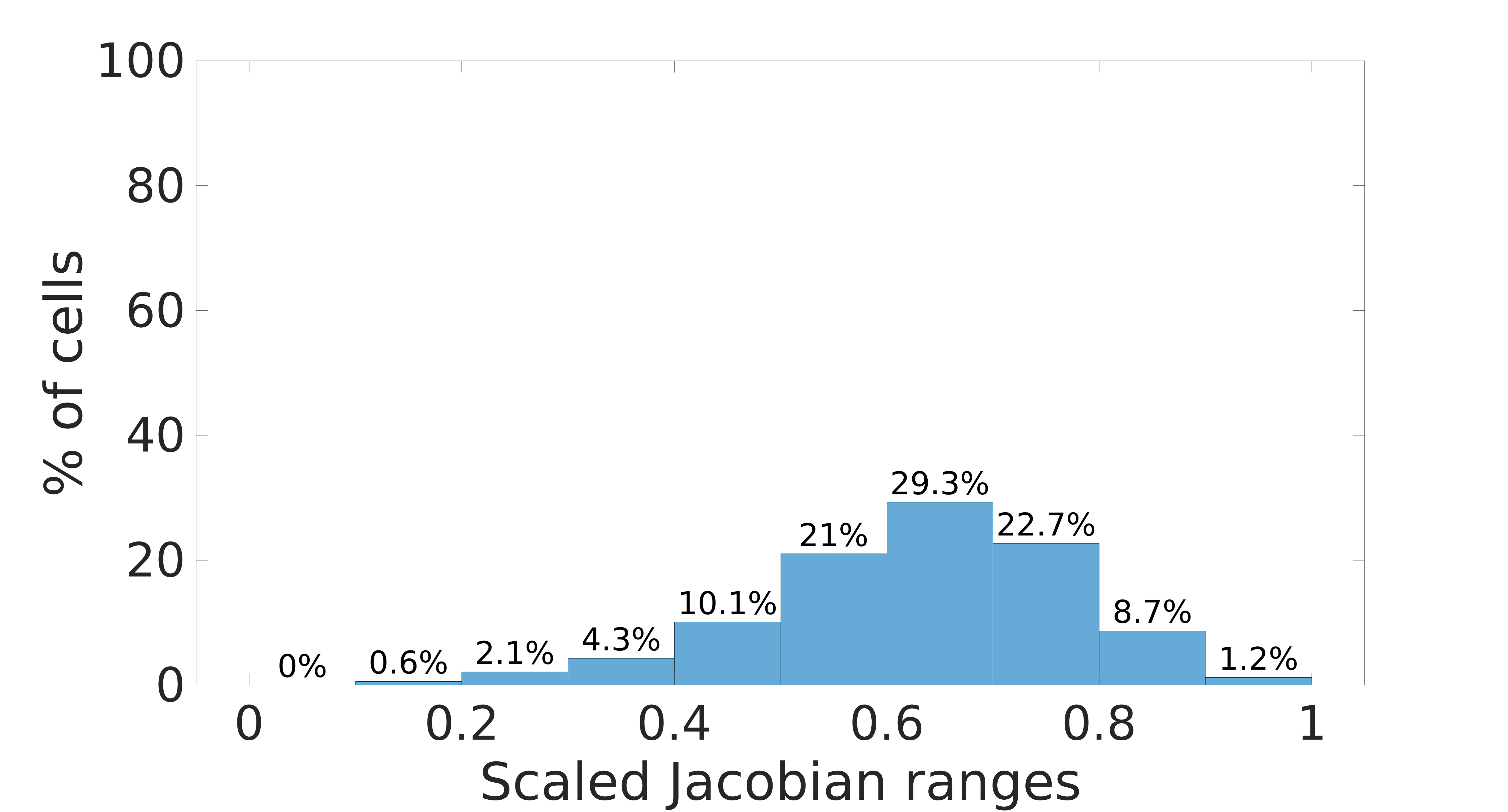}
\includegraphics[width = 5.2cm, height = 3.5cm]{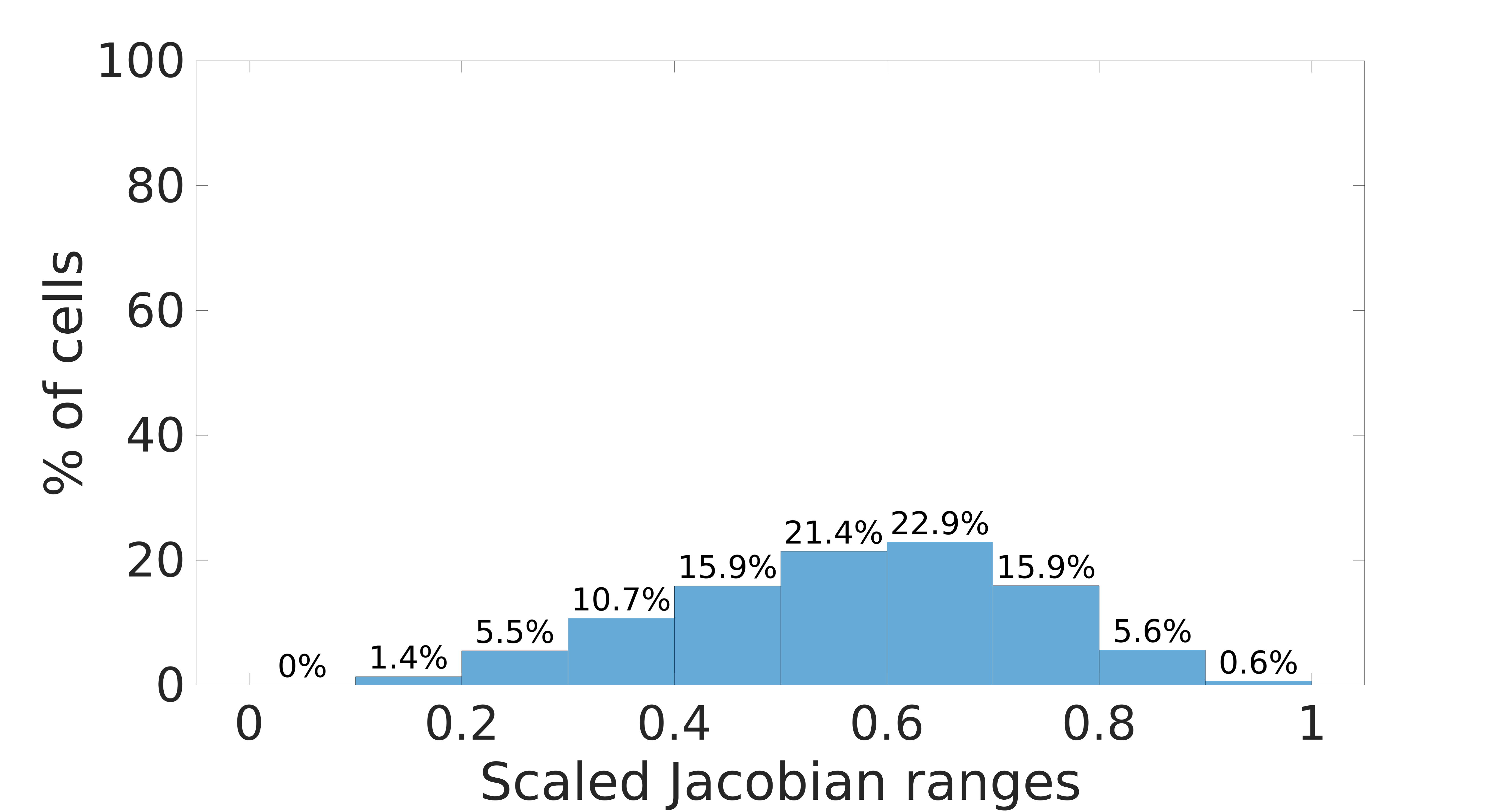}
\includegraphics[width = 5.2cm, height = 3.5cm]{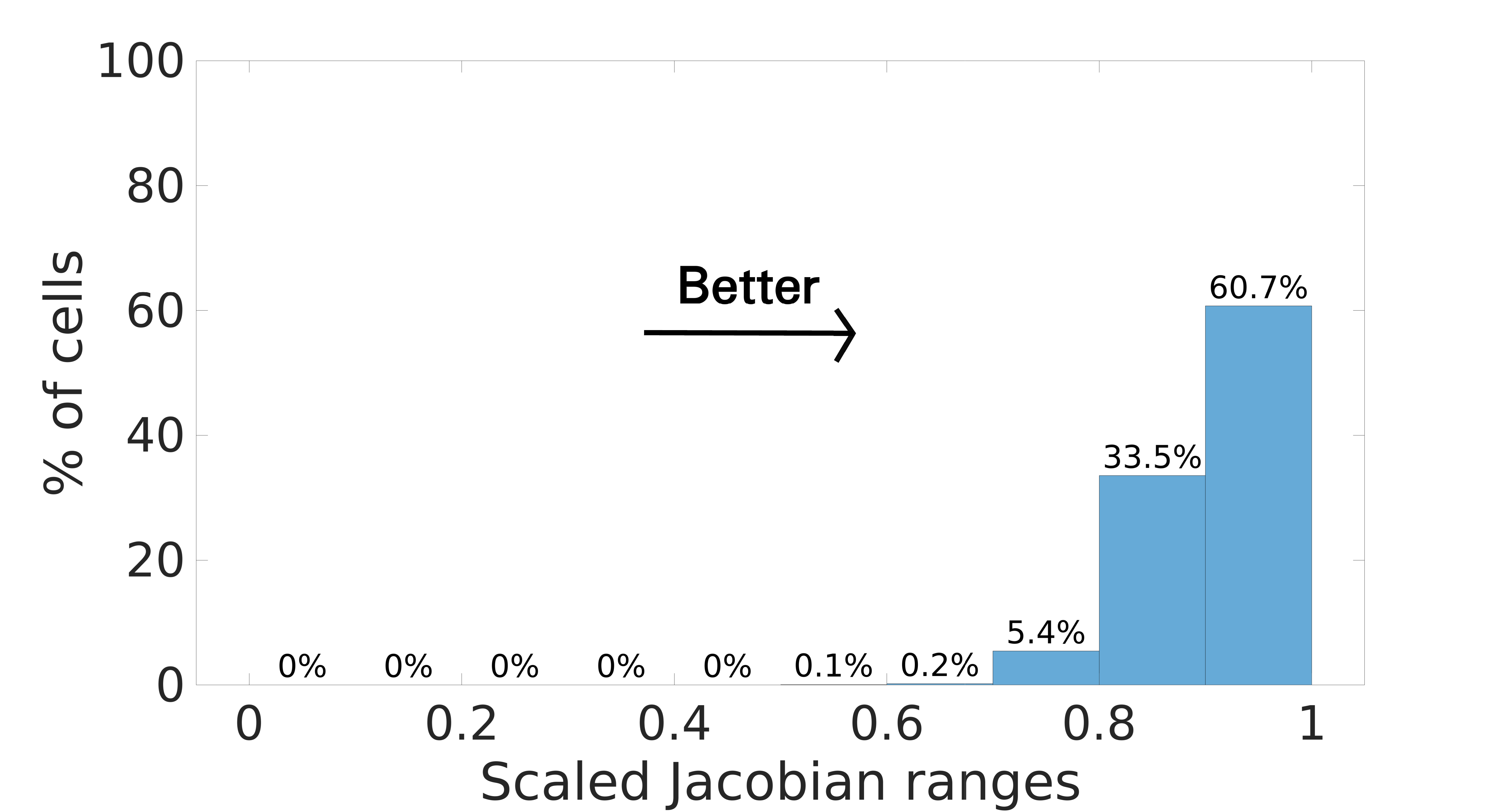}\\ \vspace{1cm}
\includegraphics[width = 5.2cm, height = 3.5cm]{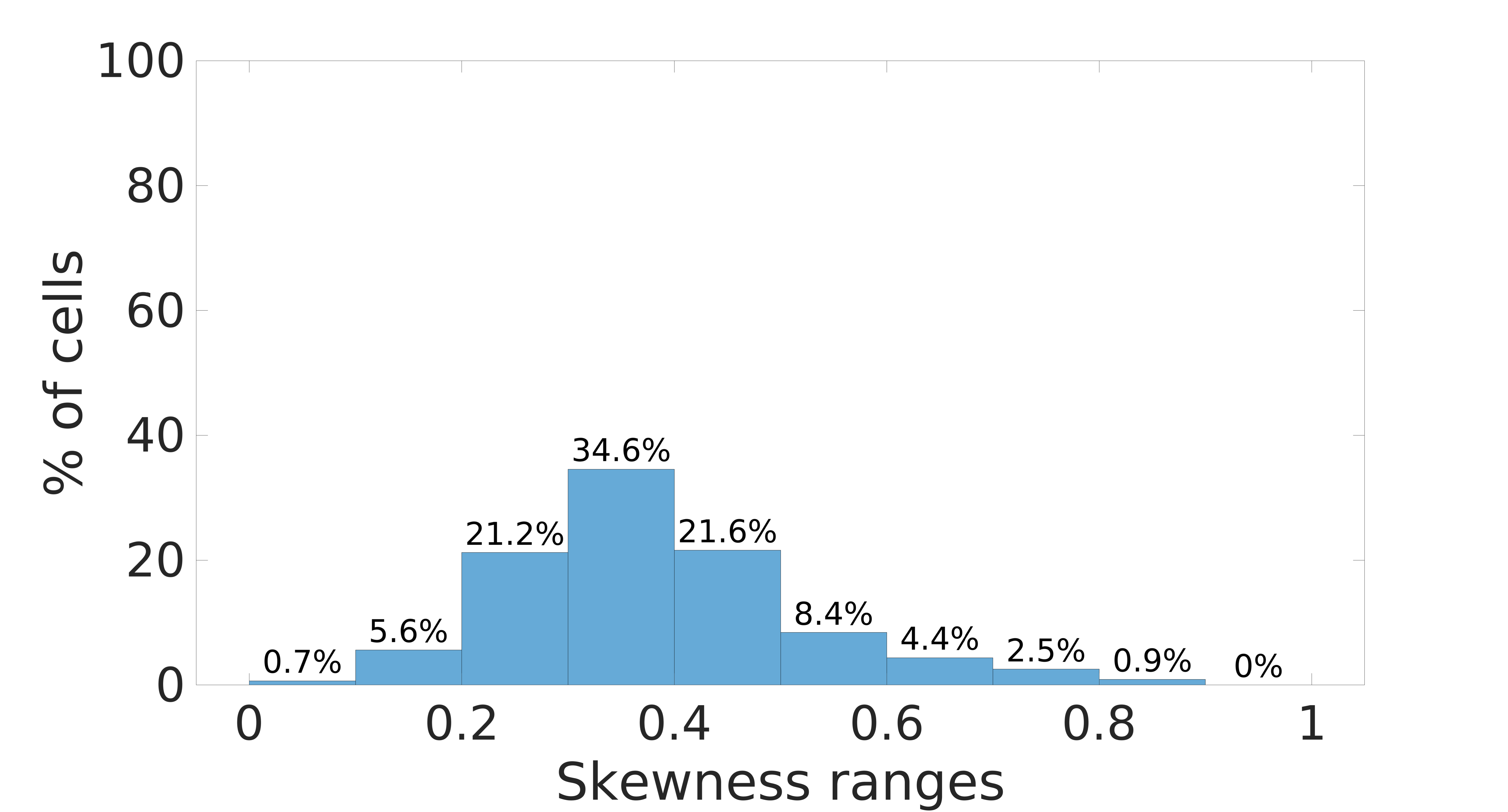}
\includegraphics[width = 5.2cm, height = 3.5cm]{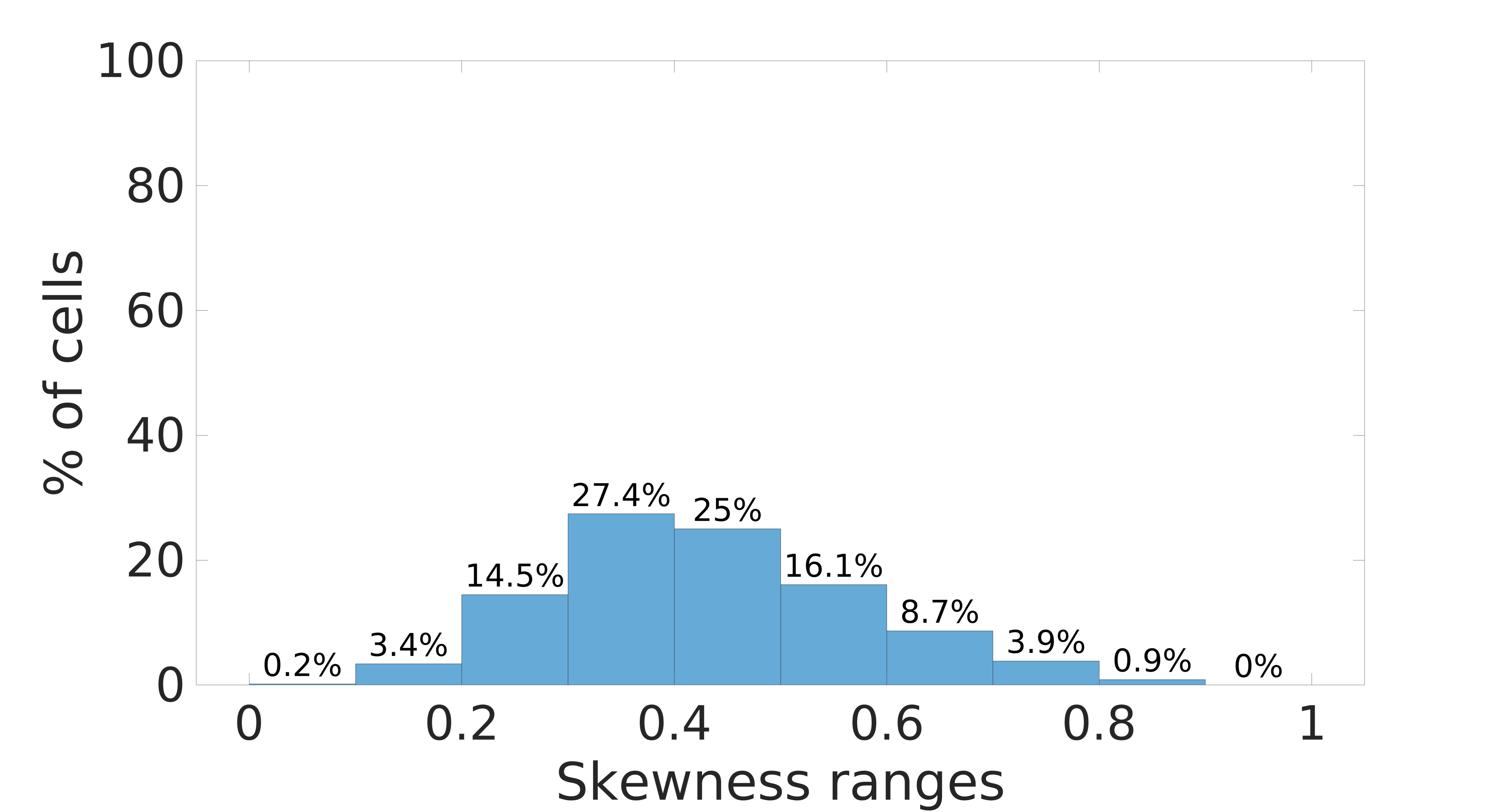}
\includegraphics[width = 5.2cm, height = 3.5cm]{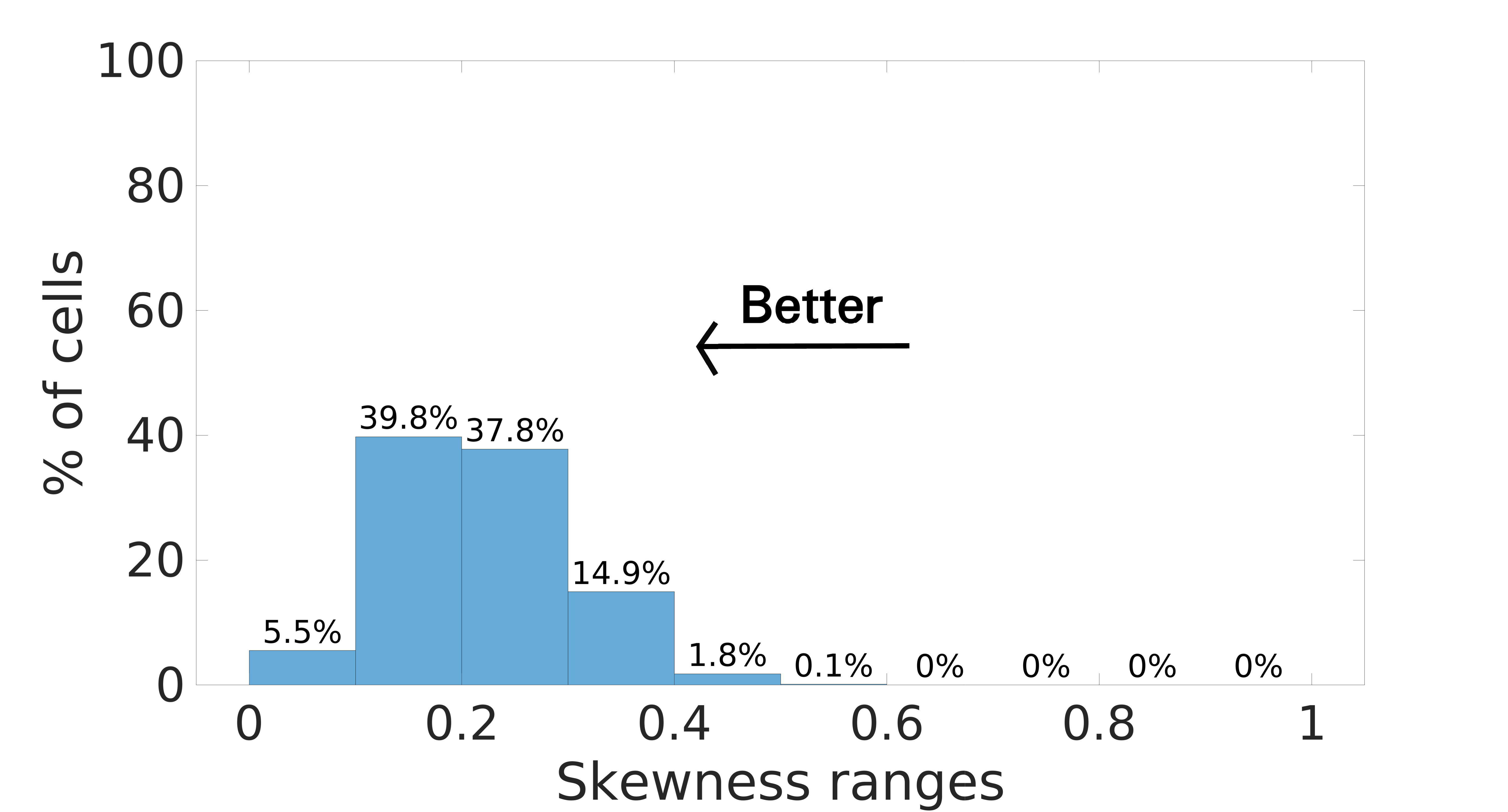}
\caption{Comparison between structured and unstructured grids for the bifurcation case. The first row depicts the values of the SJ measure.
Second row shows the percentage of cells within a certain range of the SJ measure. The third row reports the analogous results for the NES measure. }
\label{Fig: comparison_bif}
\end{figure}

\begin{table}[h!]
\caption{Comparison between our approach, VMTK and Gmsh for the single branch test case with the number of vertices in parentheses. SJ stands for \text{scaled Jacobian} and NES for \text{normalized equiangular skewness}}
\centering
\begin{tabular}{|c|c|c|c|c|c|c|}
\hline
 & \multicolumn{6}{|c|}{Single branch}  \\
\hline
 & $\text{SJ}_{\text{min}}$ & $\text{SJ}_{\text{mean}}$ & $\text{SJ}_{\text{max}}$ & $\text{NES}_{\text{min}}$ & $\text{NES}_{\text{mean}}$ & $\text{NES}_{\text{max}}$ \\
\hline 
\text{VMTK} ($\sim{550k}$) & 0.092 & 0.627 & 0.998 & 0.002 & 0.382 & 0.892  \\
\text{Gmsh} ($\sim{520k}$)& 0.073 & 0.582 & 0.998 & 0.008 & 0.413 & 0.855  \\
\text{ours} ($\sim{570k}$) & 0.785 & 0.979 & 0.999 & 0.008 & 0.109 & 0.402  \\
\hline
\end{tabular}
\label{tab: unstructured_grid_comparison_single_branch}
\end{table}

\begin{table}[h!]
\caption{Comparison between our approach, VMTK and Gmsh for the bifurcation test case with the associated number of nodes in in parentheses. SJ stands for \text{scaled Jacobian} and NES for \text{normalized equiangular skewness}}
\centering
\begin{tabular}{|c|c|c|c|c|c|c|}
\hline
 & \multicolumn{6}{|c|}{Bifurcation}  \\
\hline
 & $\text{SJ}_{\text{min}}$ & $\text{SJ}_{\text{mean}}$ & $\text{SJ}_{\text{max}}$ & $\text{NES}_{\text{min}}$ & $\text{NES}_{\text{mean}}$ & $\text{NES}_{\text{max}}$ \\
\hline 
\text{VMTK} ($\sim{64k}$) & 0.078 & 0.627 & 0.995 & 0.008 & 0.381 & 0.883  \\
\text{Gmsh} ($\sim{68k}$)& 0.08 & 0.564 & 0.996 & 0.012 & 0.431 & 0.890  \\
\text{ours} ($\sim{64k}$) & 0.488 & 0.907 & 0.999 & 0.017& 0.217& 0.634  \\
\hline
\end{tabular}
\label{tab: unstructured_grid_comparison_bifurcation}
\end{table}

\begin{figure}[h!]
\centering
\includegraphics[width = 7cm, keepaspectratio]{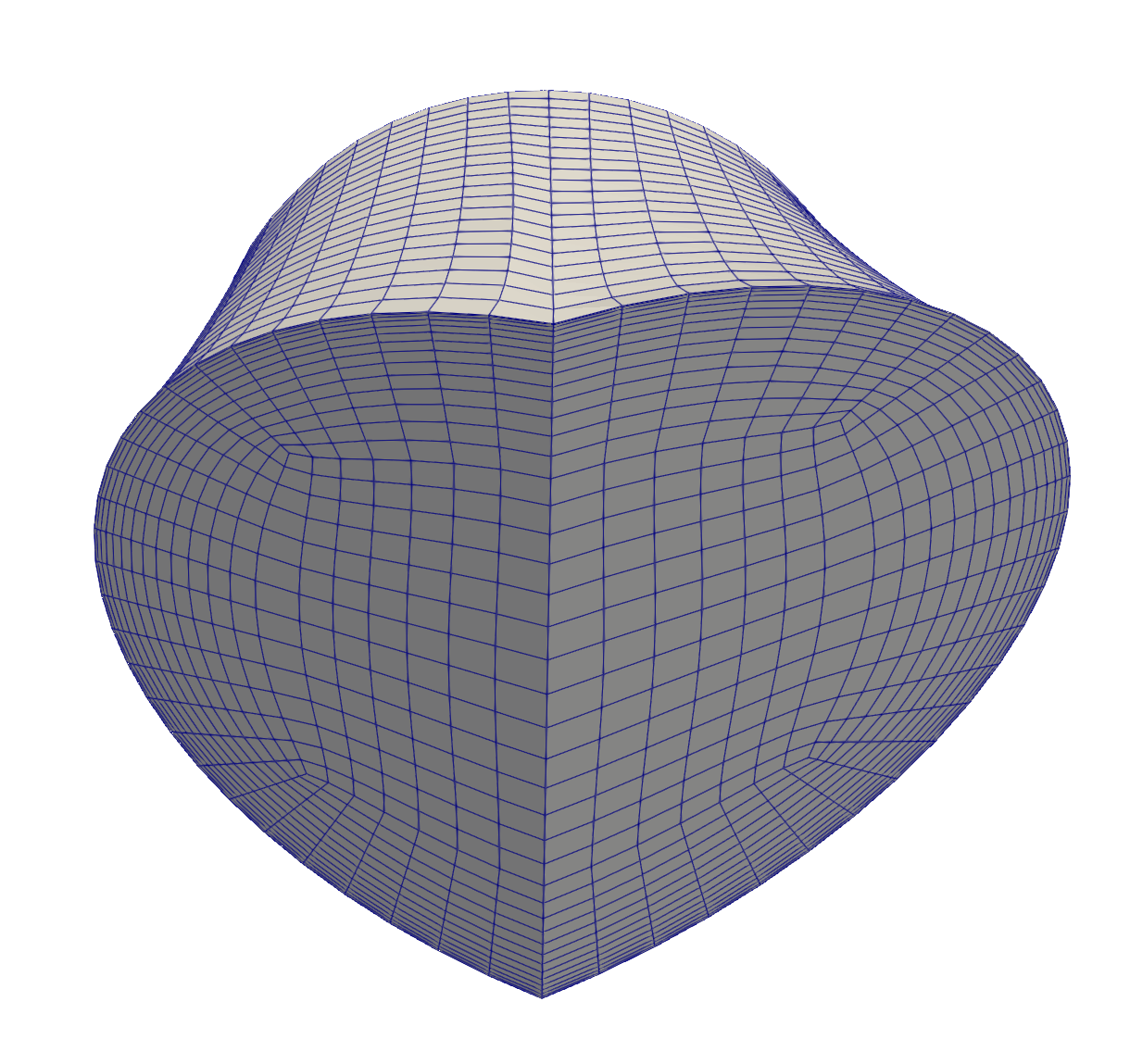}\quad
\includegraphics[width = 7cm, keepaspectratio]{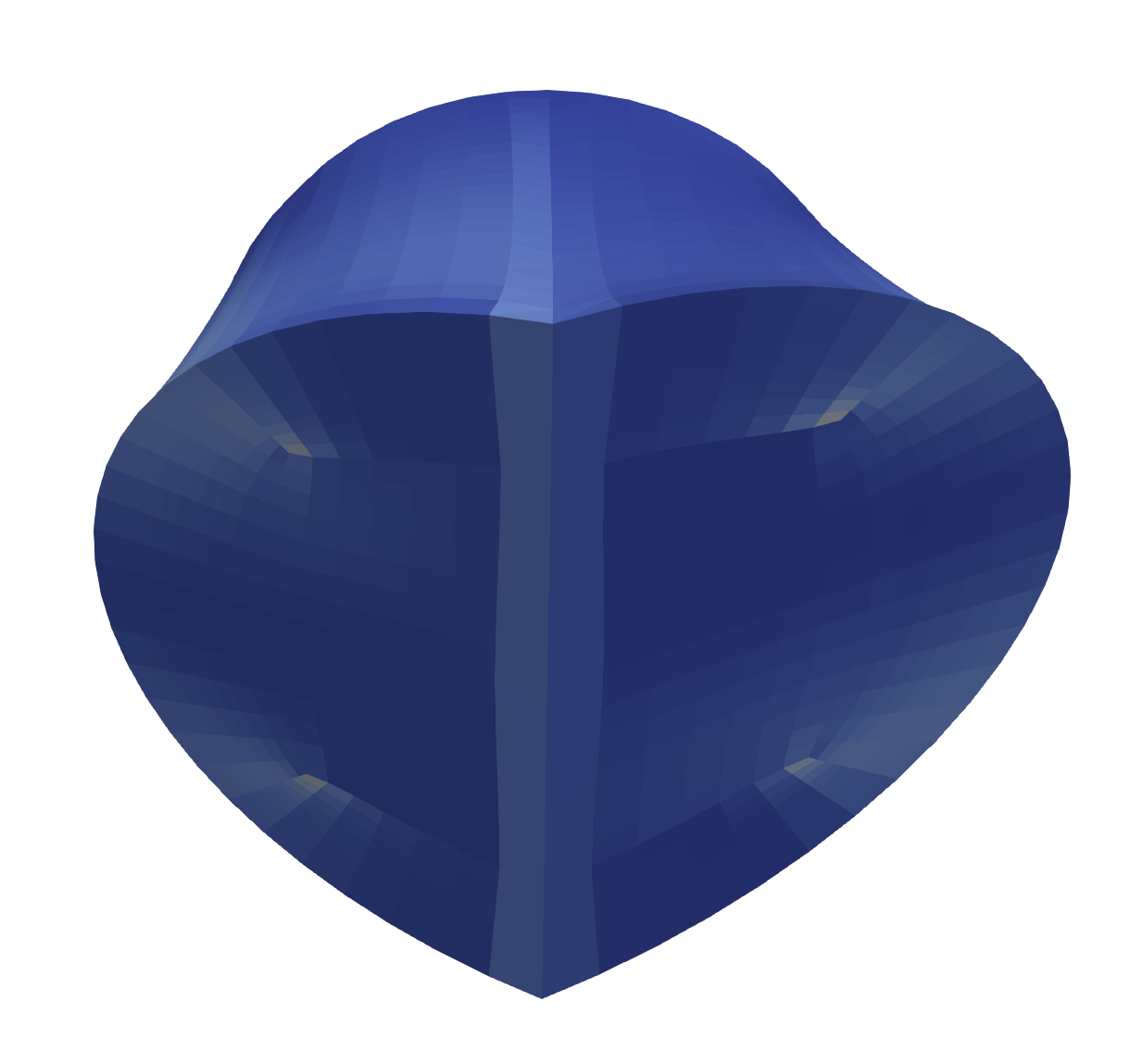}\quad
\includegraphics[width = 8cm, keepaspectratio]{images/results/first_section/scale_jacobian.pdf}
\caption{Left: Numerical mesh inside the bifurcation. Right: Scaled Jacobian.}
\label{Fig: comparison_bif}
\end{figure}

In Fig. \ref{Fig: comparison_bif}, we report the comparison against VMTK and Gmsh for the non-planar bifurcation test case.
The same trend is also visible in this test case, where the structured grid shows better results for both the NES and SJ measures.
The presence of a bifurcation leads to a mesh quality deterioration in the structured approach. Conversely, a comparison to the single branch results demonstrates that the unstructured approaches are only marginally affected. Concerning the structured setting, Fig.\ \ref{Fig: comparison_bif} shows that the mesh quality is bottlenecked by the central patch's four lateral external vertices whose adjacent cells display the lowest values of the scaled Jacobian. This is not surprising given the inherent distortion of the reference grid (see Fig. \ref{Fig: harmonic_map}) in the vicinity of the extraordinary vertices. Additionally, the butterfly structure inherently introduces further distortion to the cells. We emphasize that, despite this limitation, the structured approach still significantly outperforms both unstructured approaches.

As previously done for the single branch test case, in Tab.\ \ref{tab: unstructured_grid_comparison_bifurcation}, we report the statistics of the metrics for the bifurcation case. 

In summary, in both cases (single branch and bifurcation), the structured mesh demonstrates superiority over the unstructured approach in both the NES and SJ measures, despite equal input surface information. 
Additionally, it is worth noting that our approach enables easy control over the total number of grid nodes.

Not surprisingly, incorporating the boundary layers (see Sec.\ \ref{sec: BL}) inherently introduces more cell distortion resulting in deterioration of the SJ measure. Here, it is important to find a good trade-off between the desired grid density by the boundary and undesirable side effects, such as increased grid skewness. 

\subsubsection{Comparison against state of the art for structured meshes}\label{sec: IGA} 

\begin{figure}[t]
\centering
\includegraphics[width = 6.5cm, keepaspectratio]{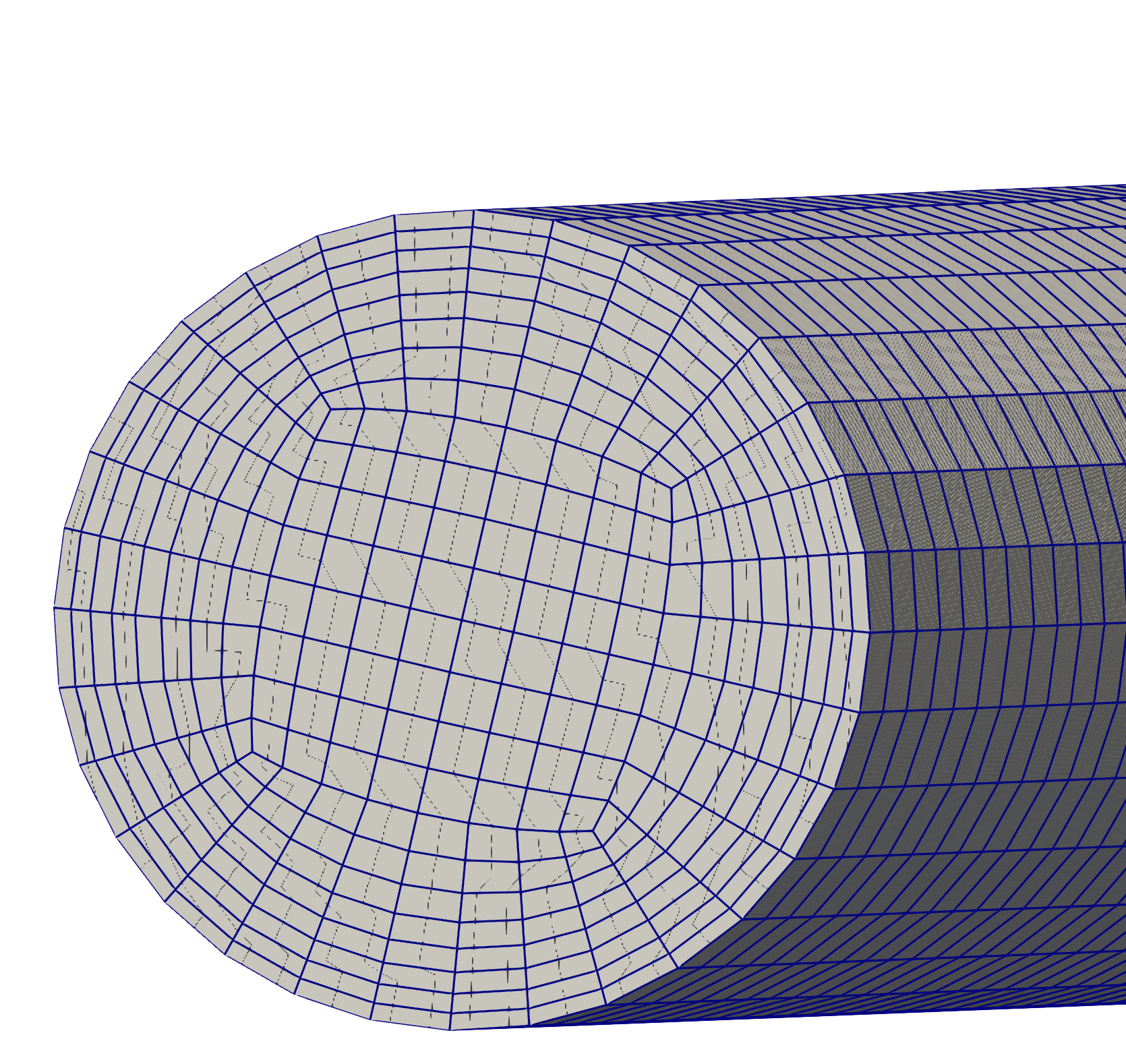}\quad
\includegraphics[width = 6.5cm, keepaspectratio]{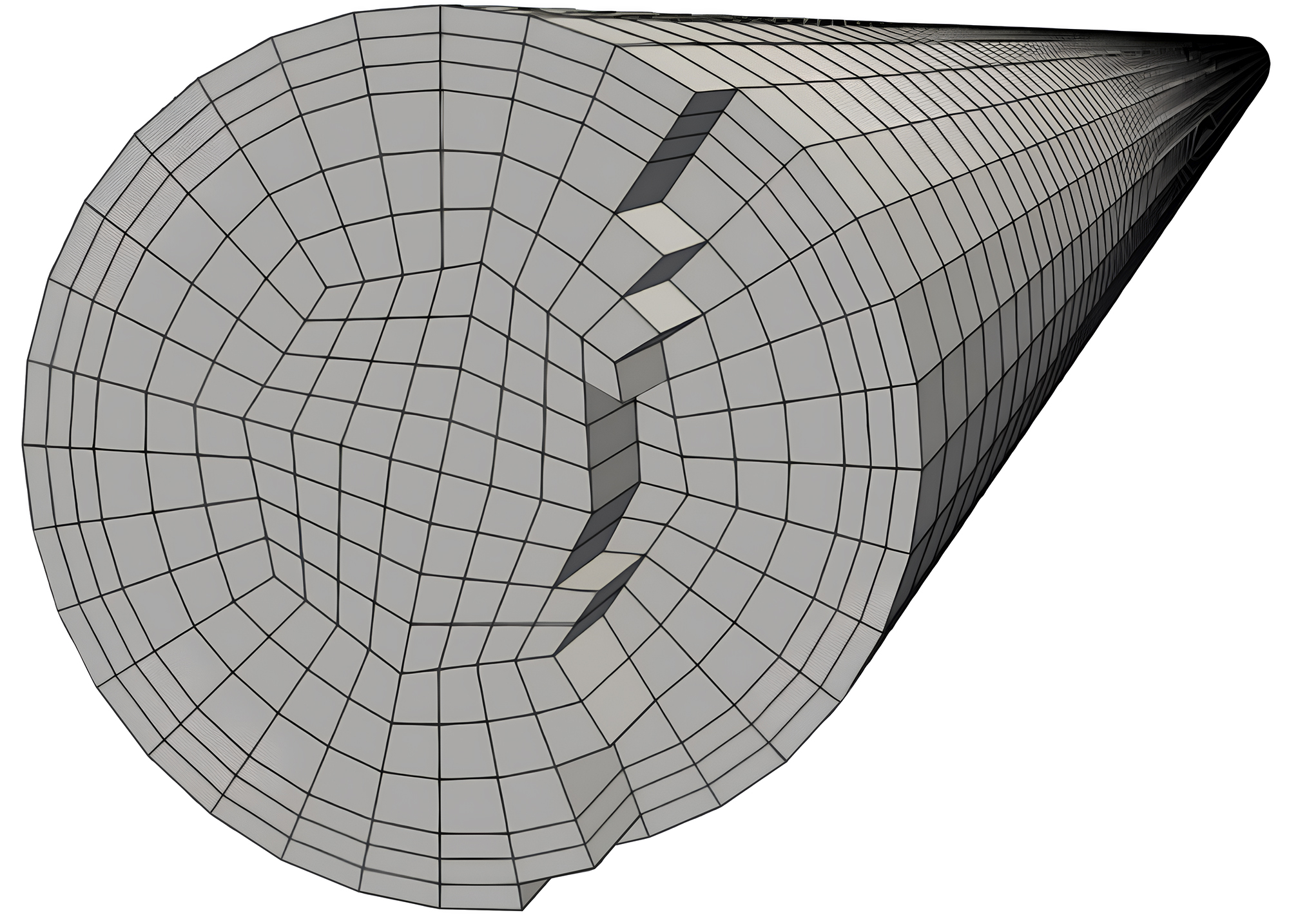}
\caption{Left: Our coarse mesh. Right: Coarse mesh from \cite{decroocq2023modeling}.}
\label{Fig: structured_comparison}
\end{figure}

In this section, we proceed our comparison with other mesh generation methods. In particular, we compare our approach with one of the most recent works on the generation of structured grids in numerical h{\ae}modynamics. In \cite{decroocq2023modeling}, the authors evaluate the performance of the mesh generator by means of a test case represented by a cylindrical geometry with a diameter of $2.5$ mm and a length of $200$ mm. For a fair comparison, we adopt the same test case and compare our approach to the one proposed in \cite{decroocq2023modeling}. We compare different levels of refinement: coarse, medium and fine while monitoring the grid quality as measured by SJ. In Fig. \ref{Fig: structured_comparison}, we depict the coarsest mesh from \cite{decroocq2023modeling} to highlight the differences in the mesh generation approach.
Tab. \ref{tab: structured_grid_comparison}, shows a complete comparison of the SJ values. As seen from the values, with our spline-based approach, we are able to maintain a minimum SJ whose value exceeds the corresponding value from \cite{decroocq2023modeling} by $20\%$, $17\%$, and $10\%$ in the coarse, medium, and fine grids, respectively. Consequently, our mesh generator exceeds the corresponding average SJ value by more than $3\%$, while the values of the maximum SJ are comparable.
Regrettably, we do not currently have access to a readily reproducible test case for conducting a comparative analysis of the bifurcation mesh.

\begin{table}[h!]
\centering
\caption{Mesh comparison with \cite{decroocq2023modeling}. In parenthesis the number of nodes of each grid. SJ stands for \text{scaled Jacobian}.}
\begin{tabular}{|c|c|c|c|c|c|c|c|c|c|}
\hline
 &\multicolumn{3}{|c|}{Coarse ($\sim{105}k$)} & \multicolumn{3}{c|}{Medium ($\sim{870}$k)} & \multicolumn{3}{c|}{Fine ($\sim{8800}$k)} \\
\hline
& $\text{SJ}_{\text{min}}$ & $\text{SJ}_{\text{mean}}$ & $\text{SJ}_{\text{max}}$ & $\text{SJ}_{\text{min}}$ & $\text{SJ}_{\text{mean}}$ & $\text{SJ}_{\text{max}}$ & $\text{SJ}_{\text{min}}$ & $\text{SJ}_{\text{mean}}$ & $\text{SJ}_{\text{max}}$\\
\hline 
\text{ours} & 0.861 & 0.981 & 0.999 & 0.795 & 0.991 & 0.999 & 0.767 & 0.990 & 0.999 \\
\cite{decroocq2023modeling} & 0.707 & 0.960 & 0.995 & 0.707 & 0.971 & 0.999 & 0.707 & 0.976 & 0.999 \\
\hline
\end{tabular}
\label{tab: structured_grid_comparison}
\end{table}

Regarding the methodology, both techniques provide the flexibility to incorporate an arbitrary number of vessel sections along the centerline. However, it is noteworthy that our method generates a spline volume mesh in both the radial and longitudinal directions. Specifically, the points on the vessel section are not radially projected from the center to the vessel surface, unlike in previous studies \cite{de2010patient, de2011patient,ghaffari2017large,decroocq2023modeling}. Our approach enables a lighter representation of the mesh, as it is encoded in a relatively low number of spline control points. The framework allows for a computationally inexpensive sampling of a finer mesh from the spline mapping operator, see Sec.\ \ref{sec: single_branch}.
Moreover, our method does not necessitate the information on the vessel surface in any step of the mesh generation, neither for the single branch (\cite{zhang2007patient, de2011patient} need the vessel surface as starting point of the mesh generator) nor for the bifurcation (\cite{decroocq2023modeling} needs the information on vessel surfaces for computing their intersection in the context of bifurcation). 
Differently from \cite{de2010patient, de2011patient}, our strategy is readily generalized to non-planar bifurcations, see Sec. \ref{sec: bifurcation}.
We also propose a novel generalization to non-planar intersecting branches, which, to the best of our knowledge, has not been previously introduced in the literature, see Sec.\ \ref{sec: n_furcations}.  
A summary of the main differences with respect other studies in the context of mesh generation for structured grid has been already reported in Tab. \ref{tab: state_of_the_art}.

\begin{figure}[t]
\centering
\includegraphics[width = 15cm, keepaspectratio]{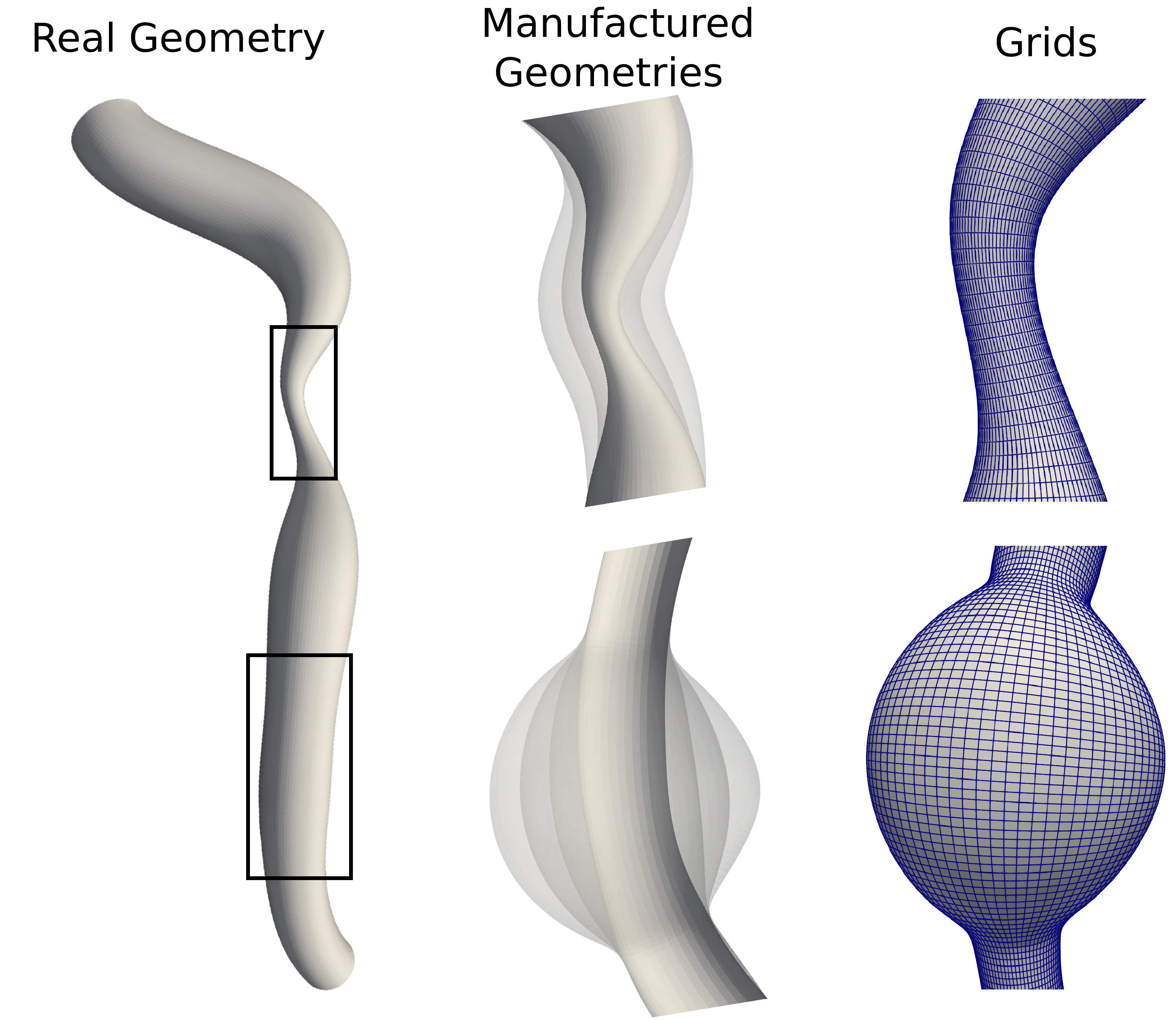}
\caption{Geometry manipulation. The original patient-specific geometry is altered in two different scenario. 
In top part, a stenotic condition is restored to a healthy condition showing two different degree of modification (medium stenosis, healthy condition). In the bottom tract, an aneurysm is created showing three different levels of severity (small, medium and big).}
\label{Fig: geometry_manipulation}
\end{figure}

\subsubsection{Geometry manipulation: incorporating pathologies}\label{sec: topology}

Due to the spline parameterization of both the centerline and radius data, the geometry of the vessel can be modified without altering its topology, i.e.\ the number of total vessel sections and the number of nodes inside each vessel section.
This approach provides the flexibility to adjust the curvature and reshape the entire or part of the vessel structure, as \cite{decroocq2023modeling} already show in their work. We emphasize that the mesh connectivity remains invariant throughout these geometrical modifications. 

\begin{figure}[t!]
\centering
\includegraphics[width = 5.3cm, keepaspectratio]{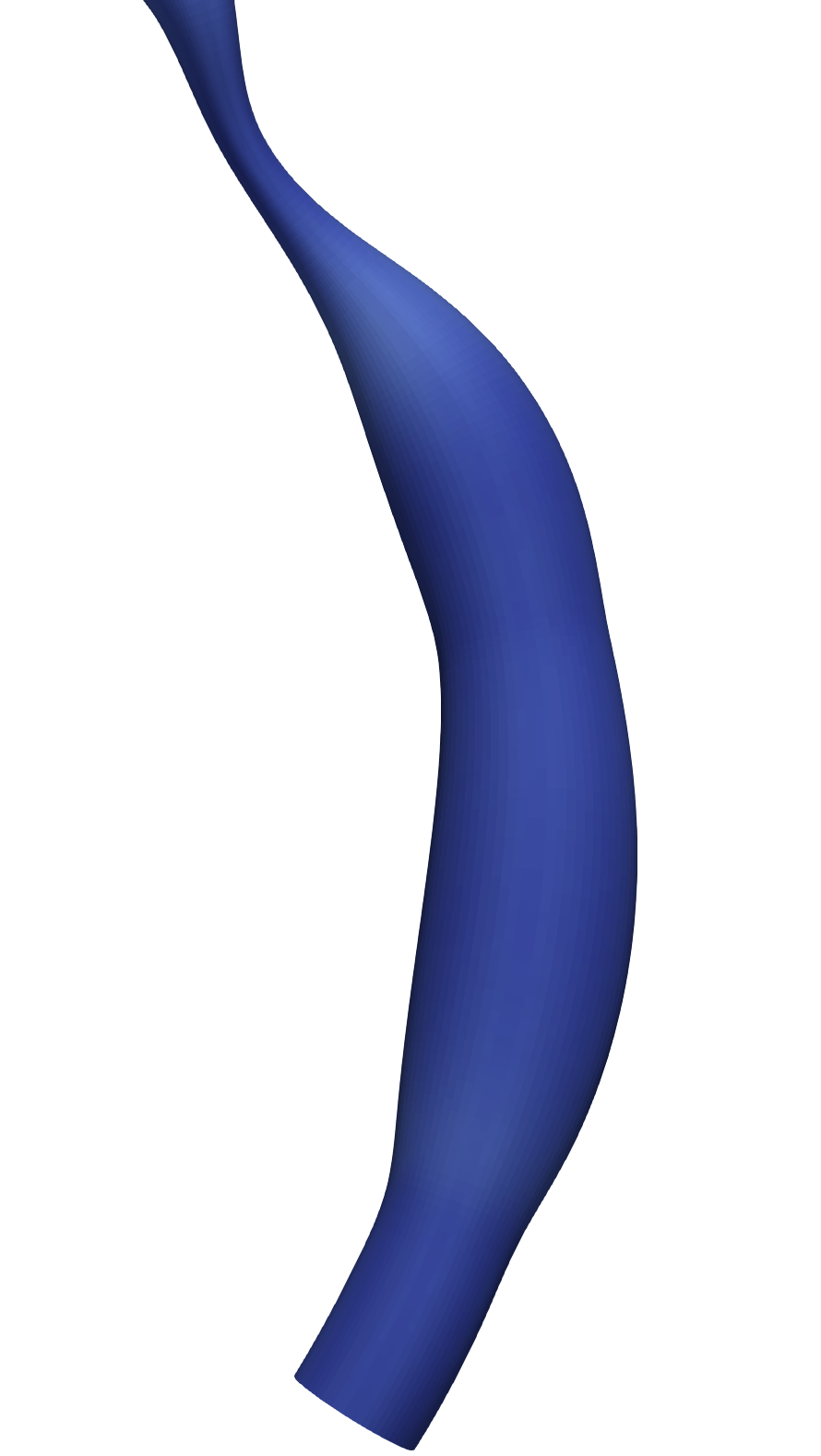}
\includegraphics[width = 5.3cm, keepaspectratio]{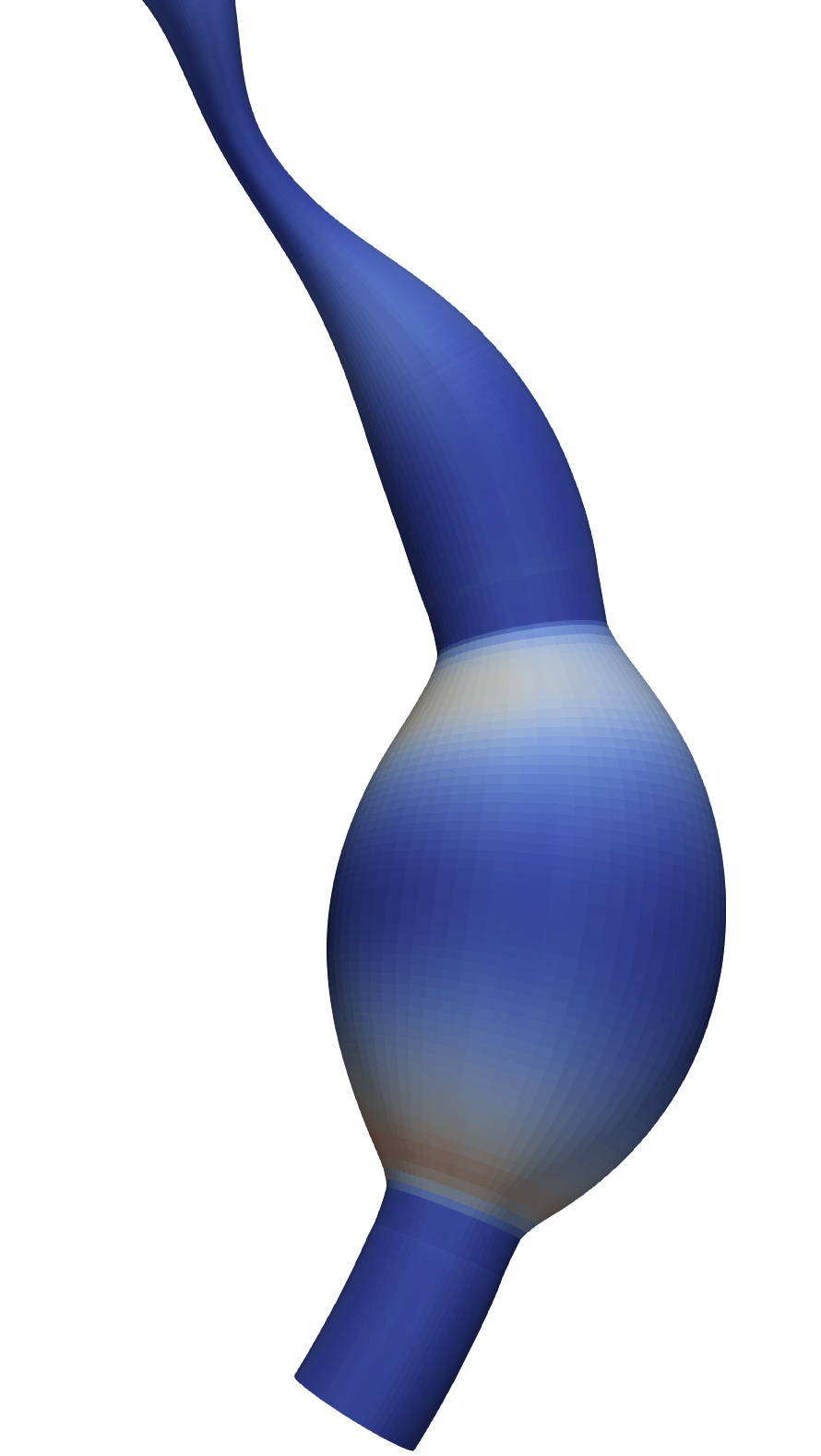}
\includegraphics[width = 5.3cm, keepaspectratio]{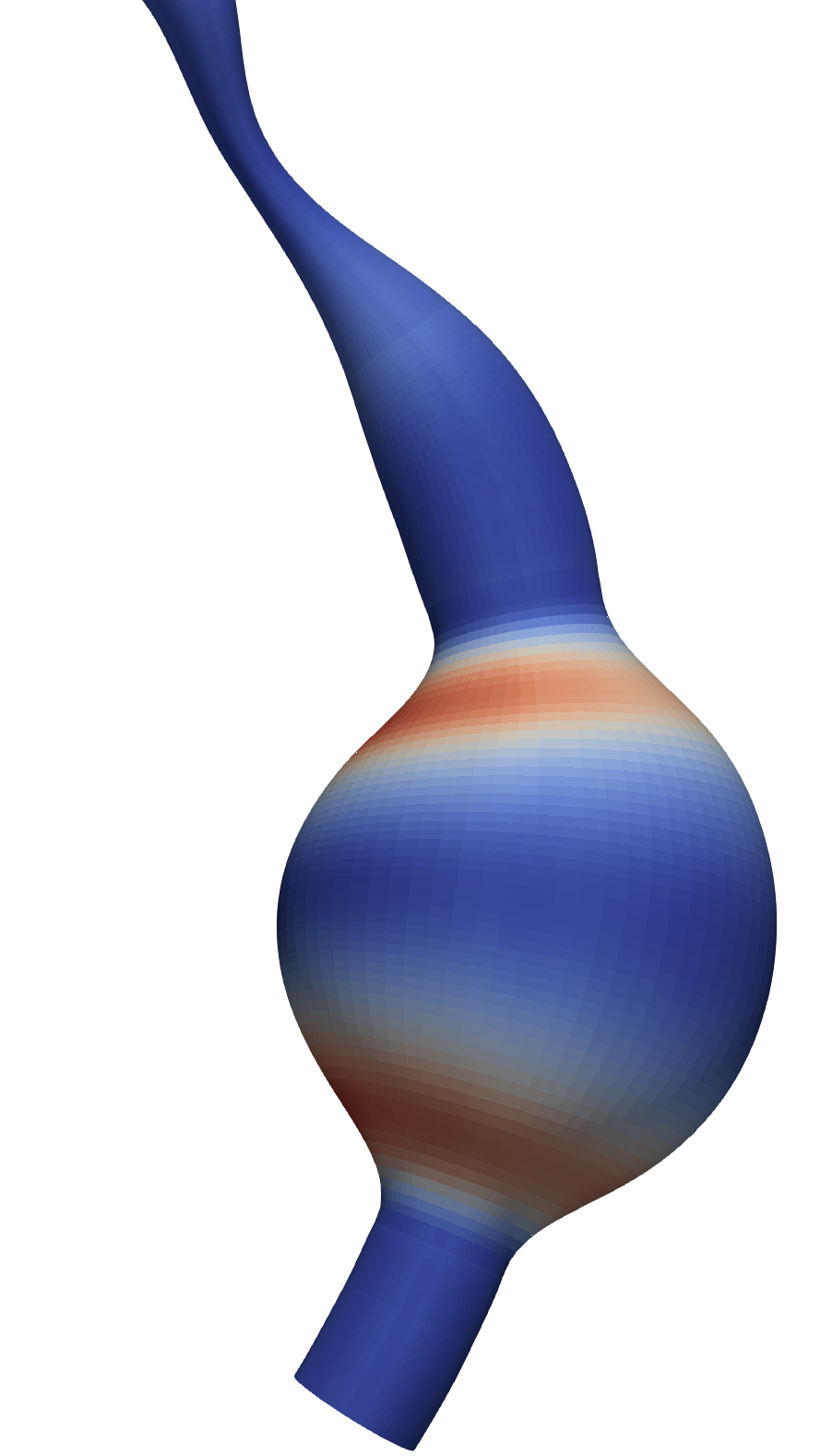}
\includegraphics[width = 8.5cm, keepaspectratio]{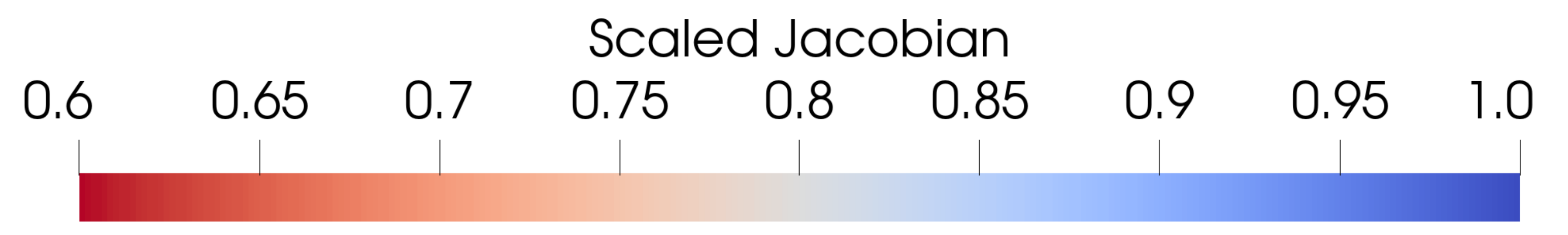} 
\caption{Scaled Jacobian visualization for the small, medium and big aneurysm.}
\label{Fig: geometry_manipulation_aneurysm_jacobian}
\end{figure}

Fig.\ \ref{Fig: geometry_manipulation} shows a patient-specific coronary artery, where in the top part a stenosis is adjusted to restore a healthy condition, while in the bottom part, a healthy tract is modified into an aneurysm.
\textcolor{black}{It is worth mentioning that there are two main types of aneurysms: fusiform and saccular.
The aneurysm depicted in Fig.\ \ref{Fig: geometry_manipulation} is of the fusiform type. 
Since the proposed mesh generator is focuses on coronary arteries, which are almost never subjected to saccular aneurysms, the latter case is not studied.}
At the center of Fig.\ \ref{Fig: geometry_manipulation}, different levels of geometry modification are depicted. 
More precisely, starting from a severe stenosis in the top section, a medium stenosis and a healthy condition are constructed, while in the bottom tract of the vessel, a small, medium and big aneurysms are introduced (see manufactured geometries).
In Fig.\ \ref{Fig: geometry_manipulation_aneurysm_jacobian}, we report the $3$D SJ visualization of the different aneurysms of the bottom tract. Notice that, as expected, the biggest aneurysm generates the lowest SJ values while maintaining an excellent mesh quality nonetheless.
All the meshes are created using a boundary layer refinement of $\alpha = 0.5$ (see eq. \ref{eq: bl}).

In Tab. \ref{tab: geometry_manipulation}, we compute the SJ and the NES for the different meshes.
Note that when introducing a pathology in one vessel segment, we maintain the real geometry condition in the other, i.e., for instance, a stenosis in the top part and an unaltered bottom part.

A central feature of these modifications is the preservation of the vessel surface smoothness at the transition interface between the original and modified vessel segments. Here, failure to retain a $C^{1}$-continuous vessel radius spline,  would affect the mesh quality, thus compromising the computation of velocity-gradient based h{\ae}modynamic quantities such as the wall shear stress (WSS), see Sec.\ \ref{sec: fluid_results}.
To maintain the smoothness, we directly act on the geometry's vessel radius spline.
As introduced in Sec.\ \ref{sec: single_branch}, the radius data is described by a scalar function $r(t): [0, 1] \rightarrow \mathbb{R}^+$. 
First of all, we detect the region of interest (the vessel's tract requiring modifications), defined by a starting and ending value in the radius spline's parametric domain, $t_{\text{start}}$ and $t_{\text{end}}$, respectively, with $0 \leq t_{\text{start}} < t_{\text{end}} \leq 1$.
\begin{table}[h!]
\caption{Computation of the index quality for the different geometrical configurations. SJ stands for \text{scaled Jacobian} and NES for \text{normalized equiangular skewness}. All three meshes have a total of 328'250 vertices.}
\centering
\begin{tabular}{|c|c|c|c|c|c|c|}
\hline
 & \multicolumn{6}{|c|}{Geometry Manipulation}  \\
\hline
\text{Top Part} & $\text{SJ}_{\text{min}}$ & $\text{SJ}_{\text{mean}}$ & $\text{SJ}_{\text{max}}$ & $\text{NES}_{\text{min}}$ & $\text{NES}_{\text{mean}}$ & $\text{NES}_{\text{max}}$ \\
\hline
Healthy  & $0.836$ & $0.989$ & $0.999$ &  $0.005$ & $0.074$ & $0.367$  \\
Medium stenosis & $0.825$ & $0.988$ & $0.999$ & $0.005$ & $0.077$ & $0.367$  \\
Stenosis & $0.820$ & $0.986$ & $0.999$ & $0.005$ & $0.083$ & $0.368$  \\
\hline
\hline
Bottom Part (in Fig.\ \ref{Fig: geometry_manipulation_aneurysm_jacobian} )& $\text{SJ}_{\text{min}}$ & $\text{SJ}_{\text{mean}}$ & $\text{SJ}_{\text{max}}$ & $\text{NES}_{\text{min}}$ & $\text{NES}_{\text{mean}}$ & $\text{NES}_{\text{max}}$ \\
\hline
Small aneurysm & $0.824$ & $0.985$ & $0.999$ & $0.005$ & $0.083$ & $0.368$  \\
Medium aneurysm & $0.719$ & $0.974$ & $0.999$ & $0.005$ & $0.108$ & $0.478$  \\
Big aneurysm & $0.621$ & $0.965$ & $0.999$ & $0.005$ & $0.122$ & $0.564$  \\
\hline
\end{tabular}
\label{tab: geometry_manipulation}
\end{table}
We parameterize the to-be-modified region using a cubic spline, interpolating $r(t_{start})$, $r(t_{end})$ and
two other values of the radius that could be positioned inside any part of the parametric domain, granting complete shape modification freedom.
The modified vessel segment's radius is then represented by a cubic spline that interpolates the four provided values. By adjusting the values, various shapes can be obtained.
To smoothly recombine the original splines, i.e.\ the ones corresponding to the tracts $t \in [0,t_{begin})$ and $t \in (t_{end},1]$, with the modified tract $t \in [t_{begin}, t_{end}]$ (describing the region of interest), the three splines are densely evaluated and a re-interpolation of all the evaluation points is performed. 
\textcolor{black}{To avoid wiggles near the interface region, a smoothing parameter is employed for the new spline evaluation. 
The smoothing parameter acts on a global level, meaning that it penalizes the norm of the spline curve's second derivative. Note that, this could, in principle, create a deviating radius also in other regions, not only close to the interface. However, this is essential to maintain the smoothness.}
The resulting spline can be evaluated in the same number of points of the original geometry spline, maintaining the same global topology.


\subsubsection{\textcolor{black}{Complex coronary artery tree}}\label{sec: complex_tree}
\textcolor{black}{Our method allows to generate meshes of more complex vessel geometries, for example a coronary tree that consider the right, left and circumflex coronary arteries and their main bifurcations.}
\textcolor{black}{The presence of multiple bifurcations requires handling the mutual torsion between vessel tracts approaching a bifurcation and the vessel tract within two bifurcations. 
This is resolved by introducing an additional rotation that perfectly aligns the last section of the single branch with the first section of the approaching bifurcation.
For a branch starting from a bifurcation and extending to the end, the no twist frame is initialized with the same rotational frame as that of the corresponding part of the bifurcation.
}

\begin{figure}[h!]
\centering
\includegraphics[width = 10cm, keepaspectratio]{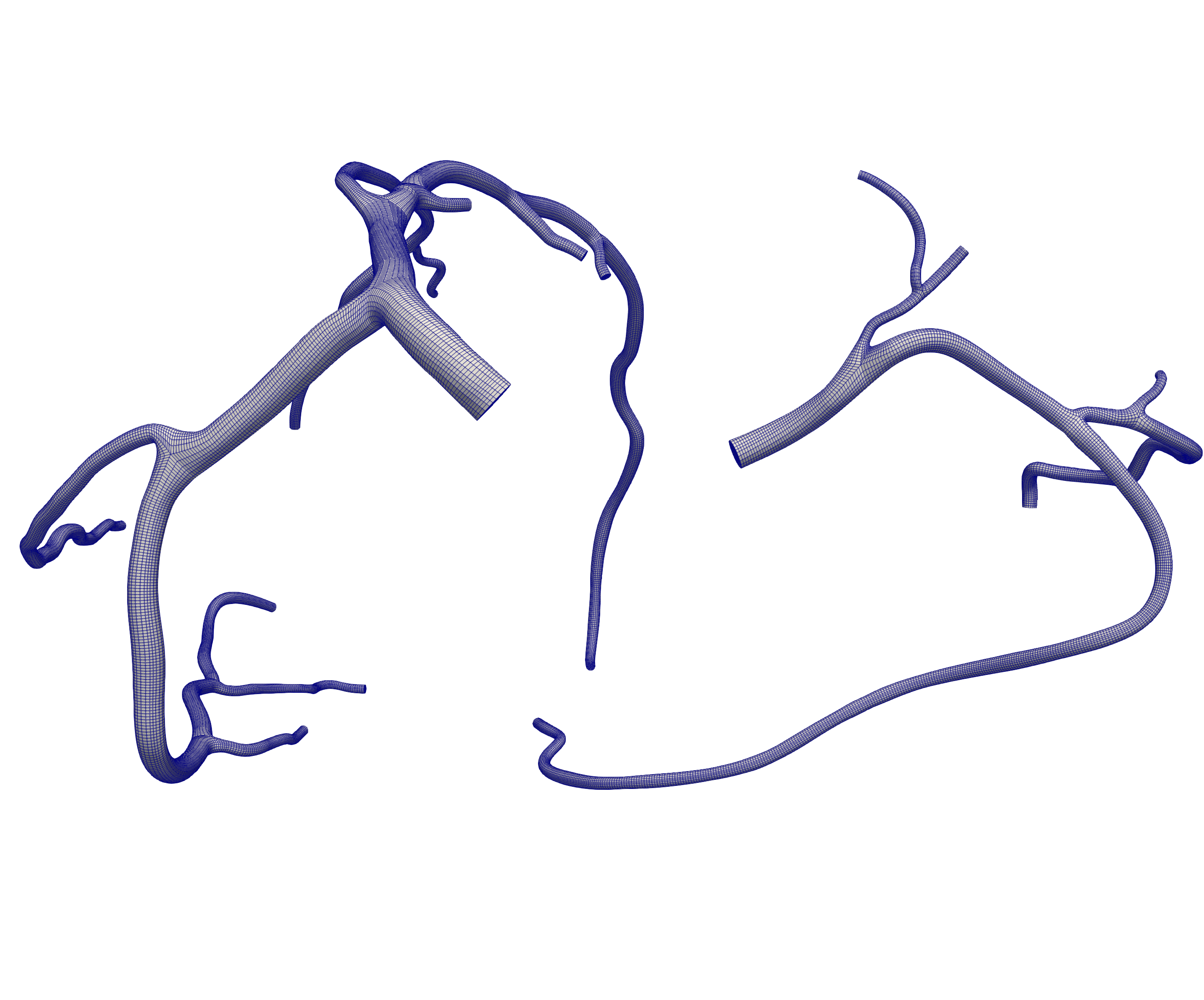}\quad
\includegraphics[width = 9cm, keepaspectratio]{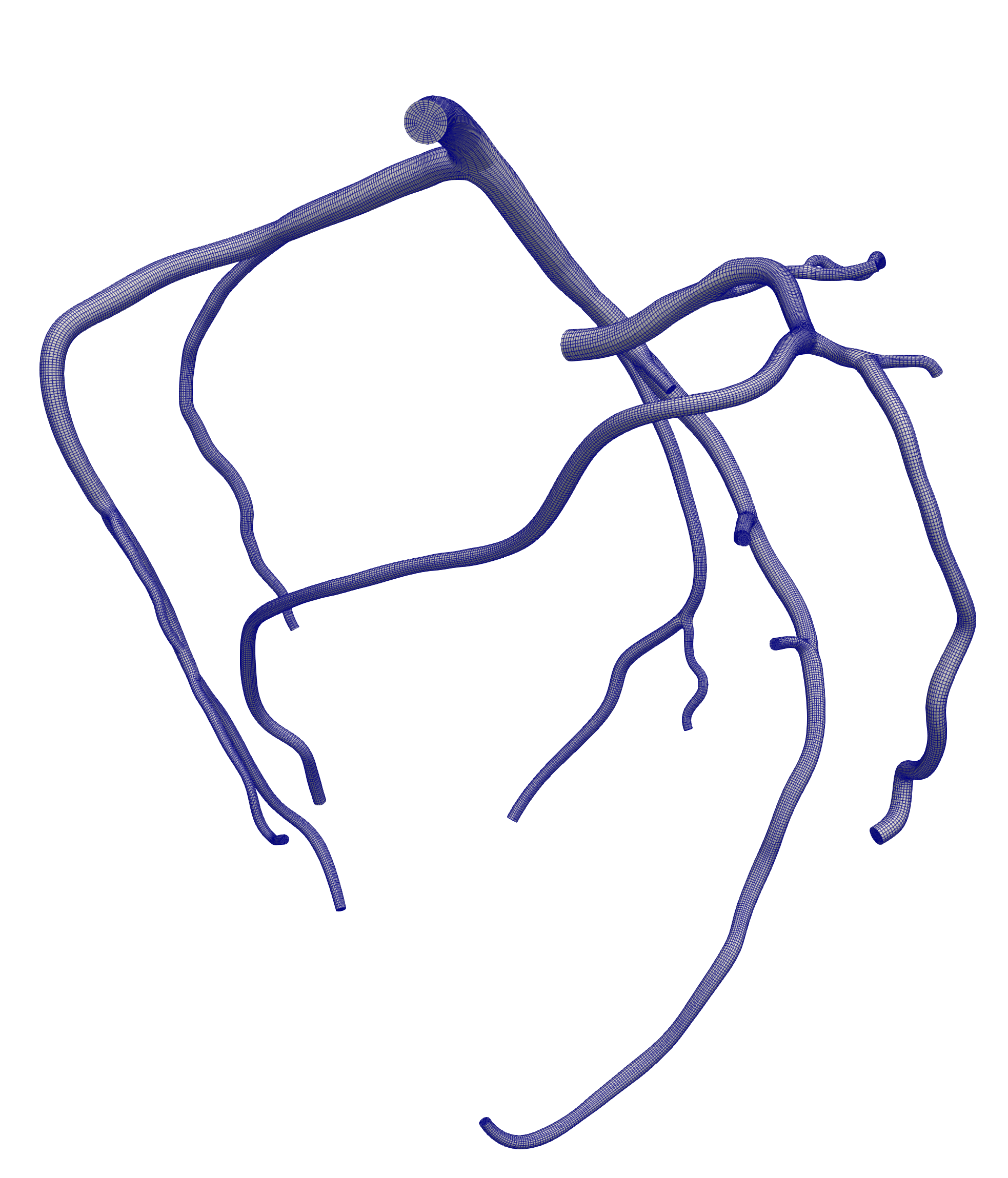}\quad
\caption{\textcolor{black}{Top: Top view of the coronary tree. Bottom: Right view of the coronary tree}}
\label{Fig: complex_mesh}
\end{figure}

\begin{figure}[h!]
\centering
\includegraphics[width = 7cm, keepaspectratio]{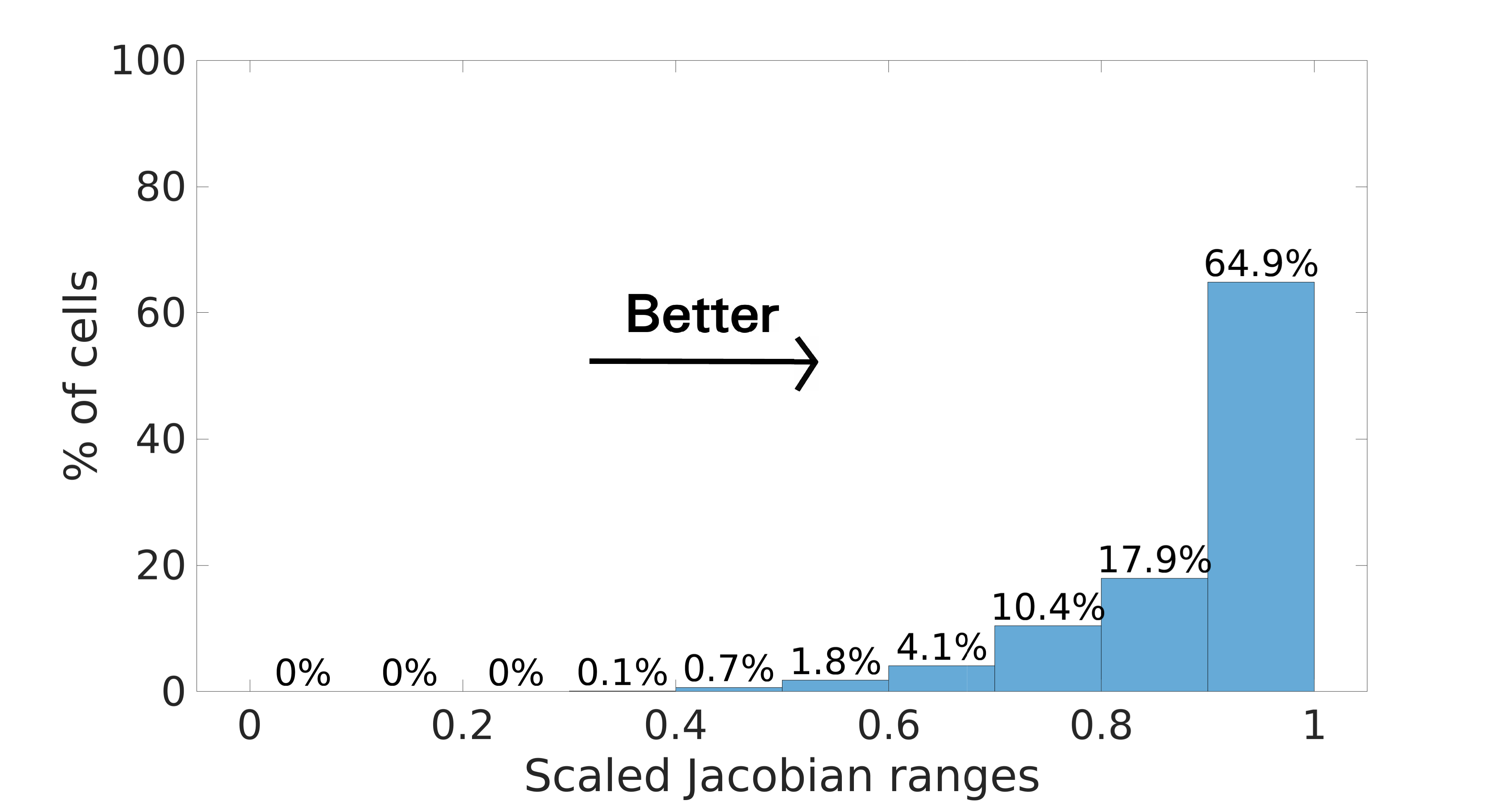}
\includegraphics[width = 7cm, keepaspectratio]{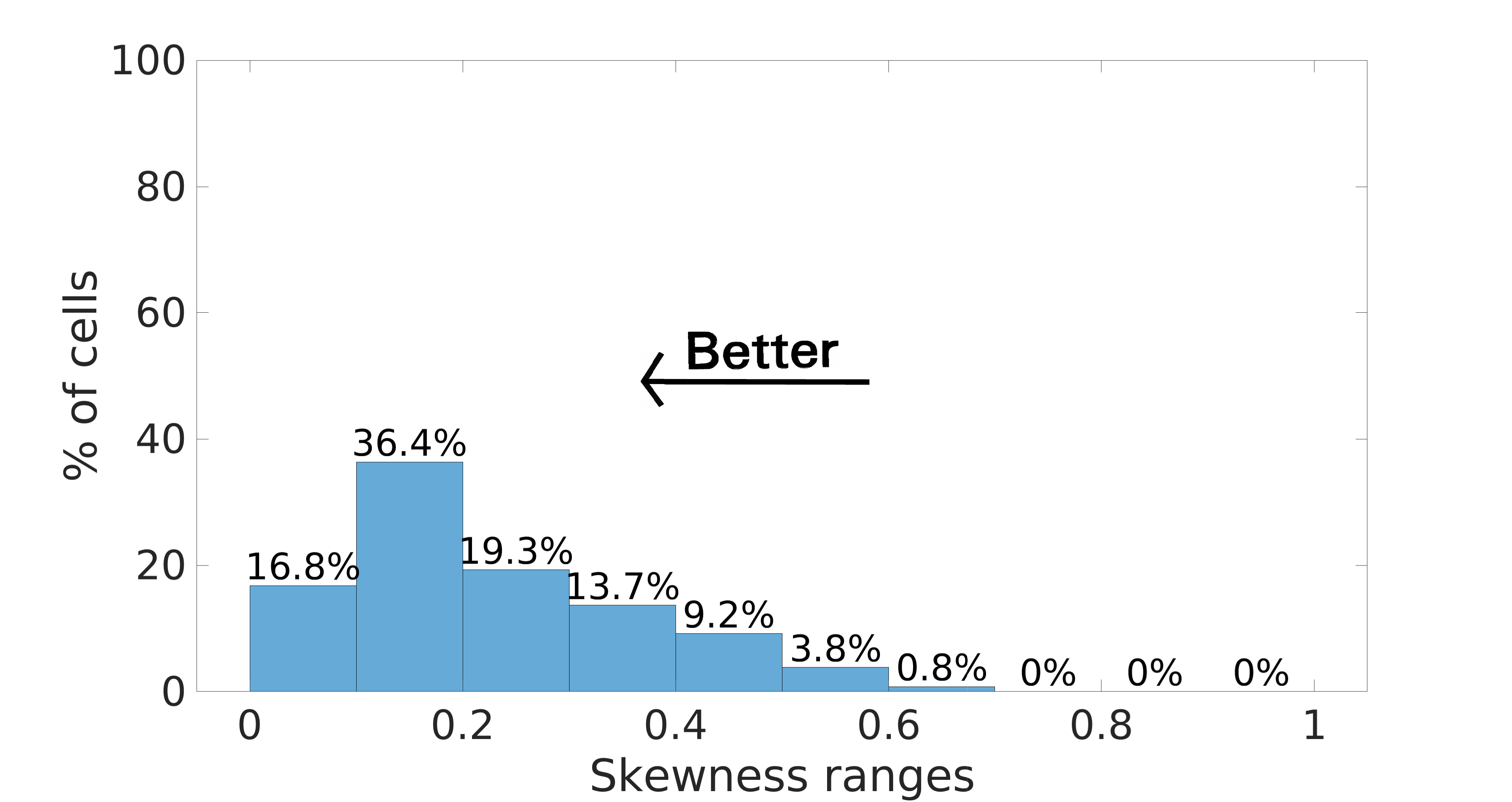}
\includegraphics[width = 7cm, keepaspectratio]{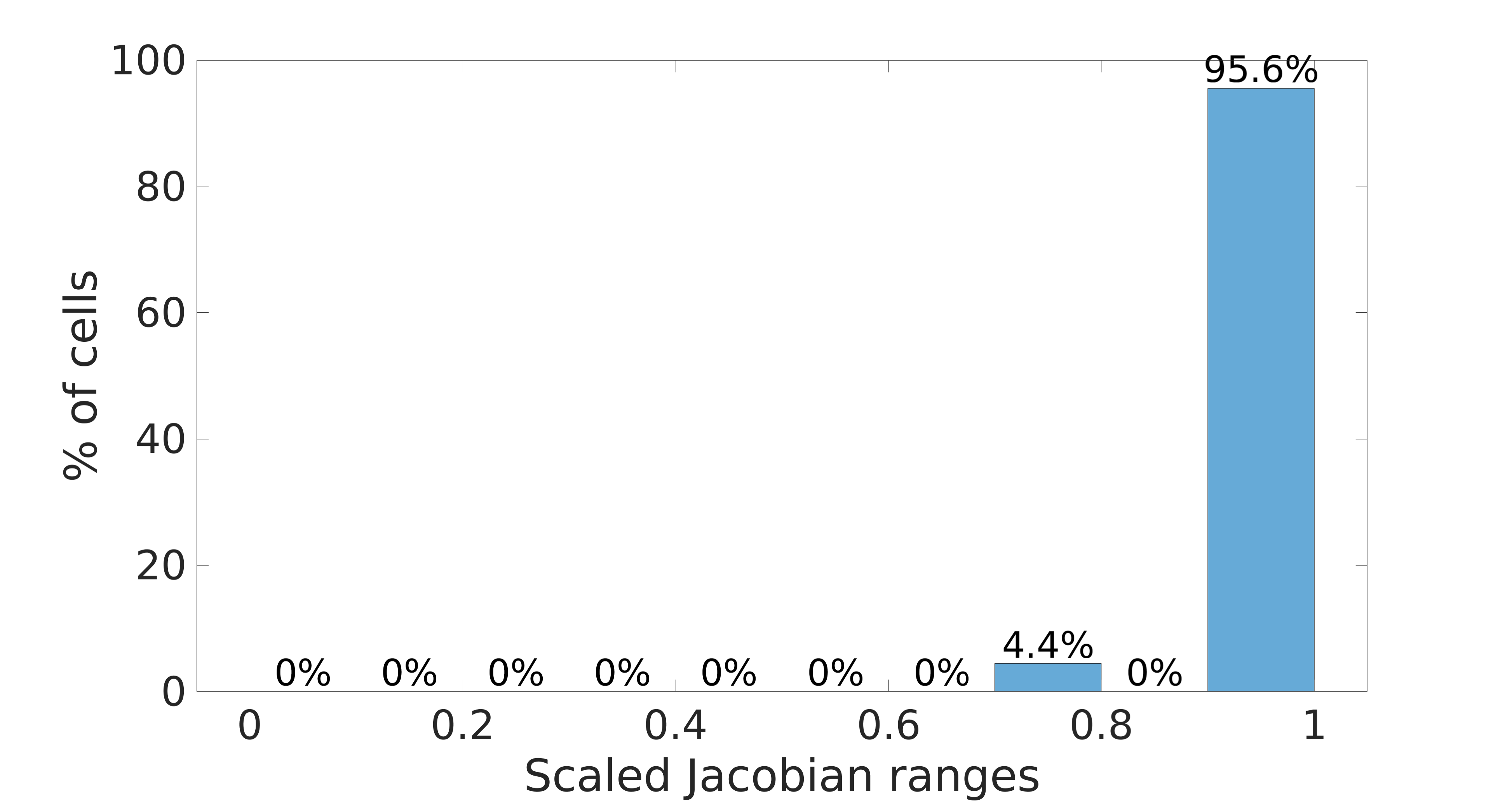}
\includegraphics[width = 7cm, keepaspectratio]{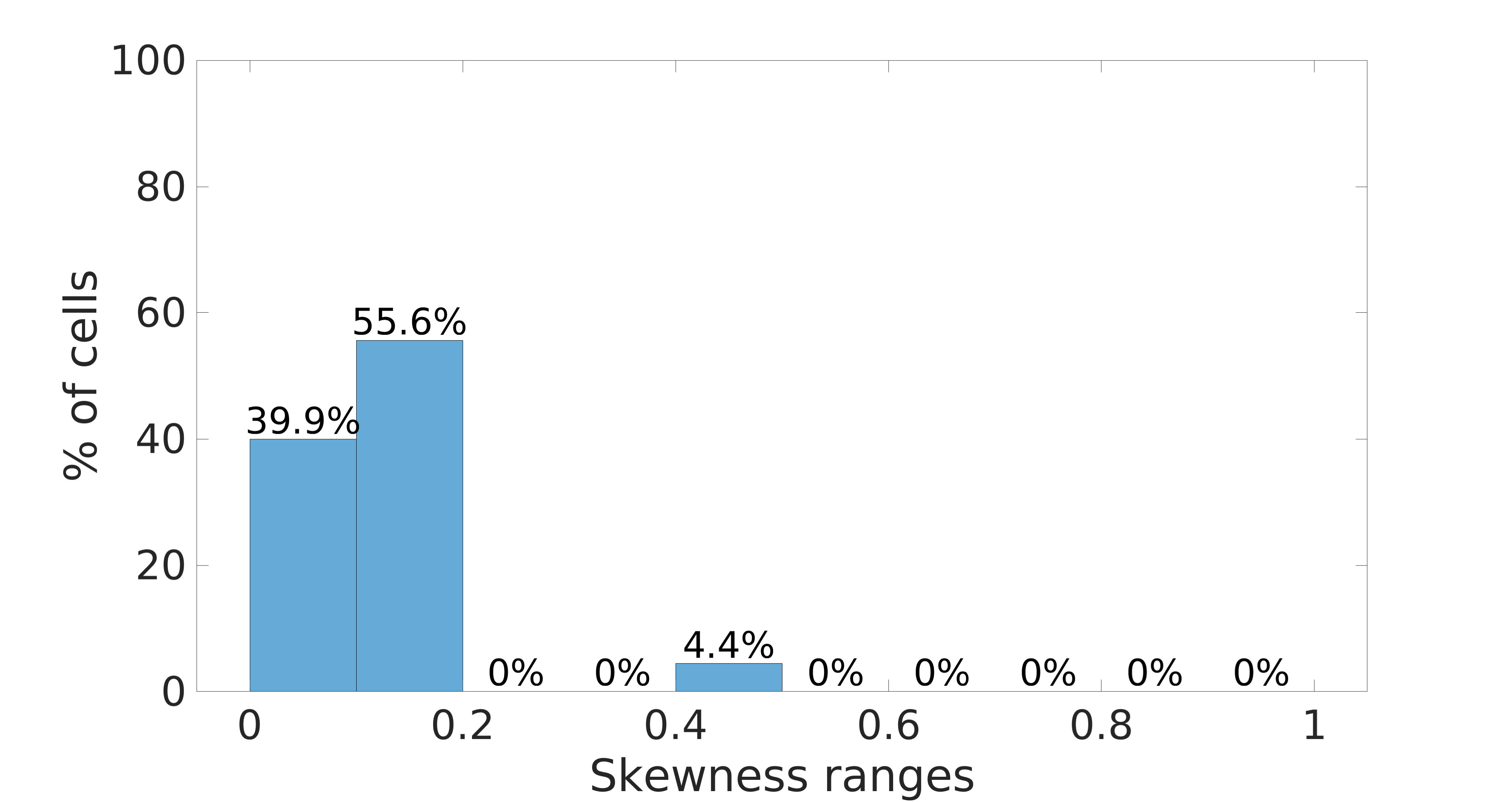}
\caption{\textcolor{black}{
The second row contains the values of the SJ and NES measures for all the bifurcations of the meshes reported in Fig.\ \ref{Fig: complex_mesh}.
The first row contains the values of the SJ and NES measures for all the single branch tracts of the mesh reported in Fig.\ \ref{Fig: complex_mesh}.
}}
\label{Fig: complex_stats}
\end{figure}

\textcolor{black}{In Fig.\ \ref{Fig: complex_mesh}, we depict the two different meshes of a coronary artery tree, while in in Fig.\ \ref{Fig: complex_stats}, as previously done in Sec.\ \ref{sec: grid_comparison}, we illustrate the percentage of cells falling within specified ranges of the SJ and NES. Tab.\ \ref{tab: complex_stats} summarizes the results for the bifurcation and single branches.
It is worth mentioning that the Hermite-based modeling of the bifurcations allows for handling various different inputs without manual intervention. In particular, all the Hermite curves (defining the bifurcation boundaries, see Sec. \ref{sec: bifurcation}) are scaled with a parameter proportional to the diameters of the three different sections, so that the curve automatically adapt to each use case.
Both quality metrics (SJ and NES) indicate that the grids show excellent quality for the majority of the cells.
}

\begin{table}[h!]
\caption{\textcolor{black}{Quality indices for the single branch tract and bifurcation of the vascular tree depicted in Fig.\ \ref{Fig: complex_mesh}. SJ stands for \text{scaled Jacobian} and NES for \text{normalized equiangular skewness}.}}
\centering
\begin{tabular}{|c|c|c|c|c|c|c|}
\hline
 & \multicolumn{6}{|c|}{Coronary trees}  \\
\hline
 & $\text{SJ}_{\text{min}}$ & $\text{SJ}_{\text{mean}}$ & $\text{SJ}_{\text{max}}$ & $\text{NES}_{\text{min}}$ & $\text{NES}_{\text{mean}}$ & $\text{NES}_{\text{max}}$ \\
\hline 
\text{Single branches} & 0.752 & 0.971 & 0.999 & 0.0017 & 0.125 & 0.45 \\
\text{Bifurcations} & 0.228 & 0.898 & 0.999 & 0.0003 & 0.227 & 0.729  \\
\hline
\end{tabular}
\label{tab: complex_stats}
\end{table}

\subsection{Fluid dynamics results}\label{sec: fluid_results}

To demonstrate the robustness of our mesh generator, we carry out numerical simulations of a coronary vessel segment that includes the Left Anterior Descending artery (LAD) and the Left Circumflex artery (LCx), see Fig.\ \ref{Fig: coronary_domain}.

\subsubsection{The geometric multiscale model}\label{sec: fluidynamics}
\begin{figure}[h!]
\centering
\includegraphics[width = 8cm, keepaspectratio]{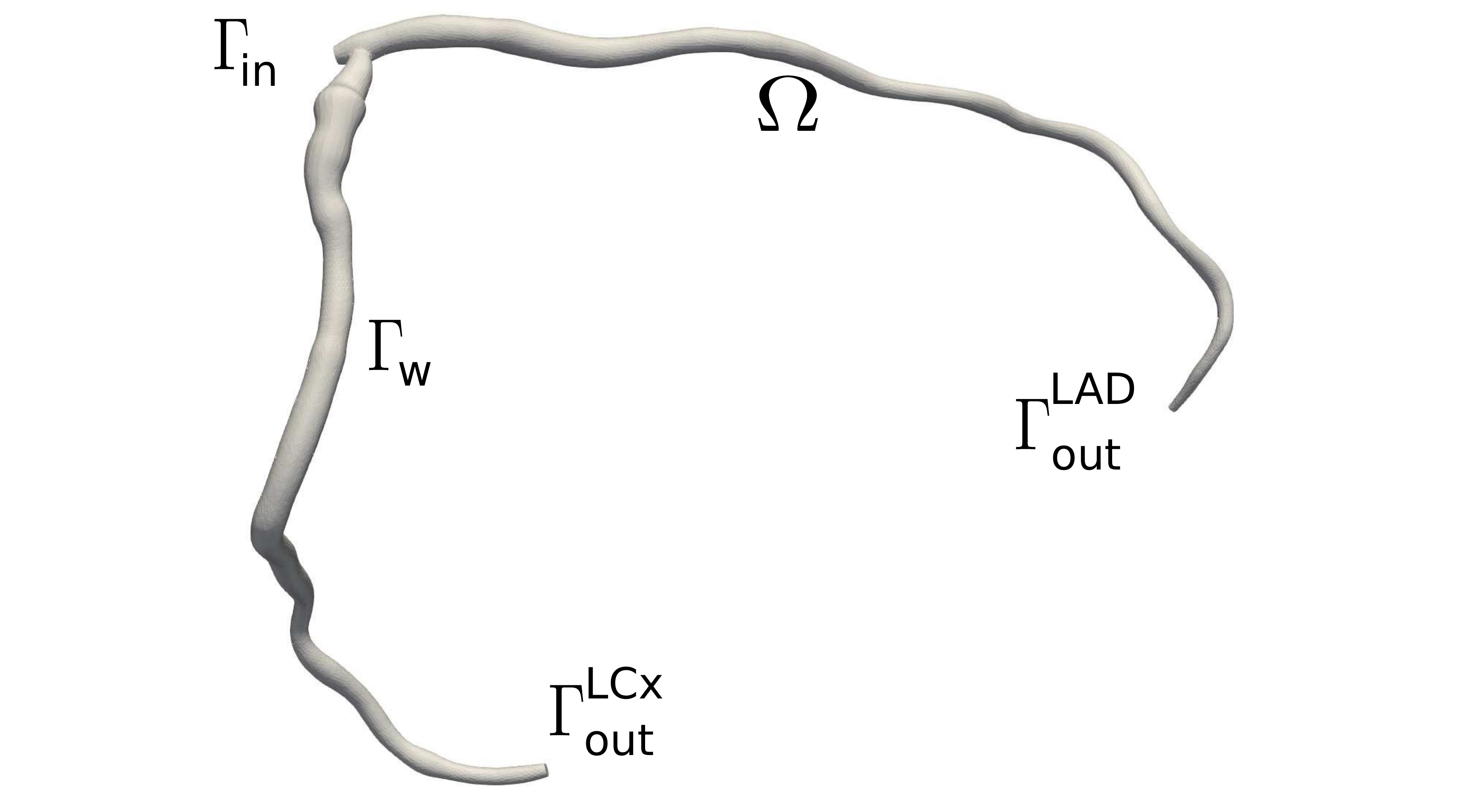}
\includegraphics[width = 8cm, keepaspectratio]{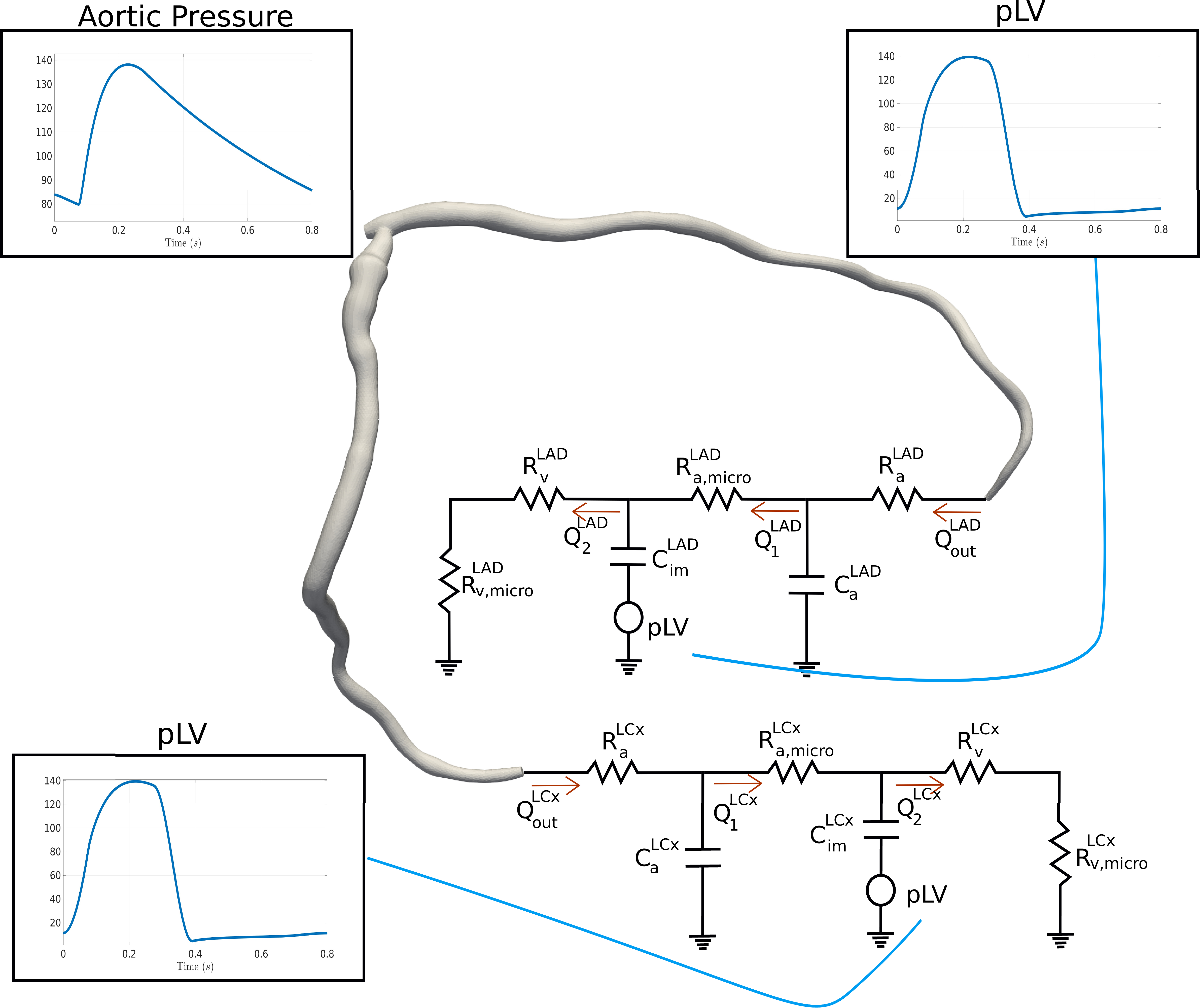} 
\caption{Left: Computational domain of the coronary artery. Right: Geometric multiscale model of the coronary arteries and the downstream circulation. In the figure, one outlet is considered as the Left Anterior Descending artery (LAD), while the other one is considered as the Left Circumflex artery (LCx).}
\label{Fig: coronary_domain}
\end{figure}

Blood is a strongly heterogeneous fluid. However, for demonstration purposes, we consider it as a homogeneous, incompressible and Newtonian fluid. 
As reported in Fig.\ \ref{Fig: coronary_domain}, we define $\Omega$ as the computational domain representing a vessel segment, 
$\operatorname{\Gamma_{w}}$ as the vessel walls while $\operatorname{\Gamma_{in}}$ is the inlet at the upstream section. The two outlets are labeled as  
$\operatorname{\Gamma}_{\operatorname{out}}^{\operatorname{LAD}}$ and $\operatorname{\Gamma}_{\operatorname{out}}^{\operatorname{LCx}}$ as outlets. 
Therefore, the fluid dynamics of the coronary arteries is modeled by means of the $3$D incompressible Navier-Stokes (NS) equations that describe the blood velocity $\mathbft{u}(\mathbft{x},t):\Omega\times\mathbb{R^+} \to \mathbb{R}^3 $
and the pressure $p(\mathbft{x},t):\Omega \times \mathbb{R^+} \to \mathbb{R}$. 
Defining $T$ as the period of the heartbeat, the formulation reads:\\
Find $\mathbft{u}$, $p$ in $t\in(0,T]$, such that:
\begin{align}
&	\displaystyle \rho \frac{\partial\mathbft{u}}{\partial t} + \rho (\mathbft{u} \cdot \nabla)\mathbft{u}
     - \nabla \cdot \sigma(\mathbft{u},p) = \mathbft{0}  && \text{in} \quad \Omega, \\
&    \nabla\cdot{\mathbft{u}}=0  &&  \text{in} \quad \Omega,
\label{eq: navier_stokes}
\end{align}

where $\sigma(\mathbft{u},p) = \mu(\nabla \mathbft{u} + \nabla \mathbft{u}^{T}) - pI$ is the Cauchy fluid stress tensor, $\mu$ stands for the constant dynamic viscosity and $\rho$ is the constant density. 

As an initial condition, we impose zero velocity and zero pressure along the entire domain.
Since we perform a rigid wall simulation, we impose the no-slip condition on the vessel walls:
\begin{equation*}
\mathbft{u}(\mathbft{x},t) = \mathbft{0}  \quad\quad\quad \text{on} \quad \operatorname{\Gamma_{w}}.
\end{equation*}

At the inlet, a Neumann condition is imposed through a pressure data generated offline by means of a closed-$0$D model of the entire cardiovascular system which imposes the aortic pressure as inlet boundary condition.
The outlet boundary conditions are imposed by means of the $0$D model introduced in \cite{kim2010patient}.
This lumped parameter model is able to account for the effects of the microvasculature in the downstream regions. All the values of the lumped parameters are taken from \cite{tajeddini2020high}, see Fig.\ \ref{Fig: coronary_domain}, right.

Concerning the numerical discretization, both the NS and $0$D model equations are discretized in time using a Backward Euler (BE) discretization \cite{kim2010patient}.
The $3$D equations are discretized in space using $P1-P1$ finite elements with the variational multiscale-streamline upwind Petrov-Galerking (VMS-SUPG) stabilization \cite{tezduyar2000finite,forti2015semi} and backflow stabilization \cite{bertoglio2014tangential}. 
The linear system arising after the space-time discretization is preconditioned using SIMPLE \cite{patankar1983calculation} where the inverse of the Schur complement is approximated using the algebraic multigrid method \cite{ruge1987algebraic}. The system is solved using the generalized minimum residual method (GMRES, \cite{saad1986gmres}) with an absolute tolerance of $10^{-12}$.
The numerical simulations are performed in $\operatorname{life}^\textsc{x}$ \cite{africa2022flexible}, a high-performance object-oriented finite element library focused on the mathematical models and numerical methods for cardiac applications.

\begin{figure}[t]
\centering
\includegraphics[width = 7.5cm, keepaspectratio]{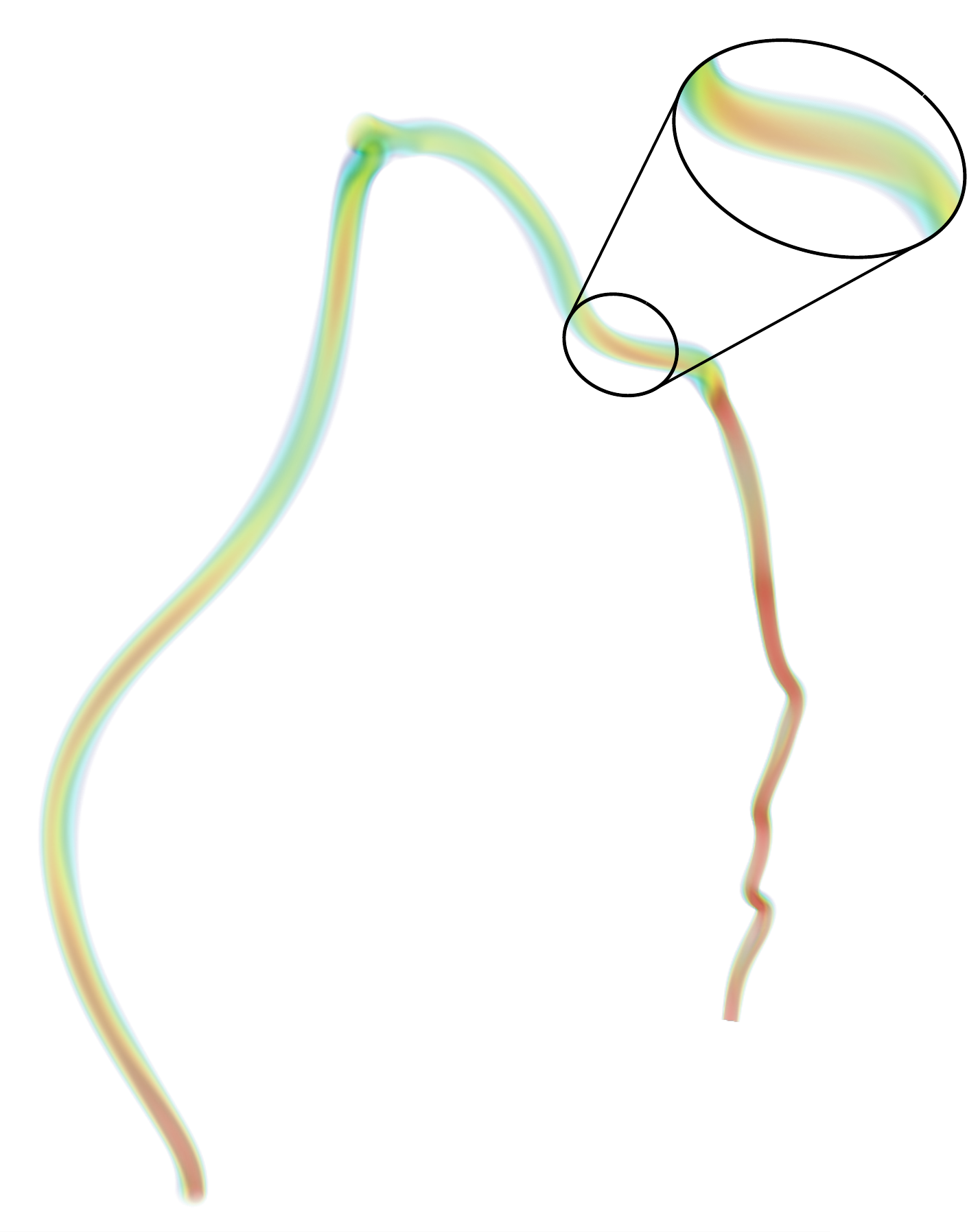} \quad\quad\quad\quad
\includegraphics[width = 6.5cm, keepaspectratio]{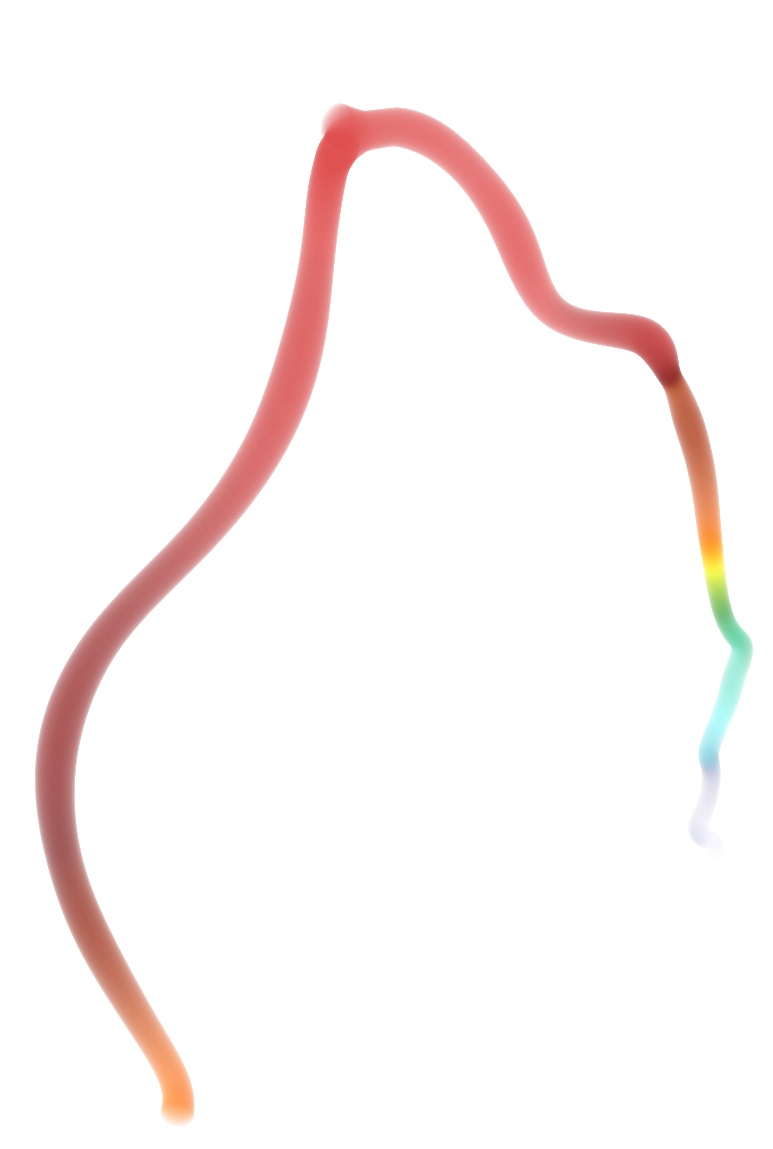}
\includegraphics[width = 6.5cm, keepaspectratio]{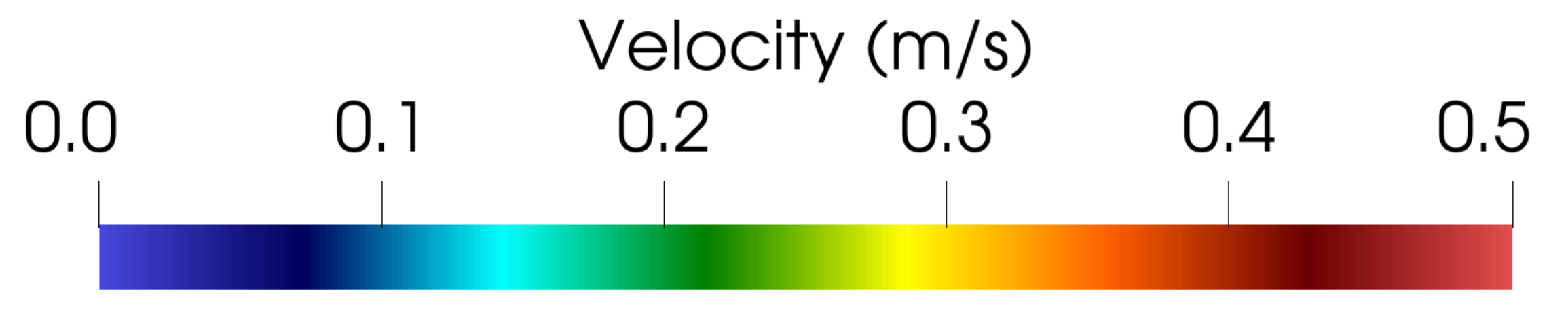} \quad\quad\quad\quad\quad\quad
\includegraphics[width = 6.5cm, keepaspectratio]{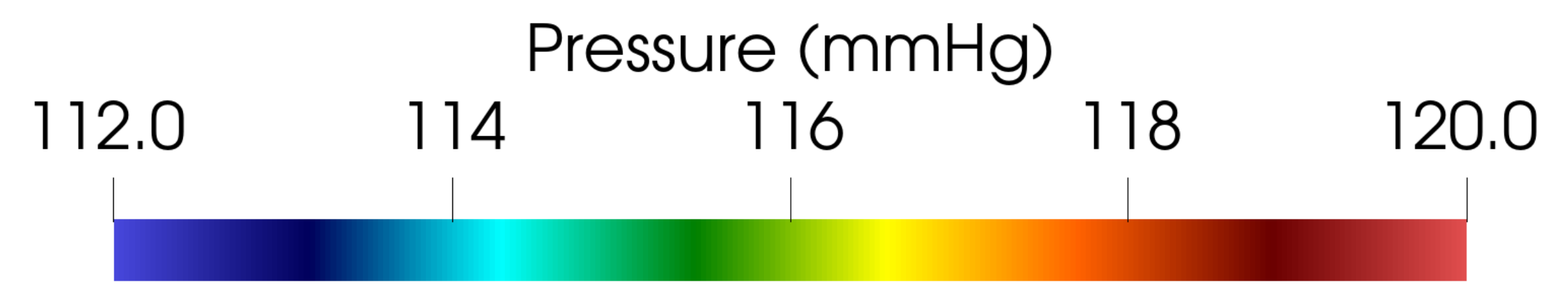}
\caption{
Left: the velocity magnitude at the systolic peak (t = 1.2s).
In the black circle: the tract of the artery considered for the mesh convergence on the WSS.
Right: Pressure distribution at the systolic peak (t = 1.2s).}
\label{Fig: result_bifurcation}
\end{figure}

\subsubsection{Mesh convergence}\label{sec: numerical_results}

\begin{figure}[h!]
\centering
\includegraphics[width = 12cm, keepaspectratio]{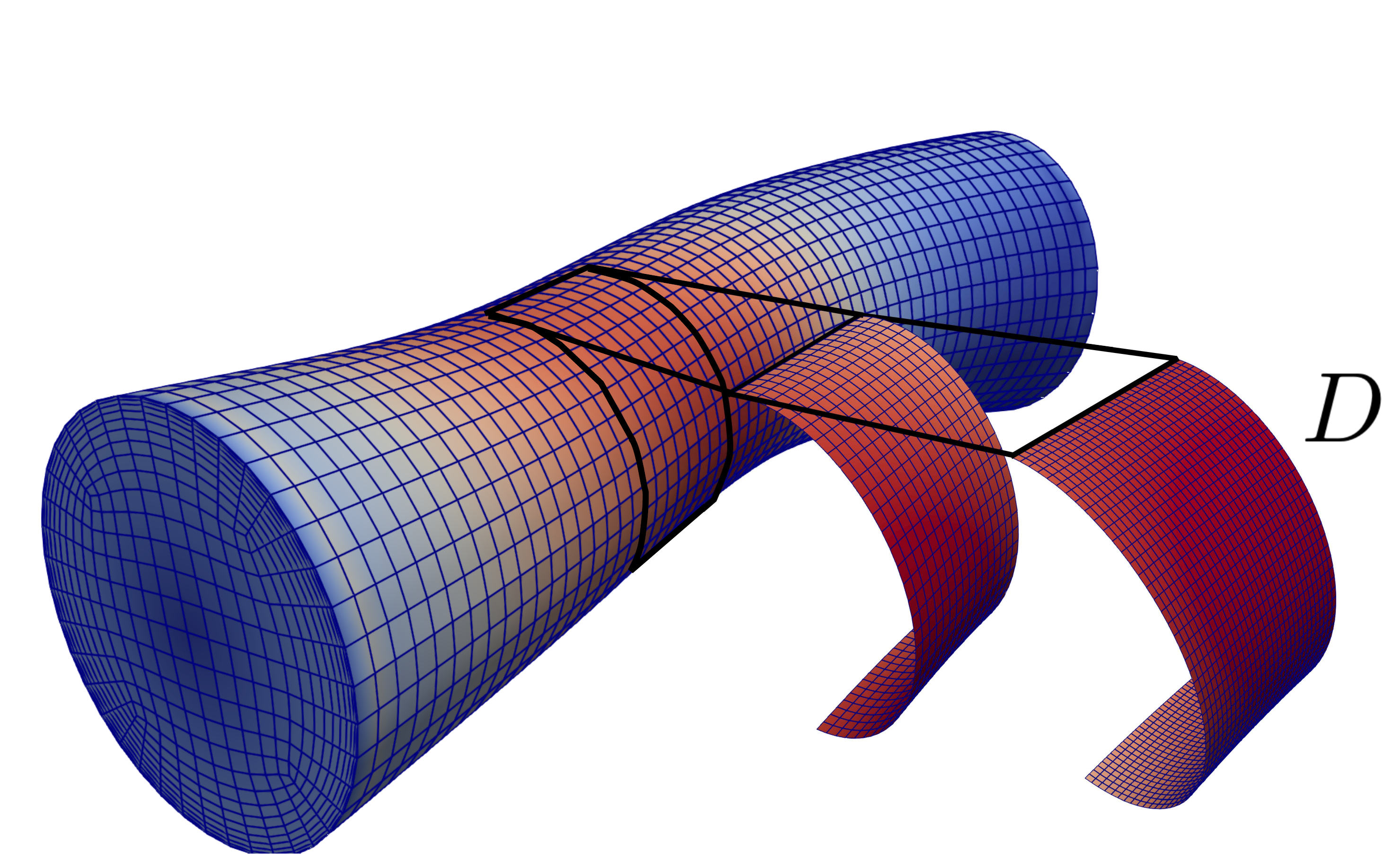}\\
\includegraphics[width = 12cm, keepaspectratio]{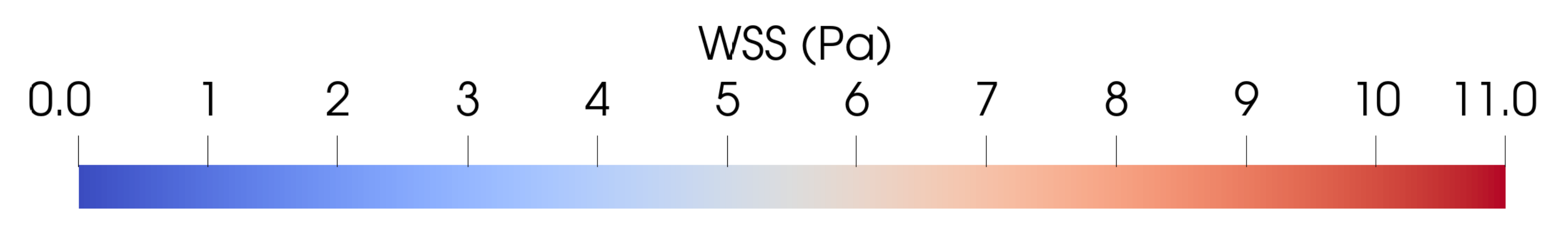}\\
\caption{
The SAWSS at the systolic peak is depicted on the 
geometry used for the test convergence. In particular, while for the coarse refinement the entire geometry is reported, for the medium and fine refinement levels only the domain $D$ is highlighted.  
}	
\label{Fig: domain_D}
\end{figure}

An important index in h{\ae}modynamics, in particular for evaluating the risk of myocardial infarction, is the Wall Shear Stress (WSS) and its derived quantities such as the Space Averaged WSS (SAWSS) and the Time Averaged WSS (TAWSS) \cite{gijsen2019expert}.
These indeces are defined as follows:
\begin{align*}
\operatorname{WSS} = \sigma\textbf{n} - (\sigma\textbf{n}\cdot\textbf{n})\textbf{n}, \quad
\operatorname{SAWSS} = \frac{1}{D}\int_{D} \operatorname{WSS}, \quad
\operatorname{TAWSS} = \frac{1}{T}\int_{0}^T \operatorname{SAWSS}, 
\end{align*}

where $\sigma$ is the Cauchy fluid stress tensor, already introduced in Sec.\ \ref{sec: fluidynamics}, and $D$ is a selected portion of the vessel segment (see the zoom from the left side in Fig.\ \ref{Fig: result_bifurcation}).
We mention that since $\sigma$ is a weak quantity, we compute the SAWSS through an integration on a small domain $D$, see Fig.\ \ref{Fig: domain_D}.
\begin{table}[t]
\caption{The relative errors with respect to the reference solution ($\operatorname{mBL}$) on TAWSS computed on the domain $D$ reported for the $6$ grids. 
$\operatorname{mBL}$ is computed as the average of the solutions using the finest grids which have a total of $2'189'100$ vertices.
Note that the TAWSS is computed until $0.1$s of the heartbeat to enable the comparison with the reference solution. In parentheses, the number of grid vertices are reported.}
\centering
{\renewcommand{\arraystretch}{2}%
\begin{tabular}{|c|c|c|c|c|c|}
\hline
 & \multicolumn{3}{|c|}{TAWSS (Pa), until $0.1$s}\\
\hline
Relative Error & Coarse (31'260) & Medium (215'410) & Fine (994'250)\\
\hline 
$\displaystyle\left|\frac{\operatorname{BL}- \operatorname{mBL}}{\operatorname{mBL}}\right| \times 100$ 
& 0.69 \% & 0.58 \% & 0.51 \%\\ 
\hline 
$\displaystyle\left|\frac{\operatorname{no BL}- \operatorname{mBL}}{\operatorname{mBL}}\right|  \times 100$ 
& 5.82 \% & 3.01 \% & 2.45 \% \\
\hline
\end{tabular}}
\label{tab: wss_convergence}
\end{table}

In Fig. \ref{Fig: result_bifurcation}, we report numerical results concerning velocity and pressure
for a patient-specific coronary artery composed by the LAD and LCx branches.
The results are in agreement with physiological values reported in the literature (see \cite{kim2010patient, guerciotti2017computational, gijsen2019expert}). This demonstrates that the framework, composed of the mesh generator and the geometric multiscale model, is able to provide meaningful results. 

To asses the quality of the boundary layer (BL), we perform a convergence test based on the WSS. This index is very sensitive to mesh refinement since it is computed using velocity gradients.
The test is carried out by selecting a tract of the artery (highlighted by a black circle in Fig.\ \ref{Fig: result_bifurcation}, left side) and then through the generation of mesh at four different levels of mesh refinement with and without BL, amounting to a total of 8 grids. Note that the grids with the same level of refinement have an equal number of vertices and the same surface. Indeed they differ only by the presence of a boundary layer. 
The boundary conditions of these simulations are taken from the main simulation of the LAD and LCx branches, see Fig.\ \ref{Fig: result_bifurcation}. 
For this, we computed the proximal flow rate and distal pressure at the corresponding sections of this reference simulation.
At the inlet, we impose the same flow rate, but we enforce it through a parabolic velocity profile. At the outlet condition, we impose a constant pressure as Neumann condition, although this is slightly
imprecise.

The finest meshes are only used to compute a reference value for the SAWSS in the region $D$ during the first 0.1 seconds. This reference value is computed as the average between the versions with and without BL, and it is defined as mBL.
In Tab.\ \ref{tab: wss_convergence}, we report the relative error of the TAWSS with respect to this reference solution until $t = 0.1$s. 
Coherently, Tab.\ \ref{tab: wss_convergence} reveals that the TAWSS of the coarse grid with BL already exhibits remarkable proximity (less than $1\%$ of difference) to our reference solution.
Note that the coarsest grid has approximately 60 times fewer cells compared to the finest one, so the presence of a BL strongly affects the computation of the WSS.


In Fig.\ \ref{Fig: domain_D}, we present the tract of the vessel of interest showing only three levels (coarse, medium and fine) of refinement and the WSS.

\begin{figure}[t]
\centering
\includegraphics[width = 11.5cm, keepaspectratio]{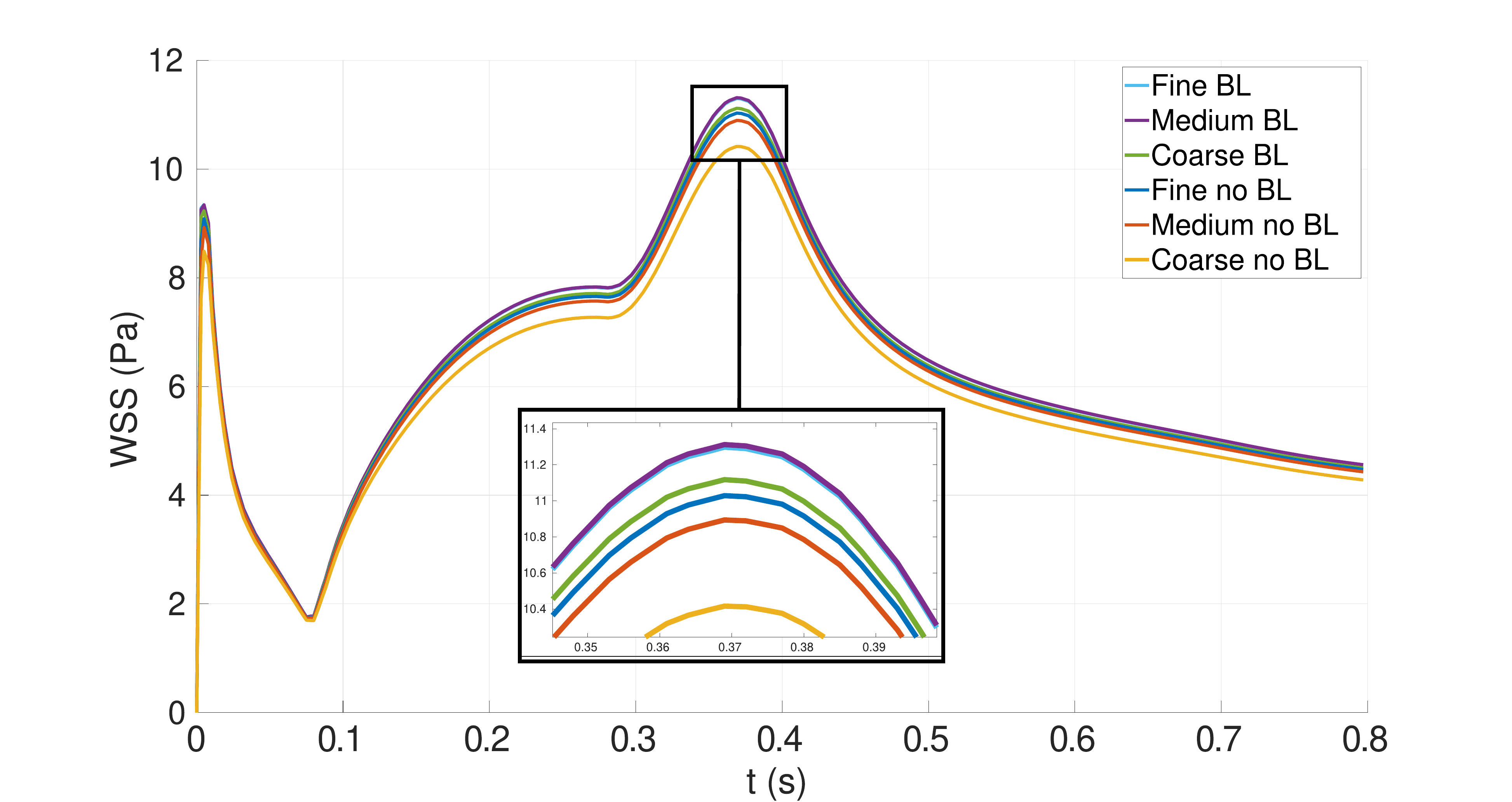}\\
\includegraphics[width = 11.5cm, keepaspectratio]{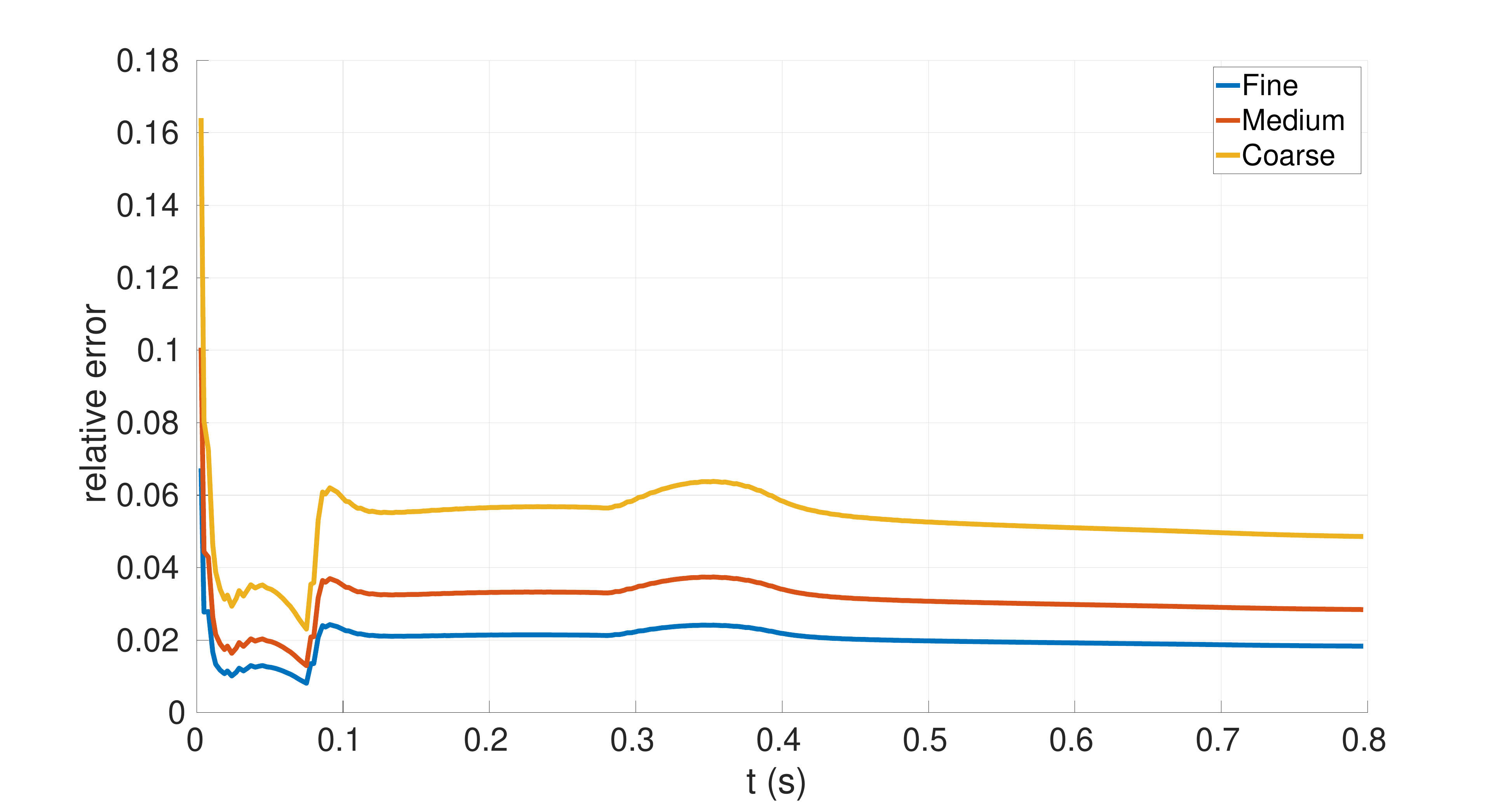}
\caption{
Top: SAWSS on the domain $D$ computed along the cardiac cycle. The coarse grid with boundary layer already shows results which are in excellent agreement with the finest grid employing 30 times fewer cells than the latter.
Bottom: Relative errors of the SAWSS along the cardiac cycle. Coherently, refining the grid leads to a smaller discrepancy between the case with boundary layer and without.  
}	
\label{Fig: wss_convergence_2}
\end{figure}

In Fig.\ \ref{Fig: wss_convergence_2}, top, the SAWSS computed on the region $D$ is plotted along the first cardiac cycle, while in the bottom row, the relative error (normalized by the WSS computed on the BL grid) between the grids is plotted. 
Fig.\ \ref{Fig: wss_convergence_2}, bottom, reveals that the discrepancy between the no BL and BL results decreases coherently when the grids become finer.
\textcolor{black}{In particular, we are plotting the relative error between the grids with and without BL, for the coarse, medium, and fine cases. Coherently, the discrepancy between the grids with BL and without BL tends to go to zero as the grid is refined.
The average error along the cardiac cycle falls from approximately 5\% to 2\%, from the coarse to the fine simulations.
}

\subsubsection{H{\ae}modynamics indices}\label{sec: FFR_vs_CFR}

A primary advantage of numerical simulations is the possibility of computing h{\ae}modynamic indices non-invasively. In clinical cardiology, the gold standard for assessing the health of the coronary artery is the Fractional Flow Reserve (FFR) \cite{pijls2012functional, johnson2014prognostic}. 
The clinicians often define it as
\begin{align*}
\operatorname{mFFR} = \frac{P_{\text{distal}}}{P_{\text{proximal}}} 
\end{align*} 
(the added 'm' stands for "measured"),
where $P_{\text{distal}}$ and $P_{\text{proximal}}$ represent the pressures downstream and upstream of the stenosis, respectively \cite{morris2013virtual}. An FFR value of less than $0.8$ indicates the presence of a significant stenosis that requires intervention.
Currently, the FFR is the clinical practice, although it is quite invasive, since it requires the insertion of a catheter to measure the average pressures upstream and downstream from the stenosis \cite{achenbach2017performing}.



The clinicians measure the FFR to have an approximation of the ratio between the flow rate in the presence ($Q_{\text{stenosis}}$) and in the absence ($Q_{\text{no stenosis}}$) of the stenosis.
Unfortunately, the latter flow rate is not measurable in vivo.
In our framework, however, by using a full spline representation of the mesh, it is possible to easily modify the geometry, as shown in Fig.\ \ref{Fig: geometry_manipulation}, perform numerical simulation on the healthy vessel, and compute $Q_{\text{no stenosis}}$. 
Therefore, we can compute the FFR in silico, more precisely, the true FFR (tFFR):
\begin{align*}
\operatorname{tFFR} = \frac{Q_{\text{stenosis}}}{Q_{\text{no stenosis}}}
\end{align*}

where, as mentioned before, $Q_{\text{stenosis}}$ and $Q_{\text{no stenosis}}$ are the flow rate in the geometry with and without stenosis, respectively.

\begin{figure}[t]
\centering
\includegraphics[width = 12cm, keepaspectratio]{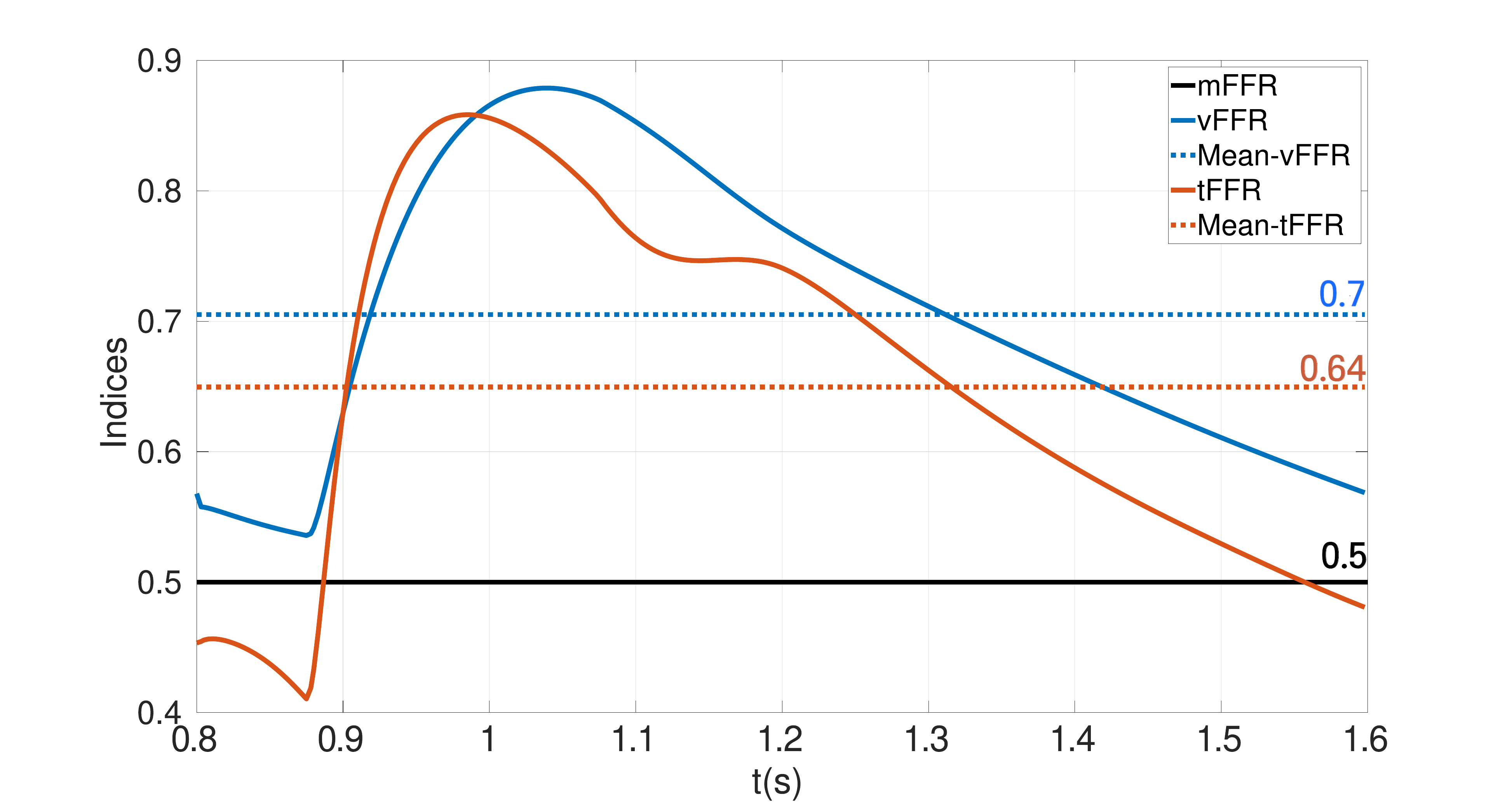}
\caption{Comparison between the mFFR (black line), the vFFR (blue line) and $\operatorname{tFFR}$ (red line).
The corresponding horizontal line represents the mean over the cardiac cycle.
}	
\label{Fig: cfr_vs_ffr}
\end{figure}

In Fig.\ \ref{Fig: cfr_vs_ffr}, the indices are compared. In particular, the mFFR computed with the numerical simulation is named vFFR \cite{morris2013virtual}. The mean values represent an average of some measurements taken during the cardiac cycle.
The results of the numerical simulations are depicted over the second cardiac cycle (i.e., when the results are at regime). 
\textcolor{black}{The discrepancies between the computed values and the measurements are due to the fact that the lumped parameters of the boundary conditions have not yet been calibrated. The boundary conditions should be patient-specific, but the focus of this section is to demonstrate that the framework can compute hemodynamic indices. Fig.\ \ref{Fig: cfr_vs_ffr} explores the relationship between these indices and, more importantly, identifies which index is more strongly correlated with myocardial infarction or other events.}

\section{Conclusion}\label{sec: conclusion}

From a methodological perspective, we proposed strategies for mesh generation of single branches and non-planar bifurcations. A key innovation of our approach is its ability to create meshes without requiring the vessel surface or any post-processing smoothing, applicable to both bifurcation and single-branch meshes. Additionally, we extended the bifurcation model to handle non-planar intersecting branches without the need for a master branch, meaning the algorithm does not rely on a primary reference branch. Leveraging spline representation for both radius and centerline data, the framework allows modifications of specific regions of interest (e.g., adding or removing stenosis or aneurysms) while maintaining consistent mesh topology across all geometric adjustments. We also discussed the generation of meshes with boundary layers.

In terms of numerical results, we compared our method against two widely adopted unstructured mesh generation tools, VMTK and Gmsh, for both single-branch and bifurcation cases. The results demonstrated the potential and superiority of our structured approach, particularly in terms of mesh quality. In a specific test case, we further compared our method to the state-of-the-art structured mesh generation techniques, achieving better results for single-branch mesh quality as evaluated by the scaled Jacobian metric.
\textcolor{black}{
We have also tested our method on a complex coronary artery tree showing its ability, thanks to the flexibility of the Hermite curves, to generate high quality mesh for the bifurcations.
} 

Moreover, we employed the generated meshes in FEM simulations, highlighting the advantages of the proposed boundary layer meshes for accurate WSS calculations with fewer cells. Various h{\ae}modynamic indices were also computed and compared with clinical measurements, underscoring the practical utility of the proposed framework.

\textcolor{black}{In a forthcoming publication, we plan to mesh a large dataset of coronary arteries and other parts of the cardiovascular system, such as the Circle of Willis and the full aortic arch.}
\textcolor{black}{Moreover, we plan to leverage the potential advantages of splines also in the simulation part, i.e.\ proposing a full spline-based framework in future work.} 

\section*{Acknowledgments}

S. Deparis and F. Marcinn\`{o} have been supported by the
Swiss National Science Foundation under project "Data-driven approximation of h{\ae}modynamics by combined reduced order modeling and deep neural networks",  n. 200021-197021.

A. Buffa and J. Hinz are partially supported by the Swiss National Science Foundation project "PDE tools for analysis-aware geometry processing in simulation science", n. 200021-215099.

We acknowledge Dr.\ T.\ Mahendiran (CHUV) and Dr.\ E.\ Andò (EPFL) for the 3D reconstructions of the arteries from the FAME 2 dataset. We thank Dr.\ Bernard De Bruyne (Cardiovascular Center OLV Aalst, Belgium) for providing permission to use the FAME 2 dataset.

\newpage
\bibliographystyle{abbrvnat}  
\bibliography{paper.bib}  

\end{document}